\documentclass[final]{siamltex} 

\usepackage{amsfonts,epsfig,cite,float,algorithmic,multirow}
\usepackage[cmex10]{amsmath} 

\newlength{\widthC}
\newlength{\widthD}
\newcounter{Lcount}

\newenvironment{adamItemize}{\begin{list}{$\bullet$}
{\setlength{\rightmargin}{0em}
\setlength{\leftmargin}{2.0em}
\setlength{\itemsep}{6pt}
\setlength{\topsep}{2pt}
\setlength{\parsep}{0pt}}}{\end{list}}

\newenvironment{adamItemize2}{\begin{list}{$\bullet$}
{\setlength{\rightmargin}{0em}
\setlength{\leftmargin}{1.0em}
\setlength{\itemsep}{4pt}
\setlength{\topsep}{2pt}
\setlength{\parsep}{0pt}}}{\end{list}}

\newcommand{\field}[1]{\mathbb{#1}}
\newcommand{\la}[1]{\mbox{$\mathbf{#1}$}}  
\newcommand{\sst}[1]{\mbox{\scriptsize{#1}}}
 
\newcommand{\Frac}[2]{{{#1}/{#2}}}  

\floatstyle{ruled}
\newfloat{algorithm}{htp}{loa}
\floatname{algorithm}{Algorithm}
\newlength{\algSkip}
\title{Simultaneously Sparse Solutions to Linear Inverse Problems with Multiple System Matrices and a Single Observation Vector\thanks{This document is a
   manuscript from July 21, 2008, unchanged except for updated references. 
   The content appears in the first author's September 2008 PhD thesis.}
}

\author{Adam C.~Zelinski\thanks{Research Laboratory of Electronics,
        Massachusetts Institute of Technology ({\tt
        zelinski@mit.edu}).}  \and Vivek K~Goyal\thanks{Research
        Laboratory of Electronics, Massachusetts Institute of
        Technology ({\tt vgoyal@mit.edu}).} \and Elfar
        Adalsteinsson\thanks{Research Laboratory of Electronics,
        Massachusetts Institute of Technology ({\tt elfar@mit.edu}).}}

\begin{document}

\maketitle

\begin{abstract}
  A linear inverse problem is proposed that requires the determination
  of multiple unknown signal vectors.  Each unknown vector passes
  through a different system matrix and the results are added to yield
  a single observation vector.  Given the matrices and lone
  observation, the objective is to find a {\em{simultaneously sparse}}
  set of unknown vectors that solves the system.  We will refer to
  this as the {\em{multiple-system single-output}} (MSSO) simultaneous
  sparsity problem.  This manuscript contrasts the MSSO problem with
  other simultaneous sparsity problems and conducts a thorough initial
  exploration of algorithms with which to solve it.

  Seven algorithms are formulated that approximately solve this
  NP-Hard problem.  Three greedy techniques are developed (matching
  pursuit, orthogonal matching pursuit, and least squares matching
  pursuit) along with four methods based on a convex relaxation
  (iteratively reweighted least squares, two forms of iterative
  shrinkage, and formulation as a second-order cone program).  While
  deriving the algorithms, we prove that {\em{seeking a single sparse
  complex-valued vector is equivalent to seeking two
  {\bf{simultaneously sparse}} real-valued vectors}}.  In other words,
  single-vector sparse approximation of a complex vector readily maps
  to the MSSO problem, increasing the applicability of MSSO
  algorithms.

  The algorithms are evaluated across three experiments: the first and
  second involve sparsity profile recovery in noiseless and noisy
  scenarios, respectively, while the third deals with magnetic
  resonance imaging radio-frequency excitation pulse design.  For
  sparsity profile recovery, the iteratively reweighted least squares
  and second-order cone programming techniques outperform the greedy
  algorithms, whereas in the magnetic resonance imaging pulse design
  experiment, the greedy least squares matching pursuit algorithm
  exhibits superior performance.  Overall, each algorithm is found to
  have particular merits and weaknesses, e.g., the iterative shrinkage
  techniques converge slowly, but because they update only a subset of
  the overall solution per iteration rather than all unknowns at once,
  they are useful in cases where attempting the latter is prohibitive
  in terms of system memory.  
\end{abstract}


\pagestyle{myheadings}
\thispagestyle{plain}
\markboth{A.~C.~ZELINSKI ET AL.}{MSSO SIMULTANEOUS SPARSE APPROXIMATION}

\section{Introduction}
\label{sec:introduction}

  In this work we propose a linear inverse problem that requires a
  {\em{simultaneously sparse}} set of vectors as the solution, i.e., a
  set of vectors where only a small number of each vector's entries
  are nonzero, and where the vectors' {\em{sparsity profiles}} (the
  locations of the nonzero entries) are equivalent (or are promoted to
  be equivalent with an increasing penalty given otherwise).

  Prior work on simultaneously sparse solutions to linear inverse
  problems involves multiple unknown vectors, a single system matrix,
  and a host of observation vectors; the $p$th observation vector
  arises by multiplying the single system matrix with the $p$th
  unknown vector \cite{Cot2005,Mal2005,Tro2006_I,Tro2006_II}.  Given
  the observation vectors and system matrix, one seeks out a
  simultaneously sparse set of unknown vectors that (approximately)
  solves the overall system.  We will refer to this as the
  {\em{single-system multiple-output}} (SSMO) simultaneous sparsity
  problem.

  Here we study a problem that is somewhat similar yet
  different from the aforementioned one.  This {\em{multiple-system
  single-output}} (MSSO) simultaneous sparsity problem still consists
  of multiple unknown vectors, but now each such vector is passed
  through a {\em{different}} system matrix and the outputs of the
  various system matrices undergo a linear combination, yielding only
  {\em{one}} observation vector.  Given the matrices and lone
  observation, one must find a simultaneously sparse set of vectors
  that (approximately) solves the system.  To date, this problem has
  been explored in a magnetic resonance imaging (MRI) radio-frequency
  (RF) excitation pulse design context
  \cite{Zel2007,Zel2008_CISS,Zel2008_TMI}, but it may also have
  applications to source localization using sensor arrays
  \cite{Joh1993, Kri1996} and signal denoising
  \cite{Don1995,Che1998,Fle2006,Ela2006}.

  Both SSMO and MSSO arise as generalizations of the single-system
  single-output (SSSO) sparse approximation problem, where there is
  one known observation vector, a known system matrix, and the
  solution one seeks is a single sparse unknown vector
  \cite{Gor1997,Che1998,Rao1998}.  Several styles of algorithms have
  been posed to solve the SSSO problem, such as forward-looking greedy
  pursuit \cite{Mal1993, Nat1995, Che1995, Dav1997, Cot1999},
  iteratively reweighted least squares (IRLS) \cite{Kar1970} (e.g.,
  FOCUSS \cite{Gor1997}), iterative shrinkage
  \cite{Don1994,Dau2004,Ela2006_TransIT,Ela2006}, and second-order
  cone programming (SOCP) \cite{Boy2004,Mal2005}.  Many of these have
  been extended to the SSMO problem
  \cite{Cot2005,Mal2005,Tro2006_I,Tro2006_II}.

  In this manuscript, we propose three forward-looking greedy
  techniques---matching pursuit (MP) \cite{Mal1993}, orthogonal
  matching pursuit (OMP) \cite{Nat1995,Che1995,Cot1999}, and least
  squares matching pursuit (LSMP) \cite{Cot1999}---and also develop
  IRLS-based, shrinkage-based, and SOCP-based algorithms to solve the
  MSSO simultaneous sparsity problem.  We then evaluate the
  performance of the algorithms across three experiments: the first
  and second involve sparsity profile recovery in noiseless and noisy
  scenarios, respectively, while the third deals with MRI RF
  excitation pulse design.

  The structure of this paper is as follows: in
  Sec.~\ref{sec:background}, we provide background information about
  ordinary (SSSO) sparse approximation and SSMO sparse approximation.
  In Sec.~\ref{sec:problem}, we formulate the MSSO problem.  Seven
  algorithms for solving the problem are then posed in
  Sec.~\ref{sec:algorithms}, while the details and results of the
  three numerical experiments appear in Sec.~\ref{sec:experiments}.
  Section~\ref{sec:discussion} highlights the strengths and weaknesses
  of the algorithms and presents ideas for future work.  Concluding
  remarks are given in Sec.~\ref{sec:conclusion}.

\section{Background}
\label{sec:background}

\subsection{Single-System Single-Output (SSSO) Sparse Approximation}

   Consider a linear system of equations $\la{d} = \la{F}\la{g}$,
   where $\la{d} \in \field{C}^{M}$, $\la{F} \in \field{C}^{M \times
   N}$, $\la{g} \in \field{C}^{N}$, and $\la{d}$ and $\la{F}$ are
   known.  It is common to use the Moore-Penrose pseudoinverse of
   $\la{F}$, denoted $\la{F}^\dagger$, to determine $\hat{\la{g}} =
   \la{F}^\dagger \la{d}$ as an (approximate) solution to the system
   of equations.  When $\la{d}$ is in the range of $\la{F}$,
   $\hat{\la{g}}$ is the solution that minimizes $\Vert \hat{\la{g}}
   \Vert_2$, the Euclidean or $\ell_2$ norm of $\hat{\la{g}}$.  When
   $\la{d}$ is not in the range of $\la{F}$, no solution exists;
   $\hat{\la{g}}$ minimizes $\Vert \hat{\la{g}} \Vert_2$ among the
   vectors that minimize $\Vert \la{d} - \la{F}\hat{\la{g}} \Vert_2$.
   
   When a sparse solution is desired, it is necessary for the analogue
   to $\hat{\la{g}}$ to have only a small fraction of its entries
   differ from zero.  We are faced with a sparse approximation problem
   of the form
   \begin{equation}\label{sparseopt_NP}
   \mathop{\mbox{min}}_{\mbox{\scriptsize{\la{g}}}}
                                \Vert \la{g} \Vert _0
                                   \text{  subject to  }
                            \Vert \la{d} - \la{F}\la{g} \Vert_2 \leq \epsilon,
   \end{equation}
   where $\Vert \cdot \Vert_0$ denotes the number of nonzero elements of a
   vector.  The subset of $\{1,\,2,\,\ldots,\,N\}$ where there are
   nonzero entries in \la{g} is called the \emph{sparsity profile}.
   For general $\la{F}$, solving (\ref{sparseopt_NP}) essentially
   requires a search over up to $2^N-1$ nonempty sparsity profiles.
   The problem is thus computationally infeasible except for very
   small systems of equations (e.g., even for $N = 30$, one may need
   to search 1,073,741,823 subsets).  Formally, the problem is
   NP-Hard \cite{Dav1994,Nat1995}.

   For problems where (\ref{sparseopt_NP}) is intractable, a large
   body of work supports a greedy search over the columns of \la{F} to
   seek out a small subset of columns that, when weighted and linearly
   combined, yields a result that is close to \la{d} in the $\ell_2$
   sense, along with a sparse \la{g} \cite{Mal1993, Nat1995, Che1995,
   Dav1997, Cot1999}.

   A second body of research supports the {\em{relaxation}} of 
   (\ref{sparseopt_NP}) to find a sparse \la{g} \cite{Che1998}:
   \begin{equation}\label{sparseopt_relaxed}
   \mathop{\mbox{min}}_{\mbox{\scriptsize{\la{g}}}}
                      \Vert \la{g} \Vert_1
                         \text{ s.t. }
                      \Vert \la{d} - \la{F}\la{g} \Vert_2 \leq \epsilon.
   \end{equation}
   This is a convex optimization and thus may be solved efficiently
   \cite{Boy2004}.  The solution of (\ref{sparseopt_relaxed}) does not
   always match the solution of (\ref{sparseopt_NP})---if it did, the
   intractability of (\ref{sparseopt_NP}) would be contradicted---but
   certain conditions on $\la{F}$ guarantee a proximity of their
   solutions \cite{Don2006b, Don2006c, Tro2006_convexIT}.
   Note that (\ref{sparseopt_relaxed}) applies an $\ell_1$ norm to
   $\la{g}$, but an $\ell_p$ norm (where $p < 1$) may also be used to
   promote sparsity \cite{Gor1997, Che1998}; this leads to a
   non-convex problem and will not considered in this paper.

   The optimization
   \begin{equation}\label{sparseopt_noisy}
   \mathop{\mbox{min}}_{\mbox{\scriptsize{\la{g}}}} \mbox{  }
                  \left \{
                        \tfrac{1}{2} \Vert \la{d} - \la{F}\la{g} \Vert _2^2
                      +     \lambda  \Vert \la{g} \Vert _1
                  \right \}
   \end{equation}
   has the same set of solutions as (\ref{sparseopt_relaxed}).  The
   first term of (\ref{sparseopt_noisy}) keeps residual error down,
   whereas the second promotes sparsity of \la{g} \cite{Tib1996,
   Che1998}.  As the control parameter, $\lambda$, is increased from
   zero to infinity, the algorithm yields sparser solutions and the
   residual error increases; sparsity is traded off with residual
   error.  In this paper we shall hereafter use formulation
   (\ref{sparseopt_noisy}) and its analogues rather than
   (\ref{sparseopt_relaxed}).

   It is important to understand that a problem of the form
   (\ref{sparseopt_noisy}) may arise as a proxy for (\ref{sparseopt_NP})
   or as the inherent problem of interest.  For example,
   in a Bayesian estimation setting, (\ref{sparseopt_noisy})
   yields the maximum a posteriori probability estimate of
   $\la{g}$ from $\la{d}$ when the observation model involves $\la{F}$
   and Gaussian noise and the prior on $\la{g}$ is Laplacian.
   Similar statements can be made about the relaxations of the
   \emph{simultaneous} sparse approximation problems posed in the
   following sections.

\subsection{Single-System Multiple-Output (SSMO) Simultaneous Sparse Approximation}

   In SSMO, we have $P$ observation vectors (``snapshots''), all of
   which arise from the same system matrix:
   \begin{equation}\label{mult_obs} 
        \la{d}_p = \la{F} \la{g}_p, \,\,\mbox{for } p = 1, \ldots, P,
   \end{equation}
   where $\la{d}_p \in \field{C}^{M}$ is known for $p = 1, \ldots, P$
   along with $\la{F} \in \field{C}^{M \times N}$.  In this scenario,
   we want to constrain the number of positions at which \emph{any} of
   the $\la{g}_p$s are nonzero.  Thus we seek approximate solutions in
   which the $\la{g}_p$s are not only sparse, but the union of their
   sparsity patterns is small; that is, a {\em{simultaneously sparse}}
   set of vectors is desired \cite{Mal2005, Tro2006_II}.
   Unfortunately, optimal approximation with a simultaneous sparsity
   constraint is even harder than (\ref{sparseopt_NP}).

   Extending single-vector sparse approximation greedy techniques is
   one way to find an approximate solution \cite{Cot2005,Tro2006_I}.
   Another approach is to extend the relaxation
   (\ref{sparseopt_noisy}) as follows:
   \begin{equation}\label{mmv1}
   \mathop{\mbox{min}}_{\mbox{\scriptsize{\la{G}}}} \mbox{  }
      \left \{
        \tfrac{1}{2}  \left \Vert \la{D} - \la{F} \la{G}
                      \right \Vert_{\mbox{\scriptsize{F}}}^2
       +     \lambda  \left \Vert \la{G} \right \Vert_{\mbox{\scriptsize{S}}}
      \right \},
   \end{equation}
   where $\la{D} = [ \la{d}_1, \ldots, \la{d}_P ] \in \field{C}^{M
   \times P}$, $\la{G} = [\la{g}_1, \ldots, \la{g}_P] \in \field{C}^{N
   \times P}$, $\Vert \cdot \Vert_{\mbox{\scriptsize{F}}}$ is the
   Frobenius norm, and
   \begin{equation}\label{S_norm}
     \Vert \la{G} \Vert_{\mbox{\scriptsize{S}}} 
        = \sum_{n=1}^{N}{ \sqrt{ \sum_{p=1}^{P}{ \vert \la{G}\left(n,p\right) \vert^2 }} } 
        = \sum_{n=1}^{N}{ \sqrt{ \sum_{p=1}^{P}{ \vert \la{g}_p[n] \vert^2 }} },
   \end{equation}
   i.e., $\Vert \la{G} \Vert_{\sst{S}}$ is the $\ell_1$ norm of the
   $\ell_2$ norms of the rows of the \la{G} matrix.\footnote{Although
   here we have applied an $\ell_1$ norm to the $\ell_2$ row energies
   of $\la{G}$, an $\ell_p$ norm (where $p < 1$) could be used in
   place of the $\ell_1$ norm if one is willing to deal with a
   non-convex objective function.  Further, an $\ell_q$ norm (where $q
   > 2$) rather than an $\ell_2$ norm could be applied to each row of
   $\la{G}$ because the purpose of the row operation is to collapse
   the elements of the row into a scalar value without introducing a
   sparsifying effect.} This latter operator is a {\em{simultaneous
   sparsity norm}}: it penalizes the program (produces large values)
   when the columns of \la{G} have dissimilar sparsity profiles
   \cite{Mal2005}.  Fixing $\lambda$ to a sufficiently large value and
   solving this optimization yields simultaneously sparse $\la{g}_p$s.
   For $P = 1$, (\ref{mmv1}) collapses to the base case of
   (\ref{sparseopt_noisy}).  Given the convex objective function in
   (\ref{mmv1}), one may then attempt to find a solution that
   minimizes the objective using an IRLS-based \cite{Cot2005}, or
   SOCP-based \cite{Mal2005} approach.  A formal analysis of the
   minimization of the convex objective may be found in
   \cite{Tro2006_II}.

\section{Multiple-System Single-Output (MSSO) Simultaneous Sparse Approximation}
\label{sec:problem}

   We outline the MSSO problem in a style analogous to that of SSMO
   systems in (\ref{mult_obs}, \ref{mmv1}) and then pose a second
   formulation that is useful for deriving several algorithms in
   Sec.~\ref{sec:algorithms}.

\subsection{MSSO Problem Formulation}

   Consider the following system:
   \begin{equation}\label{mmv2_eq1}
   \begin{split}
   \la{d} & = \la{F}_1 \la{g}_1 + \cdots + \la{F}_P \la{g}_P  \\
          & = \left [ \la{F}_1 \cdots \la{F}_P \right ]
                                  \left [ \begin{array}{c}
                                                          \la{g}_1 \\
                                                           \vdots  \\
                                                          \la{g}_P \\
                                          \end{array}
                                  \right ]
          = \la{F}_{\mbox{\scriptsize{tot}}}
            \la{g}_{\mbox{\scriptsize{tot}}},
   \end{split}
   \end{equation}
   where $\la{d} \in \field{C}^{M}$ and the $\la{F}_p \in \field{C}^{M
   \times N}$ are known.  Unlike the SSMO problem, there is now only one
   observation and $P$ different system matrices.  Here we again 
   desire an approximate solution where the $\la{g}_p$s are simultaneously sparse;
   formally,
   \begin{equation}\label{msso_nphard}
         \mathop{\mbox{min}}_{\sst{$\la{g}_1, \ldots, \la{g}_P$}} 
            \mbox{ } \Vert \la{d} - 
         \la{F}_{\mbox{\scriptsize{tot}}} \la{g}_{\mbox{\scriptsize{tot}}}
         \Vert_2 \,\,
            \mbox{ s.t. the simultaneous $K$-sparsity of the $\la{g}_p$s}.
   \end{equation}
   This is, of course, harder than the SSSO problem and thus NP-Hard.
   To keep
   the problem as general as possible, there are no constraints
   on the values of $M$, $N$, or $P$, i.e., there is no explicit requirement
   that the system be overdetermined or underdetermined.  Further,
   we have used complex-valued rather than real-valued variables.

   In the first half of Sec.~\ref{sec:algorithms}, we will pose three
   approaches that attempt to solve the MSSO problem
   (\ref{msso_nphard}) in a greedy fashion.  Another approach to solve
   the problem is to apply a relaxation similar to
   (\ref{sparseopt_noisy}, \ref{mmv1}):
   \begin{equation}\label{mmv2_eq2}
   \mathop{\mbox{min}}_{\mbox{\scriptsize{\la{G}}}} \mbox{  }
      \left \{
        \tfrac{1}{2}  \left \Vert \la{d} -
                      \la{F}_{\mbox{\scriptsize{tot}}}
                      \la{g}_{\mbox{\scriptsize{tot}}}
                      \right \Vert_2^2
       +     \lambda  \left \Vert \la{G} \right \Vert_{\mbox{\scriptsize{S}}}
      \right \},
   \end{equation}
   where \la{G} and $\Vert \la{G} \Vert_{\mbox{\scriptsize{S}}}$ are
   the same as in (\ref{mmv1}) and (\ref{S_norm}), respectively. In
   the second half of Sec.~\ref{sec:algorithms}, we will outline four
   algorithms for solving this relaxed problem.  By setting $P = 1$,
   (\ref{mmv2_eq2}) collapses to the base case of
   (\ref{sparseopt_noisy}).

\subsection{Alternate Formulation of the MSSO Problem}

   In several upcoming derivations, it will be useful to view the system
   from a different standpoint.  To begin, we construct several new variables
   that are simply permutations of the $\la{F}_p$s and $\la{g}_p$s.
   First we define $N$ new matrices:
   \begin{equation}\label{C_n}
     \la{C}_n = [ \la{f}_{1,n}, \ldots, \la{f}_{P,n} ]
              \in \field{C}^{M \times P},
              \mbox{ for } n = 1, \ldots, N,
   \end{equation}
   where $\la{f}_{p,n}$ is the $n$th column of $\la{F}_p$.
   Next we construct $N$ new vectors:
   \begin{equation}\label{h_n}
     \la{h}_n = [ \la{g}_1[n], \ldots, \la{g}_P[n] ]^{\sst{T}} 
              \in \field{C}^{P}, \mbox{ for } n = 1, \ldots, N,
   \end{equation}
   where $\la{g}_p[n]$ is the $n$th element of $\la{g}_p$ and
   $^{\sst{T}}$ is the transpose operation.
   Given the $\la{C}_n$s and $\la{h}_n$s, we now pose the
   following system:
   \begin{equation}\label{mmv2_eq3}
   \begin{split}
   \la{d} & = \la{C}_1 \la{h}_1 + \cdots + \la{C}_N \la{h}_N  \\
          & = \left [ \la{C}_1 \cdots \la{C}_N \right ]
                                  \left [ \begin{array}{c}
                                                          \la{h}_1 \\
                                                           \vdots  \\
                                                          \la{h}_N \\
                                          \end{array}
                                  \right ]
          = \la{C}_{\mbox{\scriptsize{tot}}}
            \la{h}_{\mbox{\scriptsize{tot}}}.
   \end{split}
   \end{equation}
   Due to (\ref{C_n}, \ref{h_n}), the system in (\ref{mmv2_eq3}) is 
   {\em{equivalent}} to the one in (\ref{mmv2_eq1}).  The relationship
   between the $\la{g}_p$s and $\la{h}_n$s implies that if we desire
   to find a set of simultaneously sparse $\la{g}_p$s to solve
   (\ref{mmv2_eq1}), we should seek out a set of $\la{h}_n$s where
   many of the $\la{h}_n$s equal an all-zeros vector, \la{0}, but a
   few $\la{h}_n$s are high in $\ell_2$ energy (typically with all
   elements being nonzero).  This claim is apparent if we consider the
   fact that $\la{H} = [\la{h}_1, \ldots, \la{h}_N]$ is equal to the
   transpose of $\la{G}$, and that the $\la{g}_p$s are only 
   simultaneously sparse when $\Vert \la{G}
   \Vert_{\mbox{\scriptsize{S}}}$ is sufficiently small.

   Continuing with this alternate formulation, and given our desire
   to find a solution where most of the $\la{h}_n$s are all-zero
   vectors and a few are dense, we relax the problem as follows:
   \begin{equation}\label{mmv2_eq4}
   \mathop{\mbox{min}}_{\mbox{\scriptsize{$\la{h}_{\sst{tot}}$}}}
   \mbox{  }
      \left \{
        \tfrac{1}{2}  \left \Vert \la{d} -
                      \la{C}_{\mbox{\scriptsize{tot}}}
                      \la{h}_{\mbox{\scriptsize{tot}}}
                      \right \Vert_2^2
       +     \lambda  \sum_{n=1}^{N} \Vert \la{h}_n \Vert_2
      \right \}.
   \end{equation}
   Fixing $\lambda$ to a sufficiently large value and solving this
   optimization yields many low-energy $\la{h}_n$s (each close to
   \la{0}), along with several dense high-energy $\la{h}_n$s.
   Further, because $\sum_{n=1}^{N} \Vert \la{h}_n \Vert_2$ is
   {\em{equivalent}} to $\Vert \la{G} \Vert_{\mbox{\scriptsize{S}}}$,
   this means (\ref{mmv2_eq4}) is {\em{equivalent}} to
   (\ref{mmv2_eq2}), and thus just like (\ref{mmv2_eq2}), the approach
   in (\ref{mmv2_eq4}) finds a set of simultaneously sparse
   $\la{g}_p$s.

\subsection{Differences between the SSMO and MSSO Problems}

   In the SSMO problem, we see from (\ref{mult_obs}) that there are
   many different \la{d}s and a single \la{F}.  The ratio of unknowns
   to knowns always equals $N/M$ regardless of the number of
   observations, $P$.  A large $P$ when solving SSMO is actually
   beneficial because the simultaneous sparsity of the underlying
   $\la{g}_p$s becomes more exploitable; empirical evidence of
   improved sparsity profile recovery with increasing $P$ may be found
   in both \cite{Cot2005} and \cite{Mal2005}.

   In contrast, we see from (\ref{mmv2_eq1}) that in the MSSO problem
   there is a single \la{d} and many different \la{F}s.  Here the
   ratio of unknowns to knowns is no longer constant with respect to
   $P$; rather it is equal to $\Frac{PN}{M}$.  We will show in
   Sec.~\ref{sec:experiments} that as $P$ increases, the underlying
   simultaneous sparsity of the $\la{g}_p$s is not enough to combat
   the increasing number of unknowns, and that for large $P$,
   correctly identifying the sparsity profile of the underlying
   unknown $\la{g}_p$s is a daunting task.
   
\section{Proposed Algorithms}
\label{sec:algorithms}

We
now derive algorithms to (approximately) solve
the MSSO problem defined in Sec.~\ref{sec:problem}.

\subsection{Matching Pursuit (MP)}

    To begin, we extend the single-vector case of matching pursuit
    \cite{Mal1993} to an MSSO context.  The classic MP technique first
    finds the column of the system matrix that best matches with the
    observed vector and then removes from the observation vector the
    projection of this chosen column.  It proceeds to select a second
    column of the system matrix that best matches with the
    {\em{residual observation}}, and continues doing so until either
    $K$ columns have been chosen (as specified by the user) or the
    residual observation ends up as a vector of all zeros.  This
    algorithm is fast and computationally-efficient because the
    best-matching column vector during each iteration is determined
    simply by calculating the inner product of each column vector with
    the residual observation and ranking the resulting inner product
    magnitudes.
    
    Now let us consider the MSSO system as posed in (\ref{mmv2_eq3}).  In
    the first iteration of standard MP, we seek out the single column
    of the system matrix that can best represent \la{d}.  But in
    the MSSO context, we need to seek out which of the $N$ $\la{C}_n$
    matrices can be best used to represent \la{d} when the columns of
    $\la{C}_n$ undergo an arbitrarily-weighted linear combination.
    The key difference here is that on an iteration-by-iteration
    basis, we are no longer deciding which column vector best
    represents the observation, but which {\em{matrix}} does so.  The
    intuition behind this approach is that ideal solutions should
    consist of only a few dense $\la{h}_n$s and many all-zeros
    $\la{h}_n$s.  For the $k$th iteration of the algorithm, we 
    need to select the proper index $n \in \{1, \ldots, N\}$ by solving:
    \begin{equation}\label{MP_1}
       q_k =
       \mathop{\mbox{argmin}}_{\mbox{\scriptsize{$n$}}} 
         \mbox{  }
       \mathop{\mbox{min}}_{\mbox{\scriptsize{$\la{h}_n$}}}   
         \mbox{  }
       \Vert \la{r}_{k-1} - \la{C}_n \la{h}_n \Vert_2^2,
    \end{equation}
    where $q_k$ is the index that will be selected and $\la{r}_{k-1}$
    is the current residual observation.  For fixed $n$, the solution
    to the inner minimization is obtained via the pseudoinverse,
    $\la{h}_n^{\sst{opt}} = \la{C}_n^{\dagger} \la{r}_{k-1}$, yielding
    \begin{equation}\label{MP_2}
       q_k =
       \mathop{\mbox{argmin}}_{\mbox{\scriptsize{$n$}}} 
         \mbox{  }
       \mathop{\mbox{min}}_{\mbox{\scriptsize{$\la{h}_n$}}}   
         \mbox{  }
       \Vert \la{r}_{k-1} - \la{C}_n ( \la{C}_n^{\dagger} \la{r}_{k-1} ) \Vert_2^2
	   =
       \mathop{\mbox{argmax}}_{\mbox{\scriptsize{$n$}}}
         \mbox{  }
	\la{r}_{k-1}^{\sst{H}} \la{C}_n \la{C}_n^{\dagger} \la{r}_{k-1},
    \end{equation}
    where $^{\sst{H}}$ is the Hermitian transpose.  From (\ref{MP_2})
    we see that, analogously to standard MP, choosing the best index
    for iteration $k$ involves computing and ranking a series of
    inner-product-like quadratic terms.

   {\em{Algorithm \ref{alg:MP}}} outlines the entire procedure.
   After $K$ iterations, one obtains $I_K \subset \{1, \ldots, N \}$
   (of cardinality $T \leq K$), a set indicating the chosen $\la{C}_n$
   matrices.  The weights to apply to each chosen matrix (i.e., the
   corresponding $\la{h}_n$s) are obtained via a finalization step;
   they are the best weightings in the $\ell_2$ residual error sense
   with which to linearly combine the columns of the chosen $\la{C}_n$
   matrices to best match the observation \la{d}.  Since $T$ total
   matrices end up being chosen by the process, there is no penalty in
   retuning the $T$ associated $\la{h}_n$ vectors because they are
   allowed to be dense.  The $N - T$ other $\la{C}_n$s (and
   corresponding $\la{h}_n$s) are not used.\footnote{From the
   perspective of $\la{F}_p$s and $\la{g}_p$s in (\ref{mmv2_eq1}),
   {\em{Algorithm \ref{alg:MP}}} determines weights to place along
   only $T$ rows of $\la{G}$ (leaving the other $N - T$ rows zeroed
   out) that still yields a good approximation of \la{d} in the
   $\ell_2$ residual error sense.  It is seeking out the best rows of
   $\la{G}$ which, when densely filled, yield a sound approximation of
   $\la{d}$.}

   One property of note is that if $M \leq P$, {\em{Algorithm
   \ref{alg:MP}}} stops after one iteration.  This is because
   $\la{C}_n \la{C}_n^{\dagger}$ in this case is simply an $M \times
   M$ identity matrix for all $n \in \{1, \ldots, N \}$, and thus any
   one of the $\la{C}_n$s is enough to represent \la{d} exactly when
   its columns are properly weighted and linearly combined.

\begin{algorithm}
    \caption{--- MSSO Matching Pursuit}
    \label{alg:MP}
  {\small

    {\bf{Task}:} greedily choose up to $K$ of the $\la{C}_n$s
                to best represent \la{d} via $\la{C}_1 \la{h}_1 + \cdots + \la{C}_N \la{h}_N$.\\[\algSkip]
    {\bf{Data and Parameters}:} $\la{d}$ and $\la{C}_n, n \in \{1, \ldots, N \}$ are given. 
                           $K$ iterations.\\[\algSkip]
    {\bf{Precompute}:} $\la{Q}_n = \la{C}_n \la{C}_n^{\dagger}$, for $n \in \{1, \ldots, N\}$.\\[\algSkip]
    {\bf{Initialize}:} Set $k=0$, $\la{r}_0 = \la{d}$, $I_0 = \emptyset$, $\la{S}_0 = [\mbox{  }]$.\\[\algSkip]
    {\bf{Iterate}:} Set $k=1$ and apply:
    \begin{adamItemize}

      \item{$q_k = \mathop{\mbox{argmax}}_{\,\mbox{\scriptsize{$n$}}}
               \mbox{  }
        	\la{r}_{k-1}^{\sst{H}} \la{Q}_n \la{r}_{k-1}$.}

      \item{
      \begin{algorithmic}
	\IF {$q_k \notin I_{k-1}$}
	        \STATE $I_k = I_{k-1} \cup \{ q_k \}$
		\STATE $\la{S}_k = [ \la{S}_{k-1}, \la{C}_{q_k} ]$
	\ELSE
		\STATE $I_k = I_{k-1}$
		\STATE $\la{S}_k = \la{S}_{k-1}$
	\ENDIF 
      \end{algorithmic}
      }

      \item{$\la{r}_k = \la{r}_{k-1} - \la{Q}_{q_k} \la{r}_{k-1}$.}

      \item{$k = k + 1$.  Terminate loop if $k > K$ or $\la{r}_k = \la{0}$.  
            $I_K$ ends with $T \leq K$ elements.}

    \end{adamItemize}

    {\bf{Compute Weights}:} $\la{x} = \la{S}_K^{\dagger} \la{d}$,
    unstack \la{x} into $\la{h}_{q_1}, \ldots, \la{h}_{q_T}$; set
    remaining $\la{h}_n$s to \la{0}.
  }
\end{algorithm}

\subsection{Orthogonal Matching Pursuit (OMP)}

    In single-vector MP, the residual $\la{r}_k$ always ends up
    orthogonal to the $k$th column of the system matrix, but as the
    algorithm continues iterating, there is no guarantee the residual
    remains orthogonal to column $k$ or is minimized in the
    least-squares sense with respect to the entire set of $k$ chosen
    column vectors (indexed by $q_1, \ldots, q_k$).  Furthermore, $K$
    iterations of single-vector MP do not guarantee $K$ different
    columns will be selected.  Single-vector OMP is an extension to MP
    that attempts to mitigate these problems by improving the
    calculation of the residual vector.  During the $k$th iteration of
    single-vector OMP, column $q_k$ is chosen exactly as in MP (by
    ranking the inner products of the residual vector $\la{r}_{k-1}$
    with the various column vectors), but the residual vector is
    updated by accounting for {\em{all}} columns chosen up through
    iteration $k$ rather than simply the last one \cite{Nat1995,
    Cot1999}.

    To extend OMP to the MSSO problem, we choose matrix $q_k$ during
    iteration $k$ as in MSSO MP and then in the spirit of
    single-vector OMP compute the new residual as follows:
    \begin{equation}\label{omp1} 
      \la{r}_k = \la{d} - \la{S}_{k} ( \la{S}_{k}^{\dagger} \la{d} ),
    \end{equation} 
    where $\la{S}_k = \left [ \la{C}_{q_1}, \ldots, \la{C}_{q_k}
    \right ]$ and $\la{S}_{k}^{\dagger} \la{d}$ is the best choice of
    \la{x} that minimizes the residual error $\Vert \la{d} - \la{S}_k
    \la{x} \Vert_2$.  That is, to update the residual we now employ
    all chosen matrices, weighting and combining them to best
    represent \la{d} in the least-squares sense, yielding an
    $\la{r}_k$ that is orthogonal to the columns of $\la{S}_k$ (and
    thus orthogonal to $\la{C}_{q_1}, \ldots, \la{C}_{q_k}$), which
    has the benefit of ensuring that OMP will select a new $\la{C}_n$
    matrix at each step.

    {\em{Algorithm \ref{alg:OMP}}} describes the OMP algorithm; the
    complexity here is moderately greater than that of MP because the
    pseudoinversion of an $M \times Pk$ matrix is required during each
    iteration $k$.

\begin{algorithm}
    \caption{--- MSSO Orthogonal Matching Pursuit}
    \label{alg:OMP}
  {\small
    {\bf{Task}:} greedily choose up to $K$ of the $\la{C}_n$s
                to best represent \la{d} via $\la{C}_1 \la{h}_1 + \cdots + \la{C}_N \la{h}_N$.\\[\algSkip]
    {\bf{Data and Parameters}:} $\la{d}$ and $\la{C}_n, n \in \{1, \ldots, N \}$ are given. 
                           $K$ iterations.\\[\algSkip]
    {\bf{Precompute}:} $\la{Q}_n = \la{C}_n \la{C}_n^{\dagger}$, for $n \in \{1, \ldots, N\}$.\\[\algSkip]
    {\bf{Initialize}:} Set $k=0$, $\la{r}_0 = \la{d}$, $I_0 = \emptyset$, $\la{S}_0 = [\mbox{  }]$.\\[\algSkip]
    {\bf{Iterate}:} Set $k=1$ and apply:

    \begin{adamItemize}

      \item{$q_k = \mathop{\mbox{argmax}}_{\,\mbox{\scriptsize{$n \notin I_{k-1}$}}}
               \mbox{  }
        	\la{r}_{k-1}^{\sst{H}} \la{Q}_n \la{r}_{k-1}$.}

      \item{$I_k = I_{k-1} \cup \{ q_k \}$}

      \item{$\la{S}_k = [ \la{S}_{k-1}, \la{C}_{q_k} ] $}

      \item{$\la{r}_k = \la{d} - \la{S}_k \la{S}_k^{\dagger} \la{d}$.}

      \item{$k = k + 1$.  Terminate loop if $k > K$ or $\la{r}_k = \la{0}$.  
            $I_K$ ends with $T \leq K$ elements.}

    \end{adamItemize}

    {\bf{Compute Weights}:} $\la{x} = \la{S}_K^{\dagger} \la{d}$,
    unstack \la{x} into $\la{h}_{q_1}, \ldots, \la{h}_{q_T}$; set
    remaining $\la{h}_n$s to \la{0}.
  }
\end{algorithm}

\subsection{Least Squares Matching Pursuit (LSMP)}

    Beyond OMP there exists a greedy algorithm with an even greater
    computational complexity known as LSMP\@.  In single-vector LSMP,
    one solves $N-k+1$ least squares minimizations during iteration
    $k$ in order to determine which column of the system matrix may be
    used to best represent \la{d} \cite{Cot1999}.

    Thus to extend LSMP to MSSO systems, we must ensure that during
    iteration $k$ we account for the $k-1$ previously chosen $\la{C}_n$
    matrices when choosing the $k$th one to best construct an
    approximation to \la{d}.  Specifically, index $q_k$ is
    selected as follows:
    \begin{equation}\label{LSMP_1}
       q_k =
       \mathop{\mbox{argmin}}_{\mbox{\scriptsize{$n \in \{1, \ldots, N\}, n \notin I_{k-1}$}}} 
         \mbox{  }
       \mathop{\mbox{min}}_{\mbox{\scriptsize{$\la{x}$}}}   
         \mbox{  }
       \Vert \la{S}_{k}^{(n)} \la{x} - \la{d} \Vert_2^2,
    \end{equation}
    where $I_{k-1}$ is the set of indices chosen up through iteration
    $k-1$, $\la{S}_k^{(n)} = [\la{S}_{k-1}, \la{C}_n]$, $\la{S}_{k-1}
    = [ \la{C}_{q_1}, \ldots, \la{C}_{q_{k-1}} ]$, and $\la{x} \in
    \field{C}^{Pk}$.  For fixed $n$, the solution of the inner
    iteration is $\la{x}^{\sst{opt}} = (\la{S}_k^{(n)})^{\dagger}
    \la{d}$; it is this step that ensures the residual observation
    error will be minimized by using {\em{all}} chosen matrices.
    Substituting $\la{x}^{\sst{opt}}$ into (\ref{LSMP_1}) and
    simplifying the expression yields 
    \begin{equation}\label{LSMP_2}
       q_k = \mathop{\mbox{argmax}}_{\,\mbox{\scriptsize{$n \notin I_{k-1}$}}}
         \mbox{ } \la{d}^{\sst{H}} \la{Q}_k^{(n)} \la{d},
    \end{equation}
    where $\la{Q}_k^{(n)} = (\la{S}_k^{(n)})(\la{S}_k^{(n)})^{\dagger}$.
    
    {\em{Algorithm \ref{alg:LSMP}}} describes the LSMP method.  The
    complexity here is much greater than that of OMP because $N-k+1$
    pseudoinversions of an $M \times Pk$ matrix are required during each
    iteration $k$.  Furthermore, the dependence of $\la{Q}_k^{(n)}$ on
    both $n$ and $k$ makes precomputing all such matrices infeasible
    in most cases.  One way to decrease computation and runtime might
    be to extend the projection-based recursive updating scheme of
    \cite{Cot1999} to MSSO LSMP\@.

\begin{algorithm}
    \caption{--- MSSO Least Squares Matching Pursuit}
    \label{alg:LSMP}
  {\small
    {\bf{Task}:} greedily choose $K$ of the $\la{C}_n$s
                to best represent \la{d} via $\la{C}_1 \la{h}_1 + \cdots + \la{C}_N \la{h}_N$.\\[\algSkip]
    {\bf{Data and Parameters}:} $\la{d}$ and $\la{C}_n, n \in \{1, \ldots, N \}$ are given. 
                           $K$ iterations.\\[\algSkip]
    {\bf{Initialize}:} Set $k=0$, $I_0 = \emptyset$, $\la{S}_0 = [\mbox{  }]$.\\[\algSkip]
    {\bf{Iterate}:} Set $k=1$ and apply:

    \begin{adamItemize}

      \item{$q_k = \mathop{\mbox{argmax}}_{\,\mbox{\scriptsize{$n \notin I_{k-1}$}}}
               \mbox{  }
        	\la{d}^{\sst{H}} (\la{S}_k^{(n)})(\la{S}_k^{(n)})^{\dagger} \la{d}$, 
                where
                $\la{S}_{k}^{(n)} = [ \la{S}_{k-1}, \la{C}_n ]$}

      \item{$I_k = I_{k-1} \cup \{ q_k \}$}

      \item{$\la{S}_k = [ \la{S}_{k-1}, \la{C}_{q_k} ] $}

      \item{$k = k + 1$.  Terminate loop if $k > K$ or $\la{r}_k = \la{0}$.  
            $I_K$ ends with $T \leq K$ elements.}

    \end{adamItemize}

    {\bf{Compute Weights}:} $\la{x} = \la{S}_K^{\dagger} \la{d}$,
    unstack \la{x} into $\la{h}_{q_1}, \ldots, \la{h}_{q_T}$; set
    remaining $\la{h}_n$s to \la{0}.
  }
\end{algorithm}

\subsection{Iteratively Reweighted Least Squares (IRLS)}

    Having posed three greedy approaches for solving the MSSO problem,
    we now turn our attention toward minimizing (\ref{mmv2_eq4}), the
    relaxed objective function.  Here, the regularization term
    $\lambda$ is used to trade off simultaneous sparsity with residual
    observation error.

    One way to minimize (\ref{mmv2_eq4}) is to use an IRLS-based
    approach \cite{Kar1970}.  To begin, consider manipulating the
    right-hand term of (\ref{mmv2_eq4}) as follows:
    \begin{equation}\label{irls}
      \begin{split}
        \lambda \sum_{n=1}^{N} \Vert \la{h}_n \Vert_2 &
        = \lambda \sum_{n=1}^{N} \frac{\Vert \la{h}_n \Vert_2^2}{\Vert \la{h}_n \Vert_2}
        = \lambda \sum_{n=1}^{N} \frac{|\la{h}_n[1]|^2 + \cdots + |\la{h}_n[P]|^2}{\Vert \la{h}_n \Vert_2}\\
      & \approx \frac{\lambda}{2} \sum_{n=1}^{N} \left[ \la{h}_n^{\ast}[1], \ldots, \la{h}_n^{\ast}[P]  \right]
              \left[  \begin{array}{ccc}
                 \tfrac{2}{\Vert \la{h}_n \Vert_2 + \epsilon} &        &   \\
                                                 & \ddots &   \\
                                                 &        & \tfrac{2}{\Vert \la{h}_n \Vert_2 + \epsilon} \\
                      \end{array}
              \right]
              \left[ \begin{array}{c}
                       \la{h}_n[1] \\
                       \vdots \\
                       \la{h}_n[P]
                     \end{array}
              \right] \\
         & = \frac{\lambda}{2} \sum_{n=1}^{N} \la{h}_n^{\sst{H}} \la{W}_n \la{h}_n \\
         & = \frac{\lambda}{2} \left[ \la{h}_1^{\sst{H}} \cdots \la{h}_N^{\sst{H}} \right]
                         \left[ \begin{array}{ccc}
                                \la{W}_1 &         &           \\
                                         & \ddots  &           \\
                                         &         & \la{W}_N  \\
                                \end{array} \right]
                 \left [ \begin{array}{c}
                             \la{h}_1 \\
                             \vdots  \\
                              \la{h}_N \\
                              \end{array}
           \right ] = \frac{\lambda}{2} \la{h}_{\sst{tot}}^{\sst{H}} \la{W}_{\sst{tot}} \la{h}_{\sst{tot}},
      \end{split}
    \end{equation}
    where $^{\ast}$ is the complex conjugate of a scalar, $\la{W}_n$
    is a $P \times P$ real-valued diagonal matrix whose diagonal
    elements each equal $2 / (\Vert \la{h}_n \Vert_2 + \epsilon)$,
    and $\epsilon$ is some small non-negative value introduced to
    mitigate poor conditioning of the $\la{W}_n$s.  If we fix
    $\la{W}_{\sst{tot}} \in \field{R}^{PN \times PN}$ by computing it
    using some prior estimate of $\la{h}_{\sst{tot}}$, then the
    right-hand term of (\ref{mmv2_eq4}) becomes a quadratic function
    and (\ref{mmv2_eq4}) transforms into a Tikhonov optimization
    \cite{Tik1963, Tik1977}:
    \begin{equation}\label{irls2}
    \mathop{\mbox{min}}_{\mbox{\scriptsize{$\la{h}_{\sst{tot}}$}}}
    \mbox{  }
       \left \{
         \tfrac{1}{2}  \left \Vert \la{d} -
                       \la{C}_{\mbox{\scriptsize{tot}}}
                       \la{h}_{\mbox{\scriptsize{tot}}}
                       \right \Vert_2^2
        +      \tfrac{\lambda}{2}  \la{h}_{\sst{tot}}^{\sst{H}} \la{W}_{\sst{tot}} \la{h}_{\sst{tot}}
       \right \}.
    \end{equation}
    Finally, by performing a change of variables and exploiting the
    properties of $\la{W}_{\sst{tot}}$, we can convert (\ref{irls2})
    into an expression that may be minimized using the robust and
    reliable conjugate-gradient (CG) least-squares solver LSQR
    \cite{LSQR1:82, LSQR2:82}, so named because it applies a QR decomposition
    \cite{Gol1983} when solving the system in the least-squares sense.  LSQR works better in
    practice than several other CG methods \cite{Bjo1978} because it
    restructures the input system via the Lanczos process
    \cite{Lan1950} and applies a Golub-Kahan bidiagonalization \cite{Gol1965}
    prior to solving it.

    To apply LSQR to this problem, we first construct
    $\la{W}_{\sst{tot}}^{\Frac{1}{2}}$ as the element-by-element
    square-root of the diagonal matrix $\la{W}_{\sst{tot}}$ and then
    take its inverse to obtain $\la{W}_{\sst{tot}}^{-\Frac{1}{2}}$.
    Defining $\la{q} = \la{W}_{\sst{tot}}^{\Frac{1}{2}}
    \la{h}_{\sst{tot}}$ and $\la{A} = \la{C}_{\mbox{\scriptsize{tot}}}
    \la{W}_{\sst{tot}}^{-\Frac{1}{2}}$, (\ref{irls2}) becomes:
    \begin{equation}\label{irls3}
    \begin{split}
      \mathop{\mbox{min}}_{\mbox{\scriptsize{$\la{q}$}}}
      \mbox{  }
       \left \{
                \Vert \la{d} - \la{A} \la{q} \Vert_2^2
        +      \lambda \Vert \la{q} \Vert_2^2
       \right \}\\
    \end{split}.
    \end{equation}
    This problem may be solved directly by simply providing $\la{d}$,
    $\lambda$, and $\la{A}$ to the LSQR package because LSQR is
    formulated to solve the exact problem in (\ref{irls3}).  Calling
    LSQR with these variables yields $\la{q}^{\sst{opt}}$, and the
    solution $\la{h}_{\sst{tot}}^{\sst{opt}}$ is backed out via
    $\la{W}_{\sst{tot}}^{-\Frac{1}{2}} \la{q}^{\sst{opt}}$.

    {\em{Algorithm \ref{alg:IRLS}}} outlines how one may iteratively
    apply (\ref{irls3}) to attempt to find a solution that minimizes
    the original cost function, (\ref{mmv2_eq4}).  The technique
    iterates until the objective function decreases by less than
    $\delta$ or the maximum number of iterations, $K$, is exceeded.
    The initial solution estimate is obtained via pseudoinversion of
    $\la{C}_{\sst{tot}}$ (an all-zeros initialization would cause
    $\la{W}_{\sst{tot}}$ to be poorly conditioned).  A line search is
    used to step between the prior solution estimate and the upcoming
    one; this improves the rate of convergence and ensures the
    objective decreases at each step.  This method is global in the
    sense that all $PN$ unknowns are being estimated concurrently per
    iteration.

\begin{algorithm}[t]
    \caption{--- MSSO Iteratively Reweighted Least Squares}
    \label{alg:IRLS}
  {\small
    {\bf{Task}:} Minimize $\left \{\tfrac{1}{2} \left \Vert \la{d} -
    \la{C}_{\mbox{\scriptsize{tot}}} \la{h}_{\mbox{\scriptsize{tot}}}
    \right \Vert_2^2 + \lambda \sum_{n=1}^{N} \Vert \la{h}_n \Vert_2
    \right \}$ using an iterative scheme.\\[\algSkip]
    {\bf{Data and Parameters}:} $\lambda$, $\la{d}$, $\la{C}_{\sst{tot}}$, $\delta$, and $K$ are given.\\[\algSkip]
    {\bf{Initialize}:} Set $k=0$ and $\la{h}_{\sst{tot},k=0} = (\la{C}_{\sst{tot}})^{\dagger} \la{d}$
                     (or e.g. $\la{h}_{\sst{tot},k=0} = \la{1}$).\\[\algSkip]
    {\bf{Iterate}:} Set $k=1$ and apply:

    \begin{adamItemize}

      \item{Create $\la{W}_{\sst{tot}}$ from $\la{h}_{\sst{tot},k-1}$;
            construct $\la{W}_{\sst{tot}}^{\Frac{1}{2}}$, 
            $\la{W}_{\sst{tot}}^{-\Frac{1}{2}}$, 
            and let $\la{A} = \la{C}_{\sst{tot}} \la{W}_{\sst{tot}}^{-\Frac{1}{2}}$.}

      \item{Obtain $\la{q}_{\sst{tmp}}$ by using LSQR to solve
      $\mathop{\mbox{min}}_{\mbox{\scriptsize{$\la{q}$}}}
       \left \{
                \Vert \la{d} - \la{A} \la{q} \Vert_2^2
        +      \lambda \Vert \la{q} \Vert_2^2
       \right \}$.}

      \item{Set $\la{h}_{\sst{tot},\sst{tmp}} = \la{W}_{\sst{tot}}^{-\Frac{1}{2}} \la{q}_{\sst{tmp}}$.}

      \item{Line search: find $\mu_0 \in [0, 1]$ such that $(1-\mu)\la{h}_{\sst{tot},k-1} +
                         \mu \la{h}_{\sst{tot},\sst{tmp}}$ minimizes (\ref{mmv2_eq4}).}

      \item{Set $\la{h}_{\sst{tot},k} = 
            (1-\mu_0)\la{h}_{\sst{tot},k-1} + \mu_0 \la{h}_{\sst{tot},\sst{tmp}}$.}

      \item{$k = k + 1$.  Terminate loop when $k > K$ or (\ref{mmv2_eq4}) decreases by less than $\delta$.}

    \end{adamItemize}

    {\bf{Finalize}:} Unstack the last $\la{h}_{\sst{tot}}$ solution into
    $\la{h}_1, \ldots, \la{h}_{N}$.
  }
\end{algorithm}

\subsection{Row-by-Row Shrinkage (RBRS)}

    The proposed IRLS technique solves for all $PN$ unknowns during
    each iteration, but this is cumbersome when $PN$ is large.  An
    alternative approach is to apply an inner loop that fixes $n$ and
    then iteratively tunes $\la{h}_n$ while holding the other
    $\la{h}_m$s ($m \neq n$) constant; thus only $P$ (rather than
    $PN$) unknowns need to be solved for during each inner iteration.

    This idea inspires the RBRS algorithm.  The term ``row-by-row'' is
    used because each $\la{h}_n$ corresponds to row $n$ of the
    $\la{G}$ matrix in (\ref{mmv2_eq2}), and ``shrinkage'' is used
    because the $\ell_2$ energy of most of the $\la{h}_n$s will
    essentially be ``shrunk'' (to some extent) during each inner
    iteration: when $\lambda$ is sufficiently large and many
    iterations are undertaken, many $\la{h}_n$s will be close to
    all-zeros vectors and only several will be dense and high in
    energy.

    \subsubsection{RBRS for real-valued data} Assume $\la{d}$ and the
    $\la{C}_n$s of (\ref{mmv2_eq4}) are real-valued.  We now seek to
    minimize the function by extending the single-vector sequential
    shrinkage technique of \cite{Ela2006_TransIT} and making
    modifications to (\ref{mmv2_eq4}).  Assume we have prior estimates
    of $\la{h}_1$ through $\la{h}_N$, and that we now desire to update
    only the $j$th vector while keeping the other $N - 1$ fixed.  The
    shrinkage update of $\la{h}_j$ is achieved via:
    \begin{equation}\label{rbrs1}
      \mathop{\mbox{min}}_{\mbox{\scriptsize{$\la{x}$}}}
      \mbox{  }
       \left \{ \tfrac{1}{2}
       \left \Vert \left [ \Sigma_{n=1}^N \la{C}_n \la{h}_n - \la{C}_j  \la{h}_j \right ]
                            + \la{C}_j \la{x} - \la{d} \right \Vert_2^2
                 + \lambda \left \Vert \la{x} \right \Vert_2
       \right \},
    \end{equation}
    where $\Sigma_{n=1}^N \la{C}_n \la{h}_n - \la{C}_j \la{h}_j$ forms
    an approximation of \la{d} using the prior solution coefficients,
    but discards the component contributed by the original $j$th
    vector, replacing the latter via an updated solution vector,
    \la{x}.  The left-hand term promotes a solution \la{x} that keeps
    residual error down, whereas the right-hand term penalizes
    \la{x}s that contain nonzeros.  It is crucial to note that the
    right-hand term does not promote the element-by-element sparsity
    of \la{x}; rather, it penalizes the overall $\ell_2$ energy of
    \la{x}, and thus both sparse and dense \la{x}s are penalized
    equally if their overall $\ell_2$ energies are the same.

    One way to solve (\ref{rbrs1}) is to take its derivative with
    respect to $\la{x}^{\sst{T}}$ and find $\la{x}$ such that the
    derivative equals $\la{0}$.  By doing this and shuffling terms,
    and assuming we have an initial estimate of $\la{x}$, we may solve
    for $\la{x}$ iteratively:
    \begin{equation}\label{rbrs2}
         \la{x}_i = \left [ \la{C}_j^{\sst{T}} \la{C}_j + 
                 \frac{\lambda}{\left \Vert \la{x}_{i-1} \right \Vert_2 + \epsilon}\, \la{I} \right ]^{-1} 
                        \la{C}_j^{\sst{T}} \la{r}_j,
    \end{equation}
    where $\la{r}_j = \la{d} + \la{C}_j \la{h}_j - \Sigma_{n=1}^N
    \la{C}_n \la{h}_n$, $\la{I}$ is a $P \times P$ identity matrix,
    and $\epsilon$ is a small value that avoids ill-conditioned
    results.\footnote{Eq. (\ref{rbrs2}) consists of a direct inversion
    of a $P \times P$ matrix, which is acceptable in this paper
    because all experiments involve $P \leq 10$.  If $P$ is large,
    (\ref{rbrs2}) could be solved via a CG technique
    (e.g., LSQR).}  By iterating on (\ref{rbrs2}) until (\ref{rbrs1})
    changes by less than $\delta_{\sst{inner}}$, we arrive at a
    solution to (\ref{rbrs1}), $\la{x}^{\sst{opt}}$, and this then
    replaces the prior solution vector, $\la{h}_j$.  Having completed
    the update of the $j$th vector, we proceed to update the rest of
    the vectors, looping this outer process $K$ times or until the
    main objective function, (\ref{mmv2_eq4}), changes by less than
    $\delta_{\sst{outer}}$.  {\em{Algorithm \ref{alg:RBRS}}} details
    the entire procedure; unlike IRLS, here we essentially repeatedly
    invert $P \times P$ matrices to pursue a row-by-row solution,
    rather than $PN \times PN$ matrices to pursue a solution that
    updates {\em{all}} rows per iteration.

\begin{algorithm}
    \caption{--- MSSO Row-by-Row Sequential Iterative Shrinkage}
    \label{alg:RBRS}
  {\small
    {\bf{Task}:} Minimize $\left \{\tfrac{1}{2} \left \Vert \la{d} -
    \la{C}_{\mbox{\scriptsize{tot}}} \la{h}_{\mbox{\scriptsize{tot}}}
    \right \Vert_2^2 + \lambda \sum_{n=1}^{N} \Vert \la{h}_n \Vert_2
    \right \}$ using an iterative scheme when all data is {\em{real-valued}}.\\[\algSkip]
    {\bf{Data and Parameters}:} $\lambda$, $\la{d}$, $\la{C}_n$ ($n
    \in \{1, \ldots, N\}$), $\delta_{\sst{outer}}$, $\delta_{\sst{inner}}$, $K$, and $I$ are given.\\[\algSkip]
    {\bf{Initialize}:} Set $k=0$ and $\la{h}_{\sst{tot}} = (\la{C}_{\sst{tot}})^{\dagger} \la{d}$
                     (or e.g. $\la{h}_{\sst{tot}} = \la{1}$), unstack into $\la{h}_1, \ldots, \la{h}_N$.\\[\algSkip]
    {\bf{Iterate}:} Set $k=1$ and apply:

    \begin{adamItemize}

      \item{Sweep over row vectors: set $j=1$ and apply:}

      \begin{adamItemize}

          \item[$\circ$]{Optimize a row vector: set $i=1$ and $\la{x}_0 = \la{h}_j$ and then apply:}

          \begin{adamItemize}

              \item{$\la{x}_i = \left [ \la{C}_j^{\sst{T}} \la{C}_j + 
              \frac{\lambda}{\left \Vert \la{x}_{i-1} \right \Vert_2 + \epsilon} \, \la{I} \right ]^{-1} 
                        \la{C}_j^{\sst{T}} \la{r}_j$, where
                    $\la{r}_j = \la{d} + \la{C}_j \la{h}_j - \Sigma_{n=1}^N \la{C}_n \la{h}_n$.}

              \item{$i = i + 1$.  Terminate when $i > I$ or  
                    (\ref{rbrs1}) decreases by less than $\delta_{\sst{inner}}$.}

          \end{adamItemize}

          \item[$\circ$]{Finalize row vector update: set $\la{h}_j$ to equal the final $\la{x}$.}

          \item[$\circ$]{$j=j+1$.  Terminate loop when $j > N$.}

      \end{adamItemize}

      \item{$k = k + 1$.  Terminate loop when $k > K$ or
            (\ref{mmv2_eq4}) decreases by less than $\delta_{\sst{outer}}$.}

    \end{adamItemize}

    {\bf{Finalize}:} If $\lambda$ was large enough, several
    $\la{h}_n$s should be dense and others close to $\la{0}$.
  }
\end{algorithm}

   \subsubsection{Extending RBRS to complex-valued data} If
   (\ref{mmv2_eq4}) contains complex-valued terms, we may structure
   the row-by-row updates as in (\ref{rbrs1}), but because the
   derivative of the objective function in (\ref{rbrs1}) is more
   complicated due to the presence of complex-valued terms, the simple
   update equation given in (\ref{rbrs2}) is no longer applicable.
   One way to overcome this problem is to turn the complex-valued
   problem into a real-valued one.

   Let us create several real-valued expanded vectors:
   \begin{equation}\label{rbrs_comp2}
     \widetilde{\la{d}}   = \left [  \begin{array}{c}
                                     \mbox{Re}(\la{d}) \\ 
                                     \mbox{Im}(\la{d}) \\
                                     \end{array} \right ] \in \field{R}^{2M}, \mbox{ }\,\,\,
     \widetilde{\la{h}}_n  = \left [  \begin{array}{c}
                                     \mbox{Re}(\la{h}_n) \\ 
                                     \mbox{Im}(\la{h}_n) \\
                                     \end{array} \right ] \in \field{R}^{2P},
   \end{equation}
   as well as real-valued expanded matrices:
   \begin{equation}\label{rbrs_comp3}
      \widetilde{\la{C}}_n = \left[ \begin{array}{cc}
                                    \mbox{Re}(\la{C}_n) & -\mbox{Im}(\la{C}_n) \\
                                    \mbox{Im}(\la{C}_n) &  \mbox{Re}(\la{C}_n) \\
                                    \end{array} \right ] \in \field{R}^{2M \times 2P}.
   \end{equation}
   Due to the structure of (\ref{rbrs_comp2}, \ref{rbrs_comp3}) and the fact
   that $\Vert \la{h}_n \Vert_2$ equals $\Vert \widetilde{\la{h}}_n
   \Vert_2$, the following optimization is {\em{equivalent}} to
   (\ref{mmv2_eq4}):
   \begin{equation}\label{rbrs_comp4}
      \mathop{\mbox{min}}_{\mbox{\scriptsize{$\widetilde{\la{h}}_1, \ldots, \widetilde{\la{h}}_N$}}} 
      \mbox{  }
      \left \{
        \tfrac{1}{2}  \left \Vert \widetilde{\la{d}} -
                      \sum_{n=1}^N \widetilde{\la{C}}_n
                      \widetilde{\la{h}}_n
                      \right \Vert_2^2
       +     \lambda  \sum_{n=1}^{N} \Vert \widetilde{\la{h}}_n \Vert_2
      \right \}.
   \end{equation}
   This means we may apply RBRS to complex-valued scenarios by
   substituting the $\widetilde{\la{h}}_n$s for the $\la{h}_n$s and
   $\widetilde{\la{C}}_n$s for the $\la{C}_n$s in (\ref{rbrs1},
   \ref{rbrs2}) and {\em{Algorithm \ref{alg:RBRS}}}.  [Eq.~(\ref{rbrs2})
   becomes an applicable update equation because (\ref{rbrs1}) will
   consist of only real-valued terms and the derivative calculated
   earlier is again applicable.]  Finally, after running the algorithm to
   obtain finalized $\widetilde{\la{h}}_n$s, we may simply restructure
   them into complex $\la{h}_n$s.

\subsection{Column-by-Column Shrinkage (CBCS)}

    Here we propose a dual of RBRS---a technique that sequentially
    updates the {\em{columns}} of $\la{G}$ (i.e., the $\la{g}_p$s) in
    (\ref{mmv2_eq1}, \ref{mmv2_eq2}) rather than its rows (the
    $\la{h}_n$s).  Interestingly, we will show that this approach
    yields a {\em{separable}} optimization and reduces the overall
    problem to simply repeated {\em{element-by-element}} shrinkages of
    each $\la{g}_p$.

    \subsubsection{CBCS for real-valued data} Assume the $\la{g}_p$s,
    $\la{F}_p$s, and $\la{d}$ in (\ref{mmv2_eq2}) are real-valued and
    that we have prior estimates of the $\la{g}_p$s.  Let us consider
    updating the $p$th vector while keeping the other $P - 1$ fixed.
    This reduces (\ref{mmv2_eq2}) to
    \begin{equation}\label{cbcs1}
      \mathop{\mbox{min}}_{\mbox{\scriptsize{$\la{x}$}}}
      \mbox{  }
       \left  \{ \tfrac{1}{2}
                 \left \Vert \la{r} - \la{F}_p \la{x}     \right \Vert_2^2 +
         \lambda \sum_{n=1}^{N} \sqrt{ (\la{x}[n])^2 + \la{b}[n] }
       \right \},
    \end{equation}
    where $\la{x}$ will be the update of $\la{g}_p$, and \la{r} and \la{b} are as follows:
    \begin{equation}\label{cbcs_extra1}
       \la{r} = \la{d} + \la{F}_p \la{g}_p - \sum_{q=1}^{P} \la{F}_q \la{g}_q,
    \end{equation}
    and
    \begin{equation}\label{cbcs_extra2}
       \la{b}[n] = -(\la{g}_p[n])^2 + \sum_{q=1}^{P} (\la{g}_q[n])^2, \mbox{ for } n = 1, \ldots, N.
    \end{equation}
    If the $\la{b}[n]$s were not present, (\ref{cbcs1}) would reduce to the
    standard problem iterated shrinkage is intended to solve
    \cite{Ela2006_TransIT,Ela2006}.

    Now let us apply a proximal relaxation \cite{Dau2004, Fig2007,
    Com2007} to (\ref{cbcs1}) and seek a solution $\la{x} \in
    \field{R}^{N}$ as a shrinkage update of $\la{g}_p$:
    \begin{equation}\label{cbcs2}
      \mathop{\mbox{min}}_{\mbox{\scriptsize{$\la{x}$}}}
      \mbox{  }
       \left  \{ \tfrac{1}{2} \left(
                 \left \Vert \la{r} - \la{F}_p \la{x}     \right \Vert_2^2 +
          \alpha \left \Vert \la{x} - \la{g}_p            \right \Vert_2^2 -
                 \left \Vert \la{F}_p (\la{x} - \la{g}_p) \right \Vert_2^2 \right) + 
         \lambda \sum_{n=1}^{N} \sqrt{ (\la{x}[n])^2 + \la{b}[n] }
       \right \},
    \end{equation}
    where $\alpha$ is chosen such that $\alpha\la{I} -
    \la{F}_p^{\sst{T}} \la{F}_p$ is positive definite (e.g., $\alpha$
    may be set to the maximum singular value of $\la{F}_p$).  The idea
    here is to replace $\la{g}_p$ with the solution \la{x} and then
    iterate this procedure, repeatedly solving (\ref{cbcs2}).  This
    ultimately yields an updated \la{x} that globally minimizes
    (\ref{cbcs1}) because the proximal method is guaranteed to arrive
    at a local minimum \cite{Dau2004, Fig2007} and (\ref{cbcs1})
    itself is convex.  Having obtained \la{x}, we perform an update,
    $\la{g}_p = \la{x}$, and then repeat the overall process for the
    next $\la{g}_p$, and so forth.  Additionally, we add a layer of
    iteration on top of this column-by-column sweep, optimizing each
    of the $P$ vectors a total of $K$ times.

    The only obstacle that remains in order for us to implement the
    entire algorithm is an efficient way to solve (\ref{cbcs2}).  We
    pursue such an approach by first expanding the terms of
    (\ref{cbcs2}):
    \begin{equation}\label{cbcs3}
      \mathop{\mbox{min}}_{\mbox{\scriptsize{$\la{x}$}}}
      \mbox{  }
        \left \{
        \la{c} + \la{v}^{\sst{T}} \la{x} 
               + \frac{\alpha}{2} \la{x}^{\sst{T}} \la{x}
               + \lambda \sum_{n=1}^{N} \sqrt{ (\la{x}[n])^2 + \la{b}[n] }
        \right \},
    \end{equation}
    where $\la{c} = \tfrac{1}{2} \la{r}^{\sst{T}}\la{r} +
    \tfrac{\alpha}{2}\la{g}_p^{\sst{T}}\la{g}_p
    -\tfrac{1}{2}\la{g}_p^{\sst{T}}\la{F}_p^{\sst{T}}\la{F}_p\la{g}_p$
    and $\la{v} = \la{F}_p^{\sst{T}}\la{F}_p\la{g}_p - \alpha\la{g}_p
    - \la{F}_p^{\sst{T}}\la{r}$.  Since $\la{c}$ is constant, we may
    ignore it in the optimization.  Upon closer inspection, we see
    that (\ref{cbcs3}) is a {\em{separable}} problem and that the
    individual scalar elements of $\la{x}$ may be optimized
    independently.  For the $n$th element of \la{x}, (\ref{cbcs3})
    simplifies to:
    \begin{equation}\label{cbcs4}
      \mathop{\mbox{min}}_{\mbox{\scriptsize{$\la{x}[n]$}}}
      \mbox{  }
        \left \{
                 \la{v}[n] \la{x}[n] + 
                 \frac{\alpha}{2} (\la{x}[n])^2 +
                 \lambda \sqrt{ (\la{x}[n])^2 + \la{b}[n] }
        \right \},
    \end{equation}
    Having burrowed down to an element-by-element problem, all that
    remains is to efficiently solve (\ref{cbcs4}).  One approach is to
    compute the derivative of its objective with respect to
    $\la{x}[n]$ and find $\la{x}[n]$ such that the derivative equals
    zero.  The derivative equals the following nonlinear scalar equation:
    \begin{equation}\label{cbcs_extra3}
         \la{v}[n] + \la{x}[n] \left ( \alpha  + 
                                        \frac{\lambda}{\sqrt{(\la{x}[n])^2 + \la{b}[n]}} 
                               \right ).
    \end{equation}
    Setting the derivative in (\ref{cbcs_extra3}) to zero and assuming we have
    an initial estimate of $\la{x}[n]$, we may solve for $\la{x}[n]$
    iteratively as follows:
    \begin{equation}\label{cbcs5}
         (\la{x}[n]_i) =  -\la{v}[n] \left(\alpha + 
           \frac{\lambda}{\sqrt{(\la{x}[n]_{i-1})^2 + \la{b}[n] + \epsilon}}\right)^{-1},
    \end{equation}
    where $\epsilon$ is simply a small value that avoids
    ill-conditioned scenarios.
    
    We may now formulate CBCS as {\em{Algorithm \ref{alg:CBCS}}}.  As
    we seek to update a fixed $\la{g}_1$, note how we iteratively tune
    its $N$ elements, one at a time, via (\ref{cbcs5}), but instead of
    moving on immediately to update $\la{g}_2$, we update $\la{g}_1$,
    $\la{r}$, $\la{v}$, and $\la{b}$, and tune over the elements of
    $\la{g}_1$ yet again, doing this repeatedly until the per-vector
    objective, (\ref{cbcs1}), stops decreasing---only then moving on
    to $\la{g}_2$.  Empirically, we find this greatly speeds up the
    rate at which the $\la{g}_p$s converge to a simultaneously sparse
    solution, but unfortunately, even with this extra loop, CBCS still
    requires excessive iterations for larger problems (see
    Sec.~\ref{sec:experiments}).  Similarly to RBRS in {\em{Algorithm
    \ref{alg:RBRS}}}, note how the inner loops are cut off when the
    objective function stops decreasing to within some small value
    $\delta$ or some fixed number of iterations has been exceeded.

\begin{algorithm}[t]
    \caption{--- MSSO Column-by-Column Sequential Iterative Shrinkage}
    \label{alg:CBCS}
  {\small
    {\bf{Task}:} Minimize $ \left \{ \tfrac{1}{2} \left \Vert \la{d} -
    \la{F}_{\mbox{\scriptsize{tot}}} \la{g}_{\mbox{\scriptsize{tot}}}
    \right \Vert_2^2 + \lambda \left \Vert \la{G} \right
    \Vert_{\mbox{\scriptsize{S}}} \right \}$ when all data is
    {\em{real-valued}}.\\[\algSkip]
    {\bf{Data and Parameters}:} $\lambda$, $\la{d}$, $\la{F}_p, p \in
    \{1, \ldots, P\}$, $\delta_{\sst{tot}}$, $\delta_{\sst{vec}}$,
    $\delta_{\sst{elem}}$, $K$, $J$, $I$ are given.\\[\algSkip]
    {\bf{Initialize}:} $\la{g}_{\sst{tot}} =
    (\la{F}_{\sst{tot}})^{\dagger} \la{d}$; split into $\la{g}_1,
    \ldots, \la{g}_P$; set $\alpha = $ max.~sing.~val.~among
    $\la{F}_p$s.\\[\algSkip]
    {\bf{Iterate}:} Set $k=1$ and apply:

    \begin{adamItemize}

      \item{Sweep over column vectors: set $p=1$ and apply:}

      \begin{adamItemize}

          \item[$\circ$]{Optimize a column vector: set $j=1$ and apply:}

          \begin{adamItemize}

              \item{Construct $\la{r} = \la{d} + \la{F}_p \la{g}_p -
              \sum_{q=1}^{P} \la{F}_q \la{g}_q$.}

              \item{Construct $\la{b}[l] = -(\la{g}_p[l])^2 +
              \sum_{q=1}^{P} (\la{g}_q[l])^2$, for $l = 1, \ldots,
              N$.}

              \item{Construct $\la{v} =
              \la{F}_p^{\sst{T}}\la{F}_p\la{g}_p - \alpha\la{g}_p -
              \la{F}_p^{\sst{T}}\la{r}$.}

              \item{Set $\la{x}_0 = \la{g}_p$.}

              \item{Sweep over column elements: set $n=1$ and apply:}

                   \begin{adamItemize}

                        \item[$\circ$]{Optimize $n$th element of \la{x}: set $i=1$ and apply:}

                        \begin{adamItemize}

                              \item{$(\la{x}[n]_i) = -\la{v}[n]
                              \left(\alpha +
                              \frac{\lambda}{\sqrt{(\la{x}[n]_{i-1})^2
                              + \la{b}[n] + \epsilon}}\right)^{-1}$.}\\

                              \item{$i = i+1$.  Stop if $i > I$ or
                              (\ref{cbcs4}) decreases by less than
                              $\delta_{\sst{elem}}$.}

                        \end{adamItemize}

                        \item[$\circ$]{$n=n+1$.  Terminate when $n > N$.}

                   \end{adamItemize}

          \item{Update column vector: set $\la{g}_p$ to equal the final $\la{x}$.}

          \item{$j = j + 1$.  Terminate when $j > J$ or  
               (\ref{cbcs1}) decreases by less than $\delta_{\sst{vec}}$.}

          \end{adamItemize}

      \item[$\circ$]{$p = p + 1$.  Terminate when $p > P$.}

      \end{adamItemize}

      \item{$k = k + 1$.  Terminate loop when $k > K$ or
            (\ref{mmv2_eq2}) decreases by less than $\delta_{\sst{tot}}$.}

    \end{adamItemize}

    {\bf{Finalize}:} If $\lambda$ was sufficiently large, $\la{g}_1,
    \ldots, \la{g}_{P}$ should be simultaneously sparse.
  }
\end{algorithm}

    \subsubsection{Extending CBCS to complex-valued data} If
    (\ref{mmv2_eq2}) contains complex-valued terms, we may structure
    the column-by-column updates as in (\ref{cbcs1}, \ref{cbcs2}),
    but the expansion and derivative of the latter equation's
    objective function does not lend itself to the simple update
    equations given in (\ref{cbcs3}, \ref{cbcs4}, \ref{cbcs5}).  One
    way to overcome this problem is to turn the complex-valued problem
    into a real-valued one.  This approach is not equivalent to the
    one used to extend RBRS to complex data.

    First we stack the target vector, \la{d}, into a real-valued vector:
    \begin{equation}\label{cbcs_comp1}
     \widetilde{\la{d}}   = \left [  \begin{array}{c}
                                     \mbox{Re}(\la{d}) \\ 
                                     \mbox{Im}(\la{d}) \\
                                     \end{array} \right ] \in \field{R}^{2M}, \mbox{ }
    \end{equation}
    and then {\em{split}}, rather than stack, the unknown vectors into $2P$ new vectors:
    \begin{equation}\label{cbcs_comp2}
       \la{g}_p^{\sst{(Re)}} = \mbox{Re}(\la{g}_p) \in \field{R}^{N}, \,\,\,
       \la{g}_p^{\sst{(Im)}} = \mbox{Im}(\la{g}_p) \in \field{R}^{N}, \mbox{ for } p = 1, \ldots, P.
   \end{equation}
   We then aggregate these vectors into $\widetilde{\la{G}} =
   [\la{g}_1^{\sst{(Re)}}, \la{g}_1^{\sst{(Im)}}, \ldots,
   \la{g}_P^{\sst{(Re)}}, \la{g}_P^{\sst{(Im)}}]$. Next, we split each
   $\la{F}_p$ into two separate matrices, for $p = 1, \ldots, P$:
   \begin{equation}\label{cbcs_comp3}
      \la{F}_p^{\sst{(A)}} = \left[ \begin{array}{c}
                             \mbox{Re}(\la{F}_p) \\
                             \mbox{Im}(\la{F}_p) \\
                             \end{array} \right ] \in \field{R}^{2M \times N}, \mbox{ }
      \la{F}_p^{\sst{(B)}} = \left[ \begin{array}{c}
                             -\mbox{Im}(\la{F}_p) \\
                              \mbox{Re}(\la{F}_p) \\
                            \end{array} \right ] \in \field{R}^{2M \times N},
   \end{equation}
   yielding $2P$ new real-valued matrices.

   Due to the structure of (\ref{cbcs_comp1}, \ref{cbcs_comp2},
   \ref{cbcs_comp3}), the following optimization is {\em{equivalent}}
   to (\ref{mmv2_eq2}):
   \begin{equation}\label{cbcs_comp4}
     \begin{split}
     \mathop{\mbox{min}}_{\mbox{\scriptsize{$\widetilde{\la{G}}$}}} \mbox{ }
       &\Bigg\{  \tfrac{1}{2}  \left \Vert \widetilde{\la{d}} -
                      \sum_{p=1}^P \la{F}_p^{\sst{(A)}} \la{g}_p^{\sst{(Re)}} -
                      \sum_{p=1}^P \la{F}_p^{\sst{(B)}} \la{g}_p^{\sst{(Im)}}
                      \right \Vert_2^2 \\
      & \,\,\,\, + \lambda  \sum_{n=1}^{N} \sqrt{  \sum_{p=1}^P (\la{g}_p^{\sst{(Re)}}[n])^2 + 
                                            \sum_{p=1}^P (\la{g}_p^{\sst{(Im)}}[n])^2   } \Bigg\}.
     \end{split}
   \end{equation}
   The equivalence arises because the first and second terms of
   (\ref{cbcs_comp4}) are equivalent to $\frac{1}{2} \Vert \la{d} -
   \la{F}_{\sst{tot}} \la{g}_{\sst{tot}} \Vert_2^2$ and $\Vert \la{G}
   \Vert_{\sst{S}}$ in (\ref{mmv2_eq2}), respectively.

   This means we may apply CBCS to complex-valued problems by
   performing column-by-column optimization over the $2P$ real-valued
   unknown vectors.  This works because CBCS will
   pursue solutions where the $\la{g}_1^{\sst{(Re)}},
   \la{g}_1^{\sst{(Im)}}, \ldots, \la{g}_P^{\sst{(Re)}},
   \la{g}_P^{\sst{(Im)}}$ vectors are simultaneously sparse, which is
   equivalent to pursuing simultaneously sparse $\la{g}_1, \ldots,
   \la{g}_P$s.  After running CBCS on the $2P$ vectors, one may simply
   restructure them into $P$ complex-valued $\la{g}_p$s.

   Finally, let us set $P = 1$ and thus consider the case of
   single-vector sparse approximation.  The above derivations show
   that {\em{seeking a single sparse complex-valued vector is
   equivalent to seeking two {\bf{simultaneously sparse}} real-valued
   vectors}}.  In other words, single-vector sparse approximation of a
   complex vector readily maps to the MSSO problem, increasing the
   applicability of algorithms that solve the latter.

\subsection{Second-Order Cone Programming (SOCP)}
\label{subsec:socp}

    We now propose a seventh and final algorithm to solve the MSSO
    problem as given in (\ref{mmv2_eq2}).  We branch away from the
    shrinkage approaches that operate on individual columns or rows of
    the \la{G} matrix and again seek to concurrently estimate all $PN$
    unknowns.  Rather than using an IRLS technique, however, we pursue
    a second-order cone programming approach, motivated by the fact
    that second-order cone programs may be solved via efficient
    interior point algorithms \cite{SeDuMi, Toh1999} and are able to
    encapsulate conic, convex-quadratic \cite{Nem2001}, and linear
    constraints.  (Quadratic programming is not an option because the
    $\la{g}_p$s, $\la{F}_p$s, and \la{d} may be complex.)

    Second-order conic constraints are of the form $\la{a} = [a_1,
    \ldots, a_N]^{\sst{T}}$ such that
    \begin{equation}\label{socp0}
       \Vert [a_1, \ldots, a_{N-1}]^{\sst{T}} \Vert_2 \leq a_N.
    \end{equation}
    The generic format of an SOC program is
    \begin{equation} \label{socp1}
     \mathop{\mbox{min}}_{\sst{\la{x}}} \mbox{ }  \la{c}^{\mbox{\scriptsize{T}}}\la{x}
     \,\,\,\text{   s.t. } \la{Ax} = \la{b}  \text{ and } \la{x} \in \la{K},
    \end{equation}
    where $\la{K} = \field{R}_{+}^{N} \times \la{L}_1 \times \cdots
    \times \la{L}_{N}$, $\field{R}_{+}^{N}$ is the $N$-dimensional
    positive orthant cone, and the $\la{L}_n$s are second-order cones
    \cite{Nem2001}.
    To convert (\ref{mmv2_eq2}) into the SOC format, we first write
    \begin{equation}\label{socp2}
    \begin{split}
    \mathop{\mbox{min}}_{\sst{\la{G}}} 
     \mbox{ } \big \{ & \tfrac{1}{2} s + \lambda \la{1}^{\sst{T}}\la{t} \big \} \\
     \text{ s.t. } \la{z} = \la{d}_{\sst{tot}} - & \la{F}_{\sst{tot}} \la{g}_{\sst{tot}}
        \text{ and } \Vert \la{z} \Vert_2^2 \leq s \,\, \\
     \text{ and }  \Vert [\mbox{Re}(\la{g}_1[n]), \mbox{Im}(\la{g}_1[n]),
                & \ldots, \,\mbox{Re}(\la{g}_P[n]),\mbox{Im}(\la{g}_P[n])]^{\sst{T}} \Vert_2 \leq t_n \\
    \end{split}
    \end{equation}
    where $n \in \{1, \ldots, N\}$ and $\la{t} = [t_1, \ldots,
    t_N]^{\sst{T}}$.  The splitting of the complex elements of
    the $\la{g}_p$s mimics the approach used when extending CBCS to
    complex data, and (\ref{socp2}) makes the objective function
    linear, as required.  Finally, in order to represent the $\Vert
    \la{z} \Vert_2^2 \leq s$ inequality in
    terms of second-order cones, an additional step is needed.  Given
    that $s = \frac{1}{4}(s+1)^2 - \frac{1}{4}(s-1)^2$, the inequality
    may be rewritten as $\la{z}^{\sst{H}}\la{z} + \frac{1}{4}(s-1)^2
    \leq \frac{1}{4}(s+1)^2$ and then expressed as a conic constraint:
    $\Vert [\la{z}^{\sst{T}}, \frac{1}{2}(s-1)]^{\sst{T}}\Vert_2 \leq \frac{1}{2}(s+1)$
    \cite{Nem2001,Mal2003}.  Applying these changes yields
    \begin{equation} \label{socp3}
    \begin{split}
    \mathop{\mbox{min}} \mbox{ } \big \{ & \tfrac{1}{2} s + \lambda \la{1}^{\sst{T}}\la{t} \big \} \\
     \text{ s.t. } \la{z} = \la{d}_{\sst{tot}} - \la{F}_{\sst{tot}} & \la{g}_{\sst{tot}}
        \text{ and } \Vert [\la{z}^{\sst{T}}, u]^{\sst{T}} \Vert_2 \leq v, \,\,\,\,\, \\
     u = \Frac{(s-1)}{2}, &\, v = \Frac{(s+1)}{2}, s \geq 0, \\
     \text{ and }  \Vert [\mbox{Re}(\la{g}_1[n]), \mbox{Im}(\la{g}_1[n]),
                & \ldots, \,\mbox{Re}(\la{g}_P[n]),\mbox{Im}(\la{g}_P[n])]^{\sst{T}} \Vert_2 \leq t_n, \\
    \end{split}
    \end{equation}
    which is a fully-defined SOC program that may be implemented and
    solved numerically.  There is no {\em{Algorithm}} pseudocode for
    this technique because having set up the variables in
    (\ref{socp3}), one may simply plug them into an SOCP solver. In
    this paper we implement (\ref{socp3}) in SeDuMi
    (Self-Dual-Minimization) \cite{SeDuMi}, a free software package
    consisting of {\tt{MATLAB}} and {\tt{C}} routines.

\section{Experiments and Results}
\label{sec:experiments}
   Our motivation for solving MSSO sparse approximation problems
   comes from MRI RF excitation pulse design.
   Due to the NP-hardness of the problem (\ref{msso_nphard}),
   there is no reasonable way to check the accuracy of approximate
   solutions to these problem instances obtained with the
   algorithms introduced here.  Thus, before turning to the MRI RF
   excitation pulse design problem in Sec.~\ref{subsec:e3},
   we present several synthetic experiments.  These experiments
   allow comparisons amongst algorithms and also empirically reveal
   some properties of the relaxation (\ref{mmv2_eq2}).
   Theoretical exploration of this relaxation is also merited
   but is beyond the scope of this manuscript.

   All experiments are performed on a Linux server with a 3.0-GHz
   Intel Pentium IV processor.  The system has 16 gigabytes of random
   access memory, ample to ensure that none of the algorithms require
   the use of virtual memory and to avoid excessive hard drive paging.
   MP, LSMP, IRLS, RBRS, CBCS are implemented in {\tt{MATLAB}},
   whereas SOCP is implemented in SeDuMi.  The runtime of any method
   could be reduced significantly by implementing it in a completely
   compiled format such as {\tt{C}}.  Note: OMP is not evaluated
   because its performance always falls in between that of MP and
   LSMP\@.

\subsection{Sparsity Profile Estimation in a Noiseless Setting}
\label{subsec:e1}

   \subsubsection{Overview} We now evaluate how well the algorithms of
   Sec.~\ref{sec:algorithms} estimate sparsity profiles when the
   underlying $\la{g}_p$s are each strictly and simultaneously
   $K$-sparse and the observation $\la{d}$ of (\ref{mmv2_eq1}) is
   known exactly and not corrupted by noise.  This corresponds to a
   high-SNR source localization scenario where the sparsity profile
   indicates locations of emitters and our goal is to find the
   locations of these emitters \cite{Joh1993, Kri1996, Mal2003,
   Mal2005}.  Our goal is to get an initial grasp of the challenges of
   solving the MSSO problem.

   We synthetically generate real-valued sets of $\la{F}_p$s and
   $\la{g}_p$s in (\ref{mmv2_eq1}), apply the algorithms, and record
   the fraction of correct sparsity profile entries recovered by each.
   We vary $M$ in (\ref{mmv2_eq1}) to see how performance at solving
   the MSSO problem varies when the $\la{F}_p$s are underdetermined
   vs.~overdetermined and also vary $P$ to see how rapidly performance
   degrades as more system matrices and vectors are employed.

   \subsubsection{Details}\label{subsubsec:e1} For all trials, we fix
   $N=30$ in (\ref{mmv2_eq1}) and $K=3$, which means each $\la{g}_p$
   vector consists of thirty elements, three of which are nonzero.  We
   consider $P \in \{1, 2, \ldots, 8\}$, and $M \in \{10, 15, \ldots,
   40\}$.  For each of the fifty-six fixed $(M, P)$ pairs, we create
   50 random instances of (\ref{mmv2_eq1}).  Each of the 2,800
   instances is constructed and evaluated as follows:
   \begin{adamItemize2}

   \item Pick a $K$-element subset of $\{ 1, \ldots, N \}$ uniformly
   at random.  This is the sparsity profile.

   \item Create $P$ total $N$-element vectors, the $\la{g}_p$s.  The
   $K$ elements of each that correspond to the sparsity profile are
   filled in with draws from a Gaussian $\sim \mathcal{N}(0, 1)$
   distribution; all other elements are set to zero.

   \item Create $P$ total $M \times N$ matrices, the $\la{F}_p$s.
   Each element of each matrix is determined by drawing from
   $\mathcal{N}(0, 1)$; each column of each matrix is normalized to
   have unit $\ell_2$ energy.

   \item Compute $\la{d} = \sum_{p=1}^{P} \la{F}_p \la{g}_p$. Shuffle
   $\la{F}_p$s and $\la{g}_p$s into $\la{C}_n$s and $\la{h}_n$s via
   (\ref{C_n}, \ref{h_n}).

   \item Apply the algorithms:

       \begin{adamItemize2}

       \item[$\circ$] MP, LSMP: iterate until $K$ elements are chosen
       or the residual approximation is $\la{0}$.  If less than $K$
       terms are chosen, this hurts the recovery score.

       \item[$\circ$] IRLS, RBRS, CBCS, SOCP: approximate a $\lambda$
       oracle: proxy for a good choice of $\lambda$ by looping over
       roughly seventy $\lambda$s in $[0, 2]$, running the given
       algorithm each time.  This sweep over $\lambda$ results in
       high-energy, dense solutions through negligible-energy,
       all-zeros solutions.  For each of the estimated
       $\widehat{\la{g}}_{\sst{tot}}$s (that vary with $\lambda$),
       estimate a sparsity profile by noting the largest $\ell_2$
       energy rows of the associated $\widehat{\la{G}}$
       matrix.\footnote{For example, if the true sparsity profile is
       $\{1, 2, 9\}$ and the largest $\ell_2$ energy rows of
       $\widehat{\la{G}}$ are $\{2, 7, 8\}$, then the fraction of
       recovered sparsity profile terms equals $\frac{1}{3}$.  Now
       suppose only two rows of $\widehat{\la{G}}$ have nonzero energy
       and the profile estimate is only $\{7, 8\}$.  The fraction
       recovered is now zero.} Remember the highest fraction
       recovered across all $\lambda$s.

       \end{adamItemize2}

   \end{adamItemize2}
   After performing the above steps, we average the results of the 50
   trials associated with each fixed $(M, P)$ to yield the average
   fraction of recovered sparsity profile elements.

   \subsubsection{Results} 
\begin{figure}
     \begin{center} \small \begin{tabular}{ccc}
       \epsfig{figure=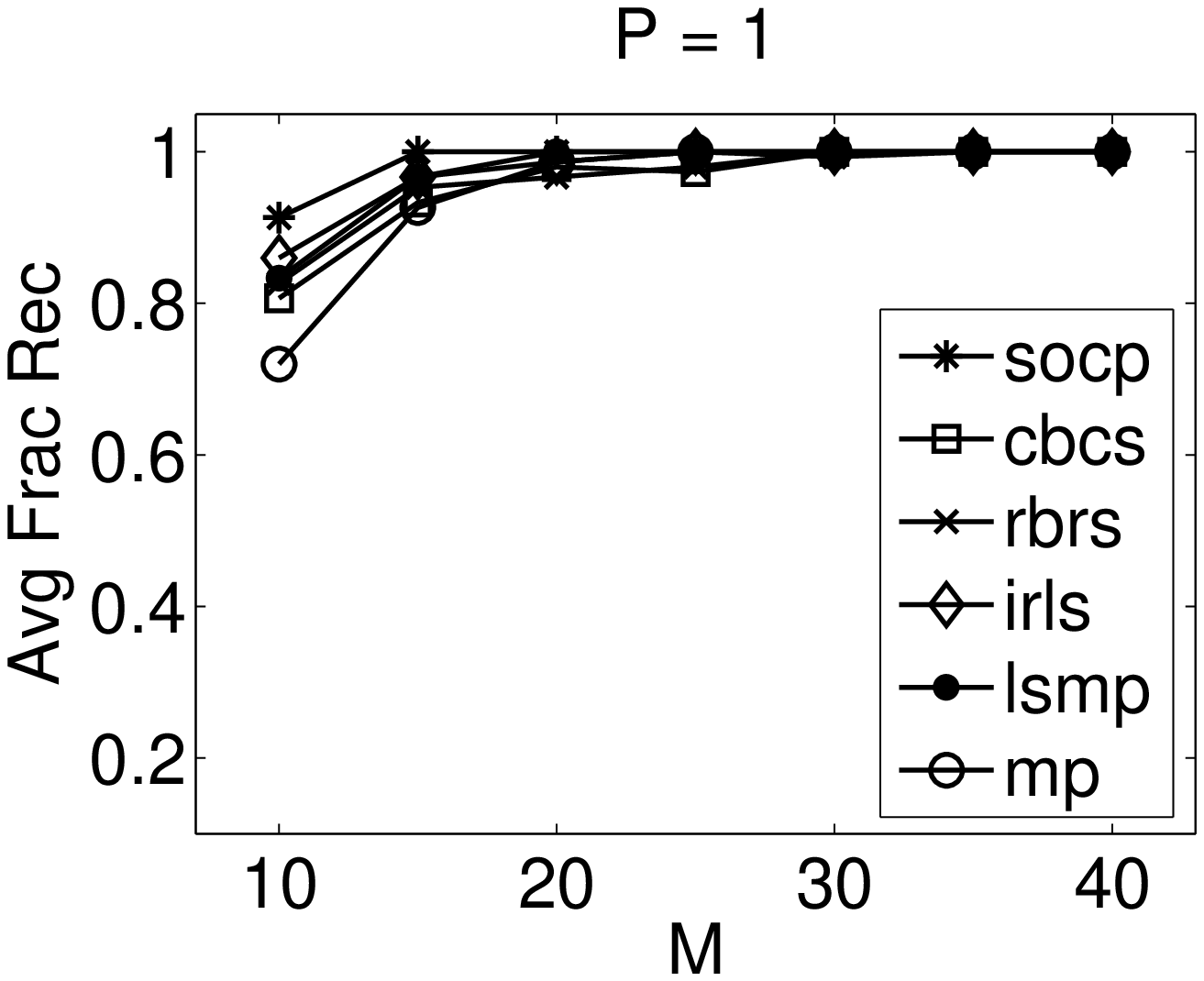,width=\widthC} &
       \epsfig{figure=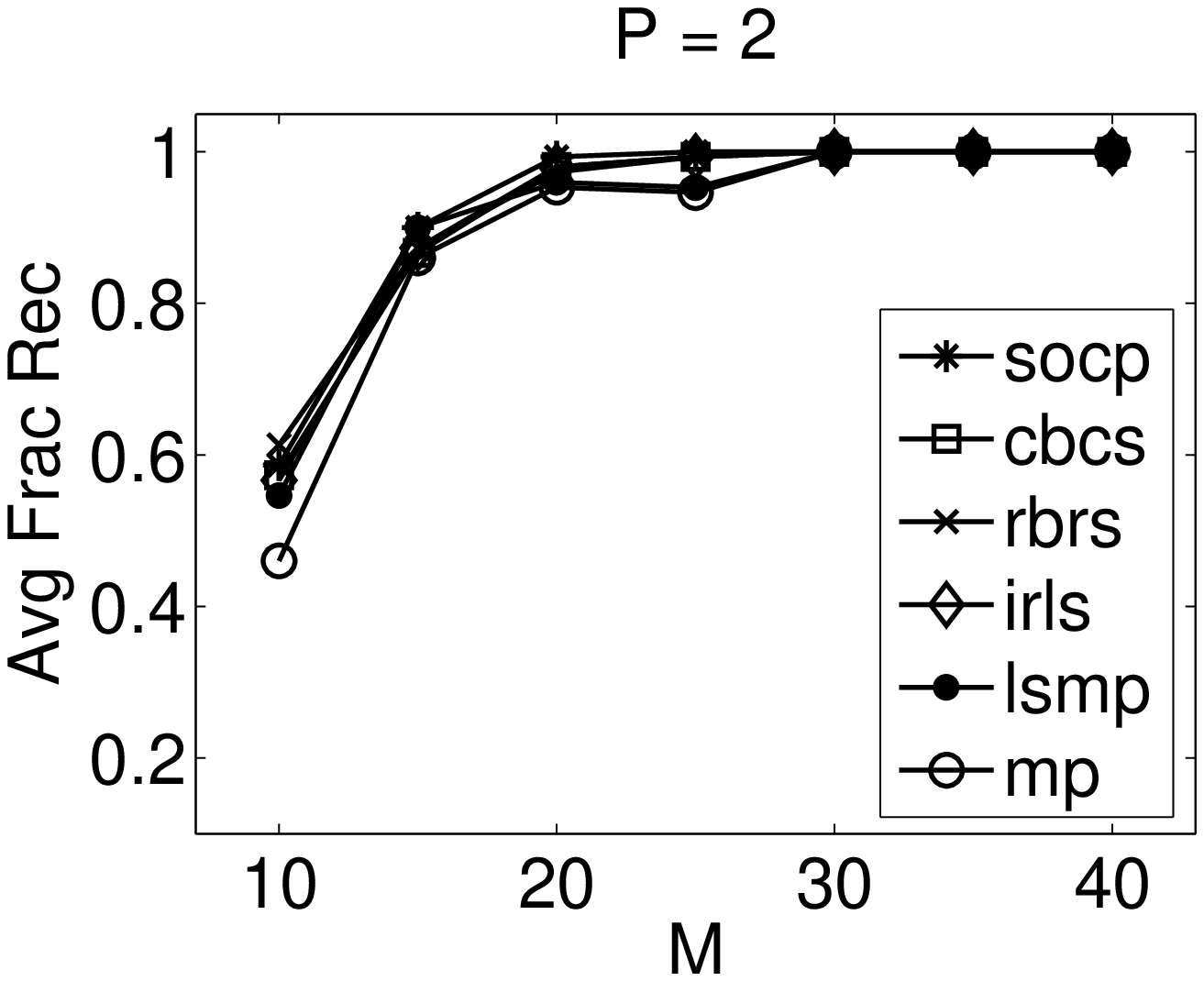,width=\widthC} &
       \epsfig{figure=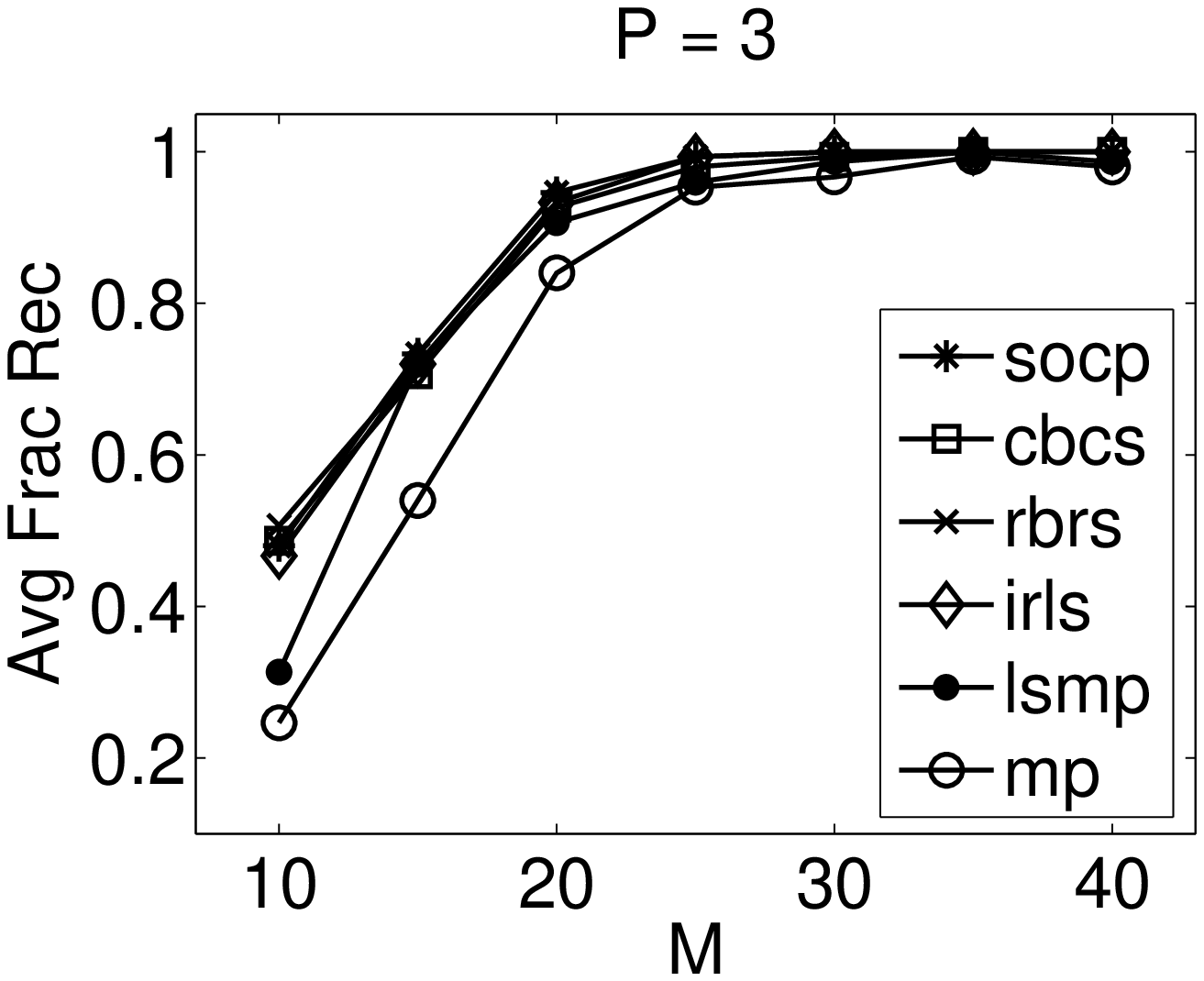,width=\widthC} 
     \end{tabular}

     \begin{tabular}{cccc}
       \epsfig{figure=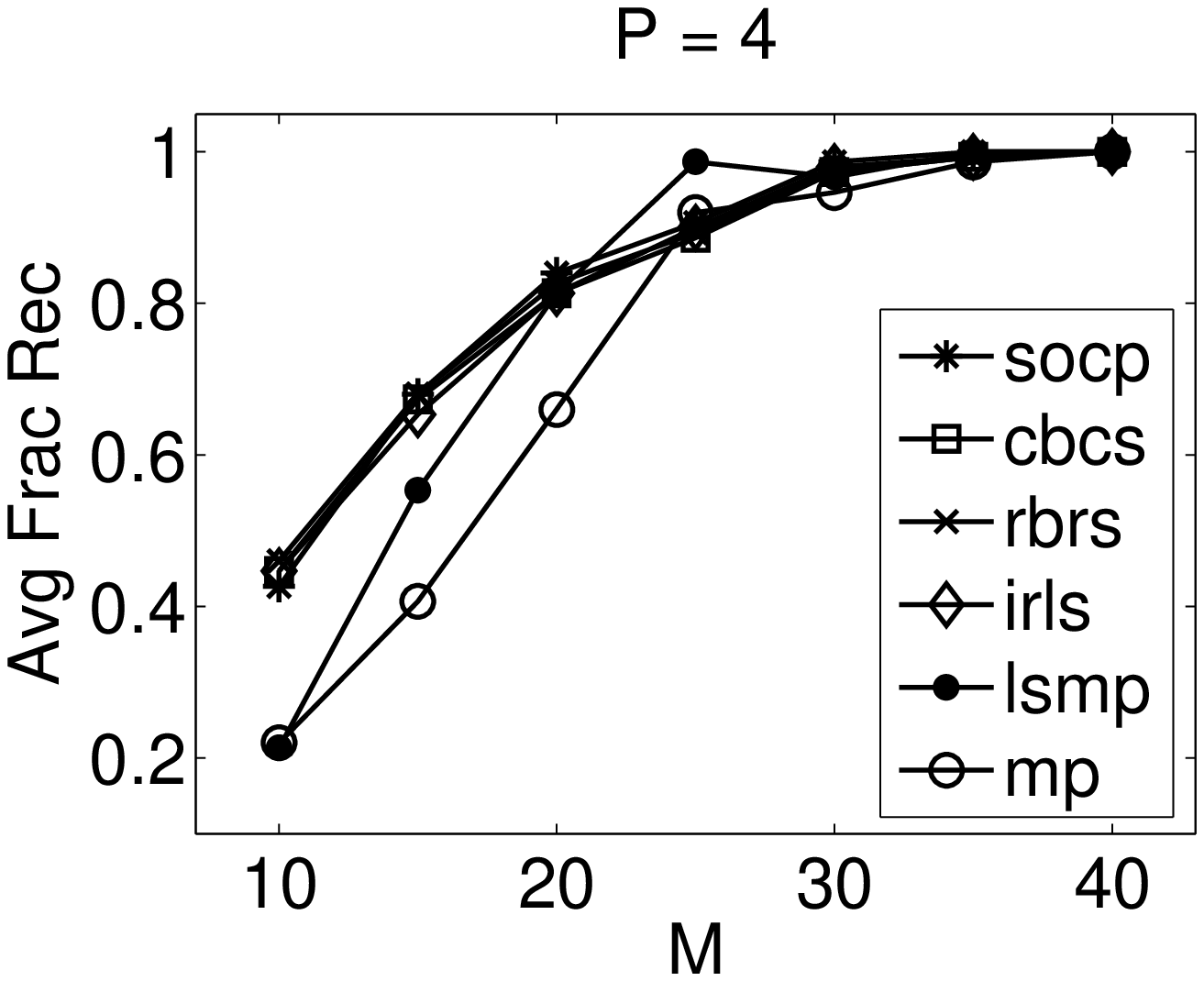,width=\widthC} &
       \epsfig{figure=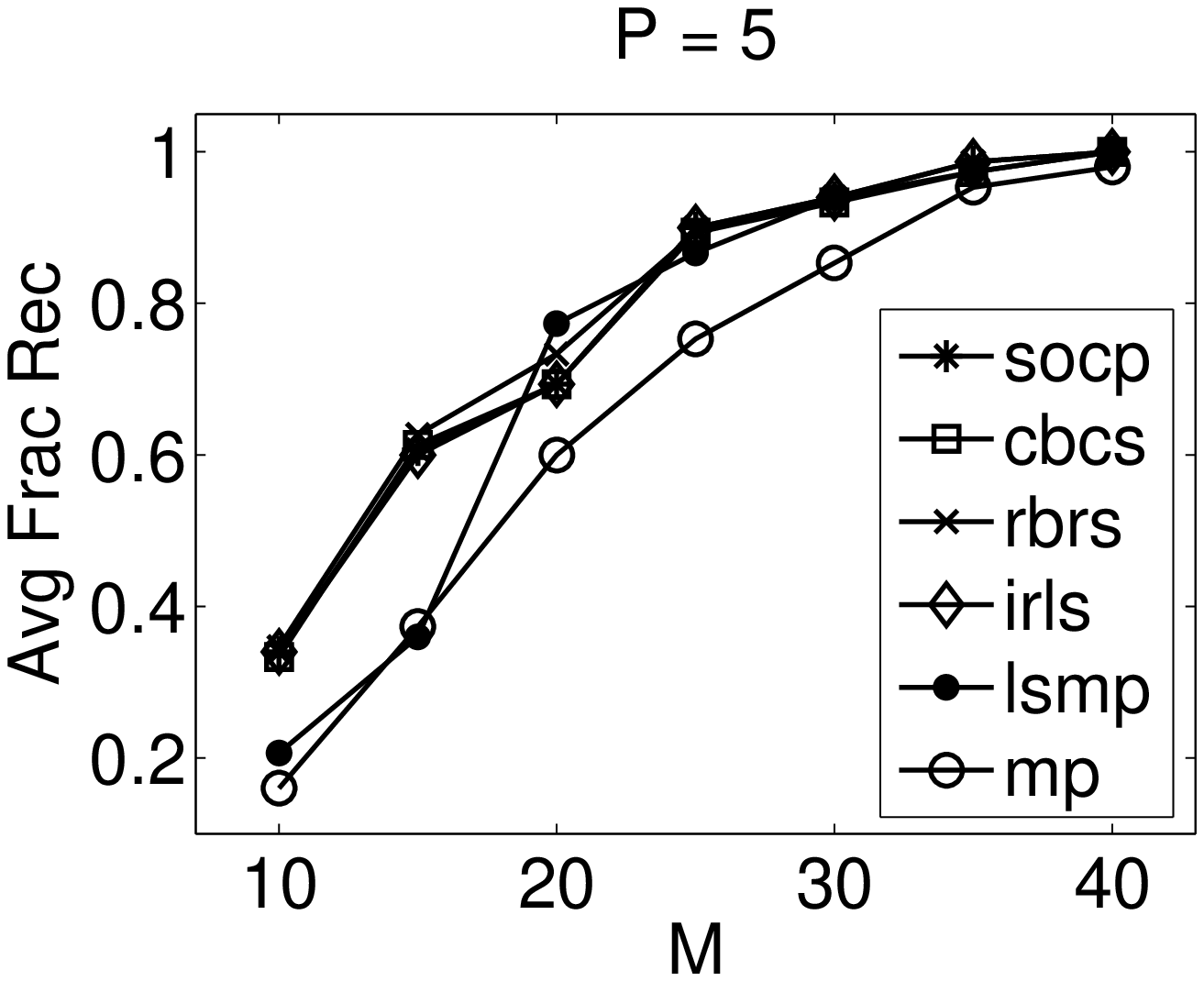,width=\widthC} &
       \epsfig{figure=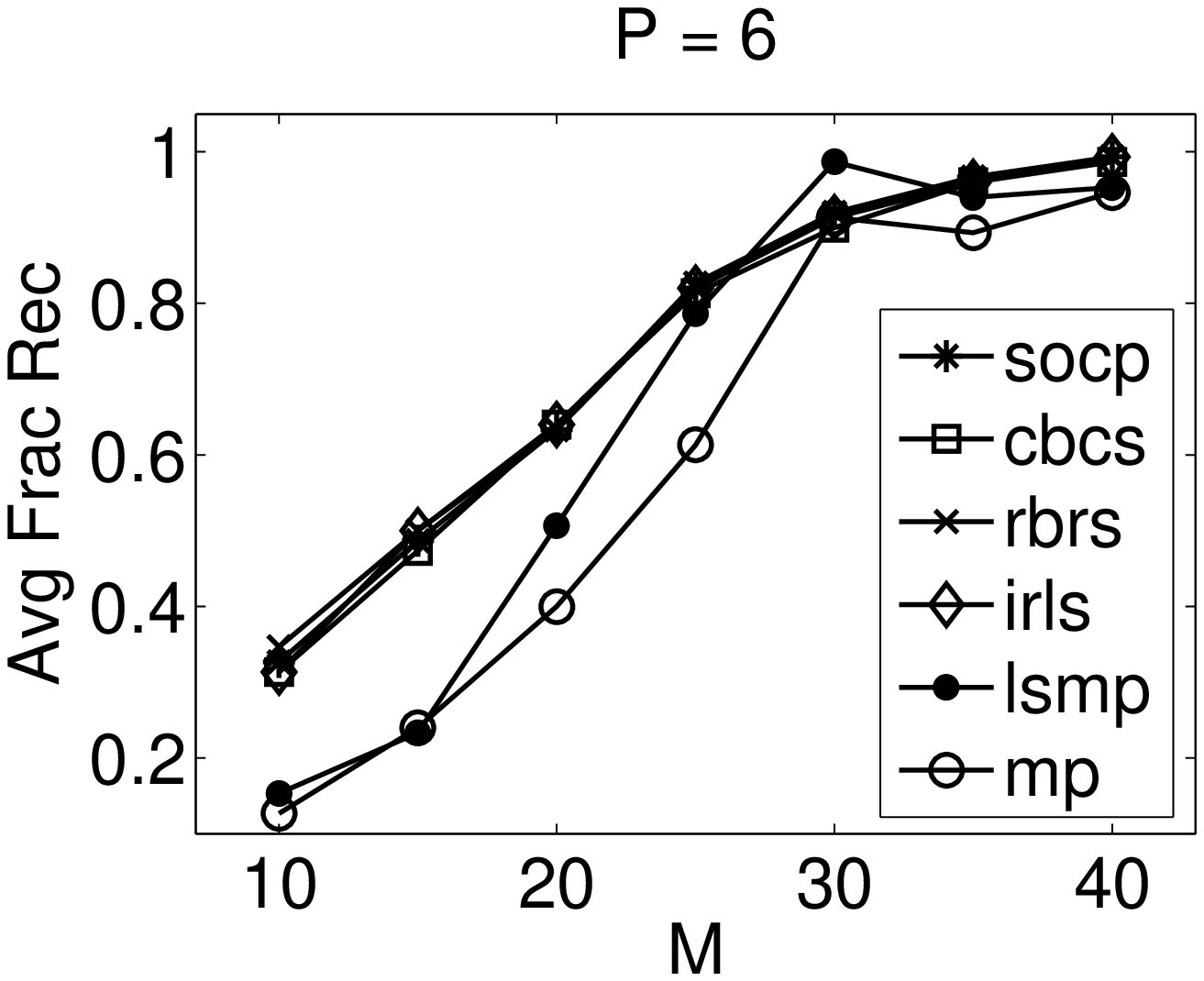,width=\widthC}
     \end{tabular}

     \begin{tabular}{cc}
       \epsfig{figure=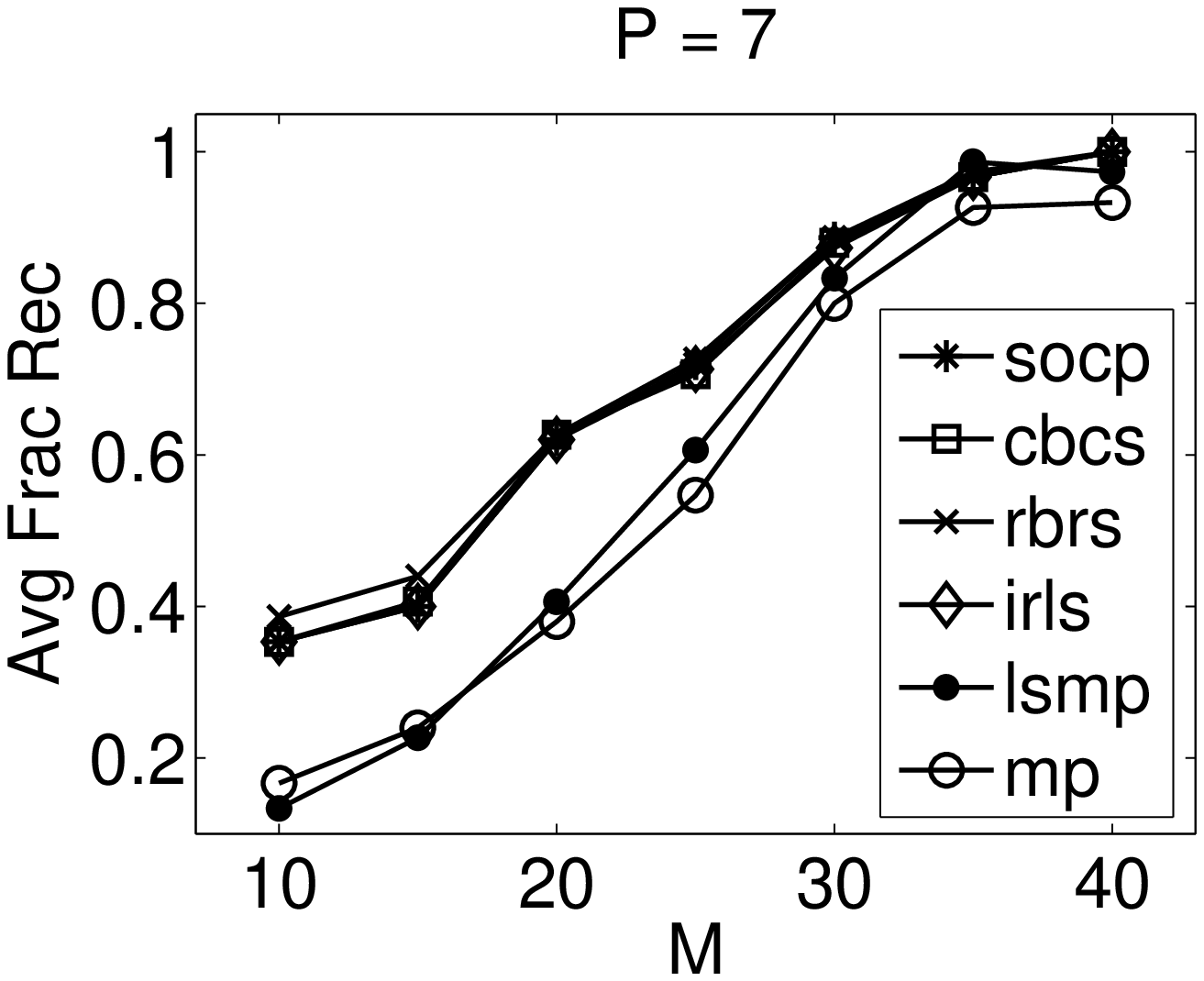,width=\widthC} &
       \epsfig{figure=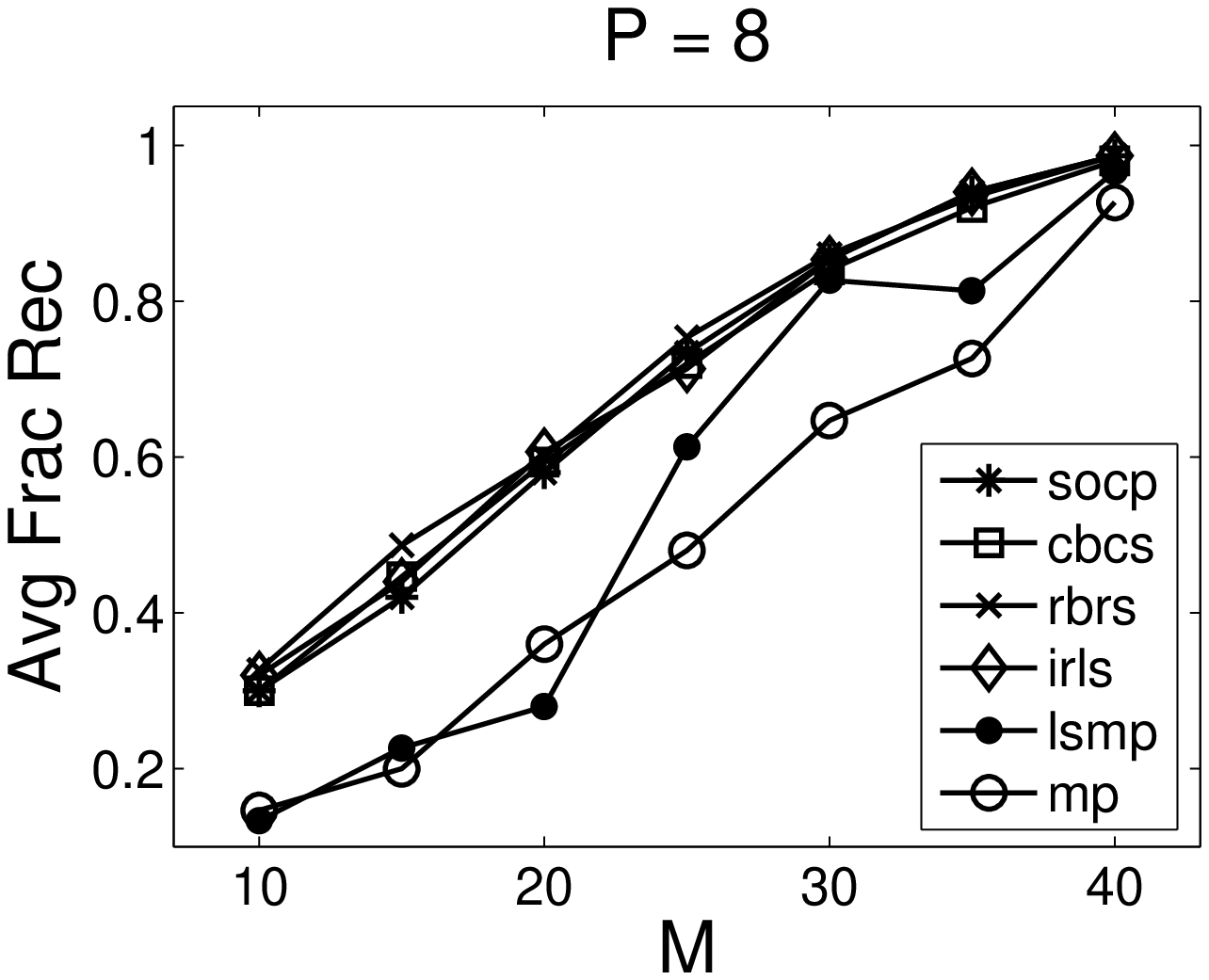,width=\widthC}
     \end{tabular}

     \caption{{\bf{Sparsity Profile Estimation in a Noiseless
     Setting}.}  Subplots depict average fraction of sparsity profile
     elements recovered over 50 trials of six algorithms as $M$ is
     varied.  $P$ is fixed per subplot, and $N=30$ and $K=3$ for all
     trials.  Data is generated as described in
     Sec.~\ref{subsubsec:e1}.  Recovery scores for IRLS, RBRS, CBCS,
     and SOCP assume a good choice of $\lambda$ is known. For large
     $M$, all algorithms exhibit high recovery rates; for large $P$,
     small $M$, or both, the algorithms that seek to minimize
     (\ref{mmv2_eq2}, \ref{mmv2_eq4}) generally outperform those that
     greedily pursue a solution.}

     \label{fig:e1_afr}
     \end{center}
   \end{figure}

Each subplot of Fig.~\ref{fig:e1_afr}
   depicts the average fraction of recovered sparsity profile elements
   versus the number of knowns, $M$, for a fixed value of $P$,
   revealing how performance varies as the $\la{F}_p \in \field{R}^{M
   \times N}$ matrices go from being underdetermined to
   overdetermined.

   {\em{Recovery Trends}.} As the number of knowns $M$ increases,
   recovery rates improve substantially, which is sensible.  For large
   $M$ and small $P$, the six algorithms behave similarly,
   consistently achieving nearly 100\% recovery.  For large $P$ and
   moderate $M$, however, sparsity profile recovery rates are
   dismal---as $P$ increases, the underlying simultaneous sparsity of
   the $\la{g}_p$s is not enough to combat the increasing number of
   unknowns, $PN$.  As $M$ is decreased and especially when $P$ is
   increased, the performance of the greedy techniques falls off
   relative to that of IRLS, RBRS, CBCS, and SOCP, showing that the
   convex relaxation approach itself is a sensible way to
   approximately solve the formal NP-Hard combinatorial MSSO
   simultaneous sparsity problem.  Furthermore, the behavior of the
   convex algorithms relative to the greedy ones coincides with the
   studies of greedy vs.~convex programming sparse approximation
   methods in single-vector \cite{Che1998, Cot1999} and SSMO contexts
   \cite{Cot2005}.  Essentially, in contrast with convex programming
   techniques, the greedy algorithms only look ahead by one term,
   cannot backtrack on sparsity profile element choices, and do not
   consider updating multiple rows of unknowns of the $\la{G}$ matrix
   at the same time.  LSMP tends to perform slightly better than MP
   because it solves a least squares minimization and explicitly
   considers earlier chosen rows whenever it seeks to choose another
   row of $\la{G}$.

   {\em{Convergence}.} Across most trials, IRLS, RBRS, CBCS, and SOCP
   converge rapidly and do not exceed the maximum limit of 500 outer
   iterations.  The exception is CBCS when $M$ is small and $P = 8$:
   here, the objective function frequently fails to decrease by less
   than the specified $\delta = 10^{-5}$.

   {\em{Runtimes}.} For several fixed $(M,P)$ pairs,
   Table~\ref{tab:r1} lists the average runtimes of each algorithm
   across the 50 trials associated with each pair.\footnote{In the
   interest of space we do not list average runtimes for all fifty-six
   $(M,P)$ pairs.}  For IRLS, RBRS, CBCS, and SOCP, runtimes are also
   averaged over the many $\lambda$ runs.  Among the convex
   minimization methods, SOCP seems superior given its fast runtimes
   in three out of four cases.  Peak memory usage is not tracked here
   because it is difficult to do so when using MATLAB for such small
   problems; it will be tracked during the third experiment where the
   system matrices are vastly larger and differences in memory usage
   across the six algorithms are readily apparent.
   \begin{table}
   \begin{center}
   \small
   \begin{tabular}{|l|r|r|r|r|}
   \hline
                          & \multicolumn{4}{c|}{$(M,P)$} \\ \hline
       {\bf{Algorithm}}     & (10,8) & (20,1) & (30,5) & (40,8) \\ \hline
    MP                    & 5.4      & 1.8     & 2.6  & 4.0  \\
    LSMP                  & 11.4     & 5.6     & 15.6 & 27.6 \\
    IRLS  & 92.6     & 10.1    & 73.2  & 175.0 \\
    RBRS  & 635.7    & 36.0    & 236.8 & 401.6 \\
    CBCS  & 609.8    & 7.1     & 191.4 & 396.3 \\
    SOCP  & 44.3     & 37.0    & 64.3  & 106.5 \\ \hline
   \end{tabular}
   \end{center}

   \caption{{\bf{Average Algorithm Runtimes for Noiseless Sparsity
   Profile Estimation.}} For several fixed $(M,P)$ pairs, each
   algorithm's average runtime over the corresponding 50 trials is
   given in units of milliseconds; $N=30$ and $K=3$ for all trials
   (runtimes of the latter four algorithms are also averaged over the
   multiple $\lambda$ runs per trial).  MP is substantially faster
   than the other techniques, as expected.  For larger problems,
   e.g. $(M,P)=(10,8)$, the runtimes of both RBRS and CBCS are
   excessive relative to those of the other convex minimization
   techniques, IRLS and SOCP\@.}  \label{tab:r1}

   \end{table}

\begin{figure}
      \begin{center}
      \small
      \begin{tabular}{cc}
        \epsfig{figure=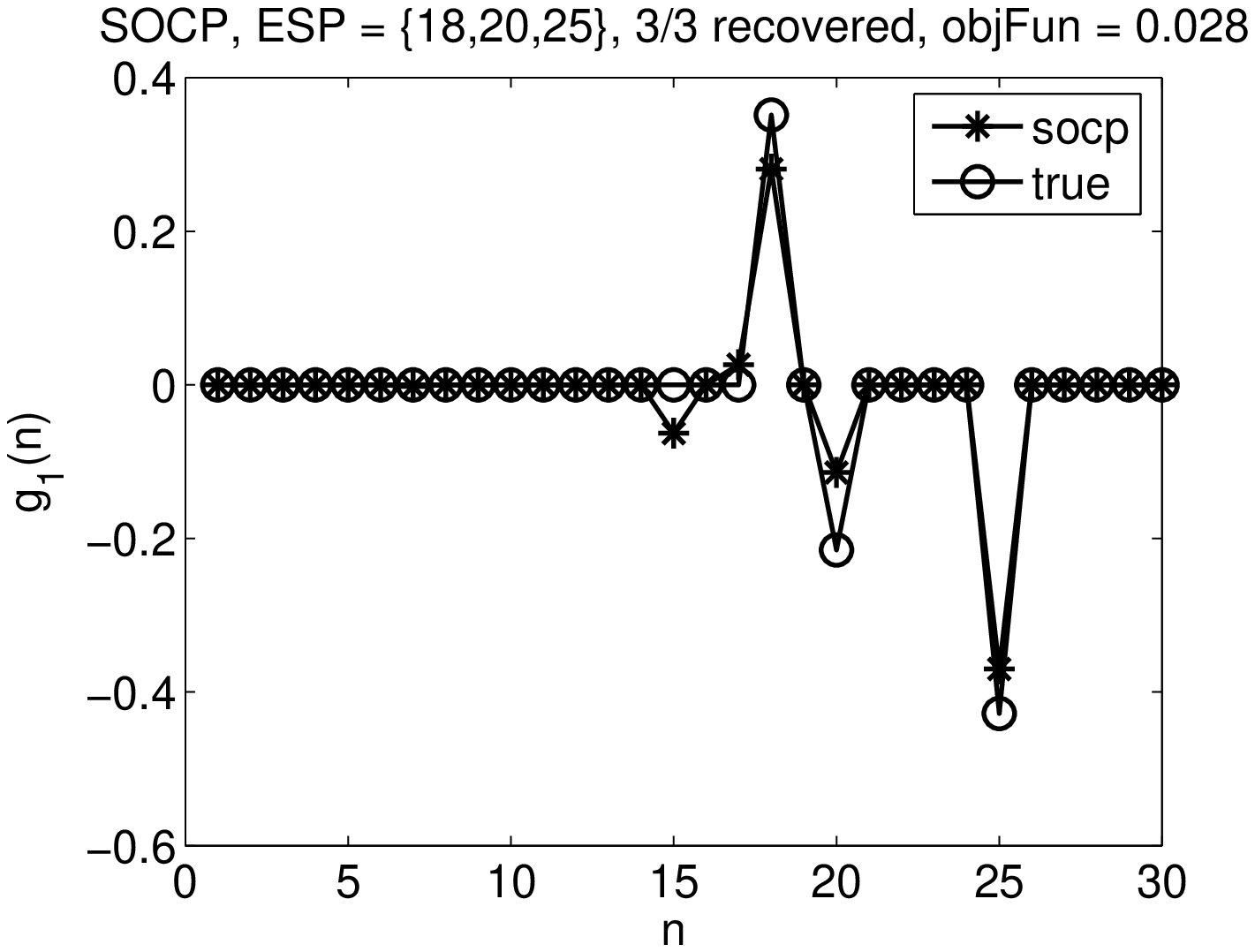,width=2in} &
        \epsfig{figure=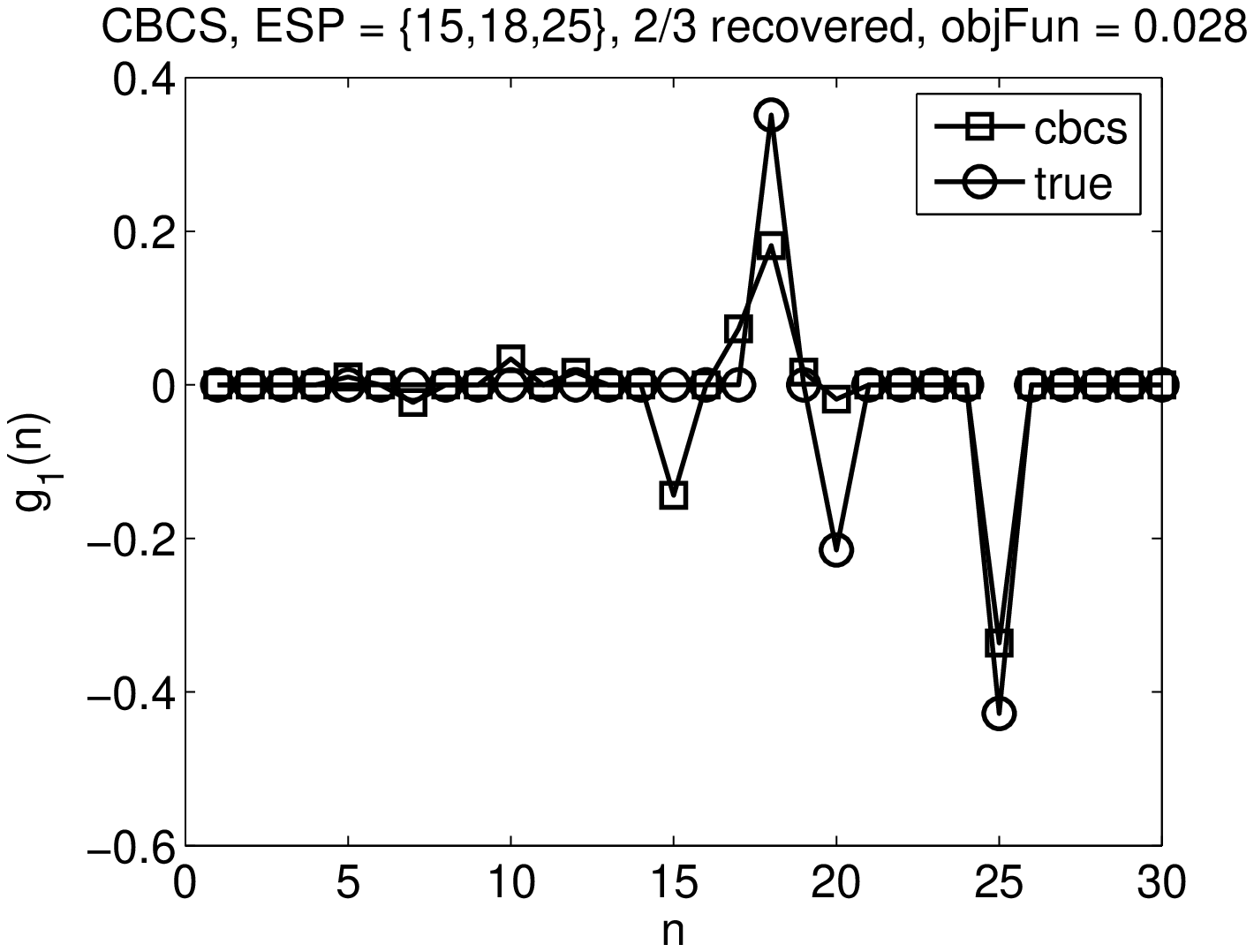,width=2in}
      \end{tabular}

      \begin{tabular}{cc}
        \epsfig{figure=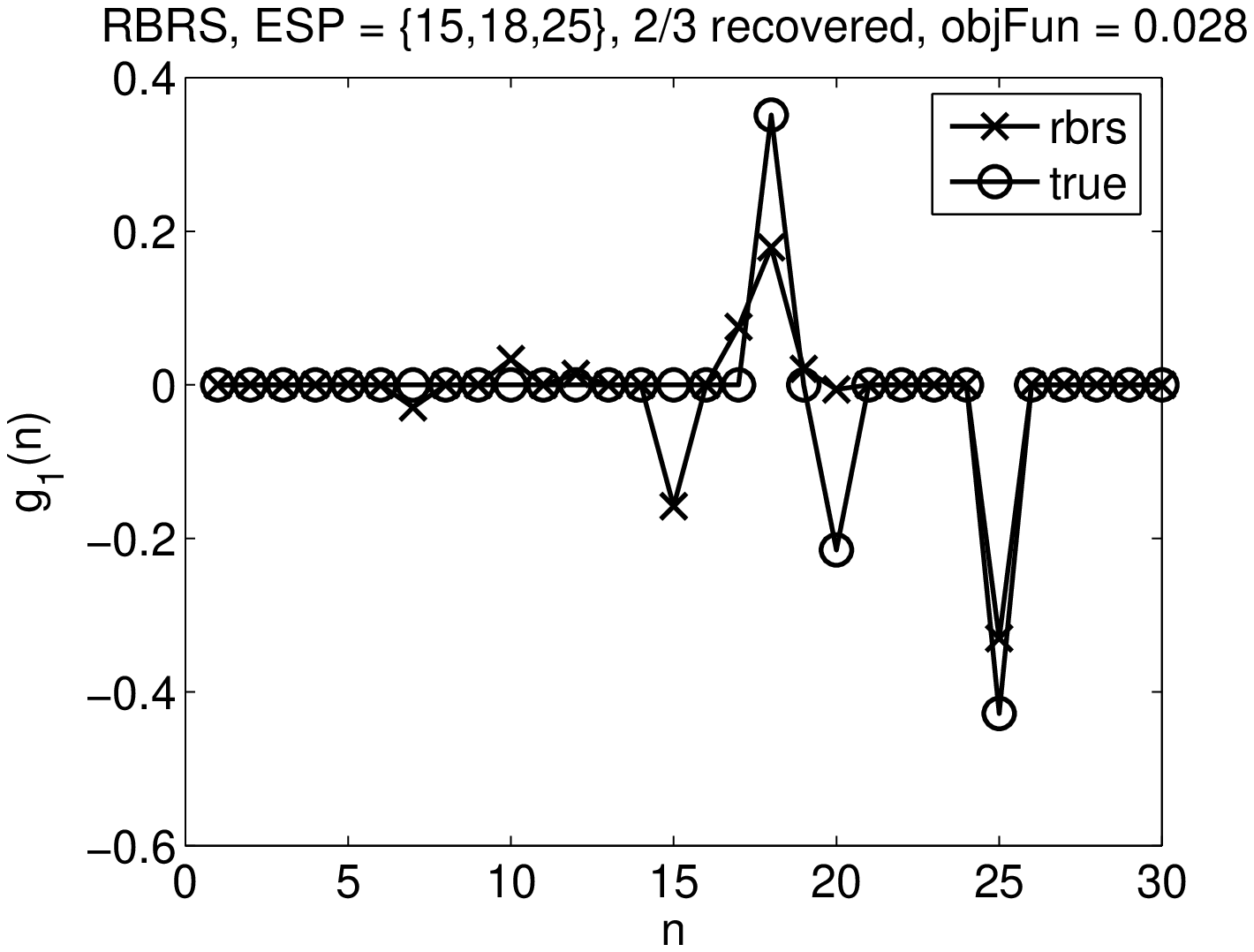,width=2in} &
        \epsfig{figure=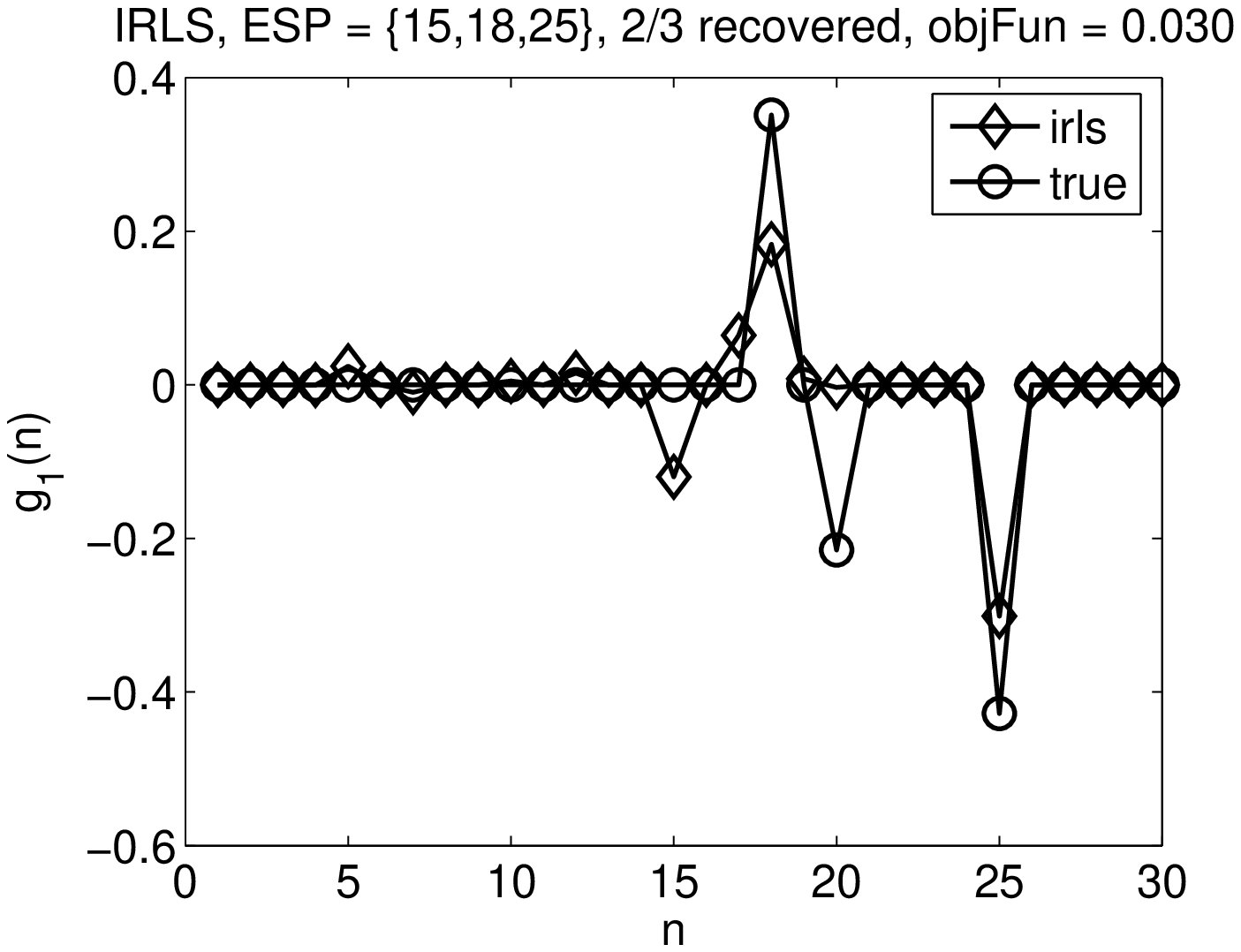,width=2in}
      \end{tabular}

      \caption{{\bf{Noiseless Sparsity Profile Estimation with IRLS,
      RBRS, CBCS, SOCP}.}  Here $M=10, N=30, P=1$, and $K=3$.  The
      algorithms are applied with $\lambda$ fixed at 0.03 and attempt
      to estimate the single unknown vector, $\la{g}_1$, along with
      the sparsity profile.  Subplots depict the elements of both the
      estimated and actual $\la{g}_1$, along with the estimated
      sparsity profile (ESP), number of profile terms recovered, and
      objective function value.  SOCP leads to a superior sparsity
      profile estimate, and SOCP, RBRS, and CBCS seem to minimize the
      convex objective given in (\ref{mmv2_eq2}, \ref{mmv2_eq4}).
      IRLS does not, but still manages to properly identify 2 out of 3 sparsity
      profile terms.}

      \label{fig:e1_zoom}
      \end{center}
   \end{figure}

   {\em{Closer Look: Solution Vectors}.}  We now observe how the
   algorithms that seek to minimize the convex objective behave during
   the 43rd trial when $K = 3$, $N = 30$, $M = 10$, and $P
   = 1$, corresponding to the base case problem of estimating one
   sparse real-valued vector, $\la{g}_1$.  Fig.~\ref{fig:e1_zoom}
   illustrates estimates obtained by SOCP, CBCS, RBRS, and IRLS when
   $\lambda = 0.03$; for each algorithm, a subplot shows elements of
   both the estimated and actual $\la{g}_1$, and lists the estimated
   sparsity profile (ESP), number of profile terms recovered, and
   value of the objective function given in (\ref{mmv2_eq2},
   \ref{mmv2_eq4}).  Although RBRS, CBCS, and SOCP yield slightly
   different solutions (among which SOCP yields the best profile
   estimate), they all yield an objective function equal to
   $0.028\pm10^{-5}$.  Convex combinations of the three solutions
   continue to yield the same value, suggesting that the three
   algorithms have found solutions among a convex set that is the
   global solution to the objective posed in (\ref{mmv2_eq2},
   \ref{mmv2_eq4}).  Given the fact that in this case SOCP outperforms
   RBRS and CBCS, we see that even the globally optimal solution to
   the relaxed convex objective does not necessarily optimally solve
   the true $K$-sparse profile recovery problem.  In contrast to the
   other methods, IRLS yields a slightly higher objective function
   value, 0.030, and its solution vector is not part of the convex
   set---it does however correctly determine 2 of the 3 terms of the true
   sparsity profile.

\begin{figure}\centering
      \epsfig{figure=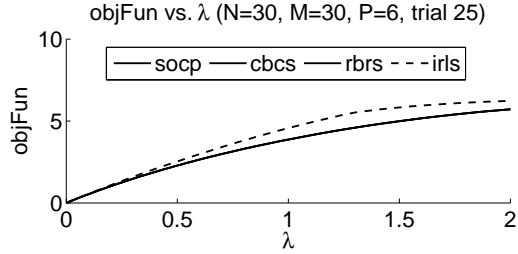,width=3.0in}

      \caption{{\bf{Noiseless Sparsity Profile Estimation: Objective
      Function Behavior}.} For the $25^{\sst{th}}$ trial of the
      $(M,P)=(30,6)$ series, SOCP, CBCS, RBRS, and IRLS are used to
      solve (\ref{mmv2_eq2}, \ref{mmv2_eq4}) for 70 values of $\lambda
      \in [0, 2]$; the value of the objective function vs.~$\lambda$
      is given above.  For $\lambda > \frac{1}{4}$, IRLS's solutions
      do not minimize the objective as well as those produced by the
      three other methods.}

      \label{fig:e1_objfun}
   \end{figure}

   {\em{Closer Look: Objective Function Behavior}.}  Concluding the
   experiment, Fig.~\ref{fig:e1_objfun} plots the objective
   vs.~$\lambda$ for the 25th trial when $M=30$ and $P=6$,
   studying how the objective (\ref{mmv2_eq2}, \ref{mmv2_eq4}) varies
   with $\lambda$ when applying SOCP, CBCS, RBRS, and IRLS\@.  For all
   seventy values of $\lambda \in [0, 2]$, SOCP, CBCS, and RBRS
   generate solutions that yield the same objective function value.
   For $\lambda < \frac{1}{4}$, IRLS attains the same objective
   function value as the other methods, but as $\lambda$ increases,
   IRLS is unable to minimize the objective function as well as SOCP,
   RBRS, and CBCS\@.  The behavior in Fig.~\ref{fig:e1_objfun} occurs
   consistently across the fifty trials of the other $(M, P)$ pairs.

\subsection{Sparsity Profile Estimation in the Presence of Noise}
\label{subsec:e2}

   \subsubsection{Overview} We now evaluate how well the algorithms of
   Sec.~\ref{sec:algorithms} estimate sparsity profiles when the
   underlying $\la{g}_p$s are each strictly and simultaneously
   $K$-sparse and the observation $\la{d}$ of (\ref{mmv2_eq1}) is
   corrupted by additive white Gaussian noise.  The
   signal-to-noise ratio (SNR) and $K$ are varied across sets of Monte
   Carlo trials in order to gauge algorithm performance across many
   scenarios.  For a given trial with a fixed SNR level in units of
   decibels (dB), the $M$ elements of the true observation vector,
   $\la{d}_{\sst{true}}$, are corrupted with independent and
   identically distributed (i.i.d.)  zero-mean Gaussian noise with
   variance $\sigma^2$, related to the SNR as follows:
   \begin{equation}\label{snr} 
      \sigma^2 = \frac{1}{M} \Vert \la{d}_{\sst{true}} \Vert_2^2 
                             \cdot 10^{-\sst{SNR}/10}
   \end{equation}
   This noise measure is analogous to that of \cite{Cot2005}.
    
   \subsubsection{Details}\label{subsubsec:e2}
    We fix $N = 30$, $M = 25$, and $P = 3$, and consider $\mbox{SNR}
   \in \{-10, -5, 0, \ldots, 25, 30\}$ and $K \in \{1, 3, 5, 7, 9\}$.
   For each fixed $(\mbox{SNR}, K)$ pair, we generate 100 noisy
   observations and apply the algorithms as follows:
   \begin{adamItemize2}

     \item Generate the sparsity profile, $\la{g}_p$s, $\la{F}_p$s,
     $\la{h}_n$s, and $\la{C}_n$s as in Sec.~\ref{subsubsec:e1}.  The
     $\la{g}_p$s are simultaneously $K$-sparse and all terms are
     real-valued.

     \item Compute $\la{d}_{\sst{true}} = \la{F}_1 \la{g}_1 + \cdots +
     \la{F}_P \la{g}_P$.

     \item Construct $\la{d}_{\sst{noisy}} = \la{d}_{\sst{true}} +
     \la{n}$ where $\la{n} \sim \mathcal{N}(0, \sigma^2 \la{I})$ and
     $\sigma^2$ is given by (\ref{snr}).

     \item Apply the algorithms by providing them with
     $\la{d}_{\sst{noisy}}$ and the system matrices:

        \begin{adamItemize2}

        \item[$\circ$] MP, LSMP: iterate until $K$ elements are
        chosen or the residual approximation is $\la{0}$.  If less
        than $K$ terms are chosen, this hurts the recovery score.

        \item[$\circ$] IRLS, RBRS, CBCS, SOCP: using a pre-determined
        {\em{fixed}} $\lambda$ (see below), apply each algorithm to
        obtain estimates of the unknown vectors and sparsity profiles.

        \end{adamItemize2}

   \end{adamItemize2} After performing the above steps, we average the
   results of the 100 trials associated with each fixed $(\mbox{SNR},
   K,\mbox{alg})$ triplet to yield the average fraction of sparsity
   profile elements that each algorithm recovers.

   {\em{Control Parameter Selection}.} The $\lambda$ mentioned in the
   list above is determined as follows: before running the overall
   experiment, we generate three noisy observations for each
   $(\mbox{SNR}, K)$ pair.  We then apply IRLS, RBRS, CBCS, and SOCP,
   tuning the control parameter $\lambda$ by hand until finding a
   single value that produces reasonable solutions.  All algorithms
   then use this hand-tuned, fixed $\lambda$ and are applied to the
   other 100 noisy observations associated with the $(\mbox{SNR},K)$
   pair under consideration.  Thus, in distinct contrast to the
   noiseless experiment, we no longer assume an ideal $\lambda$ is
   known for each denoising trial.

   \subsubsection{Results}\label{subsubsec:e2_res} 
\begin{figure}
   \begin{center}
   \small
     \begin{tabular}{ccc}
       \epsfig{figure=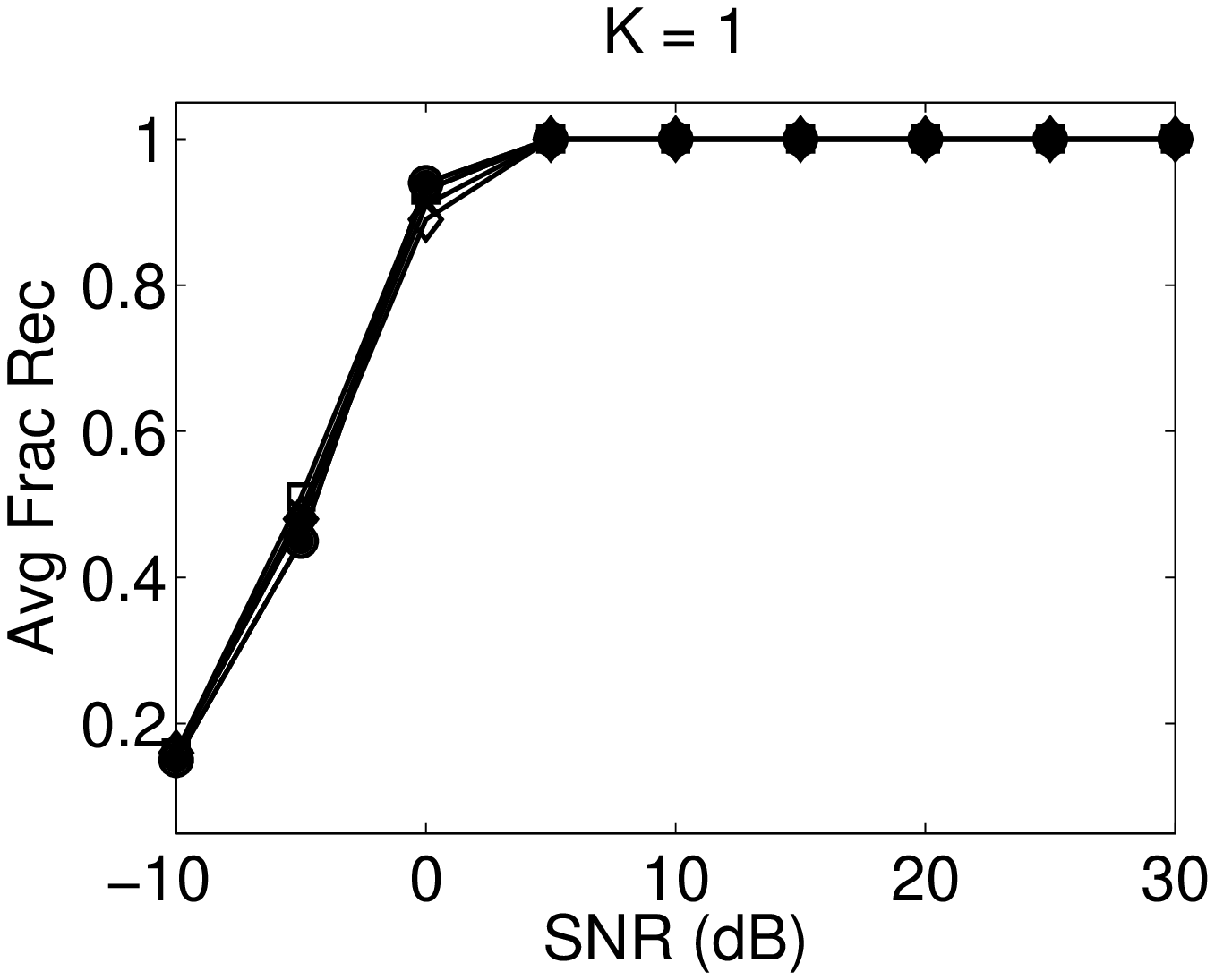,width=\widthC} &
       \epsfig{figure=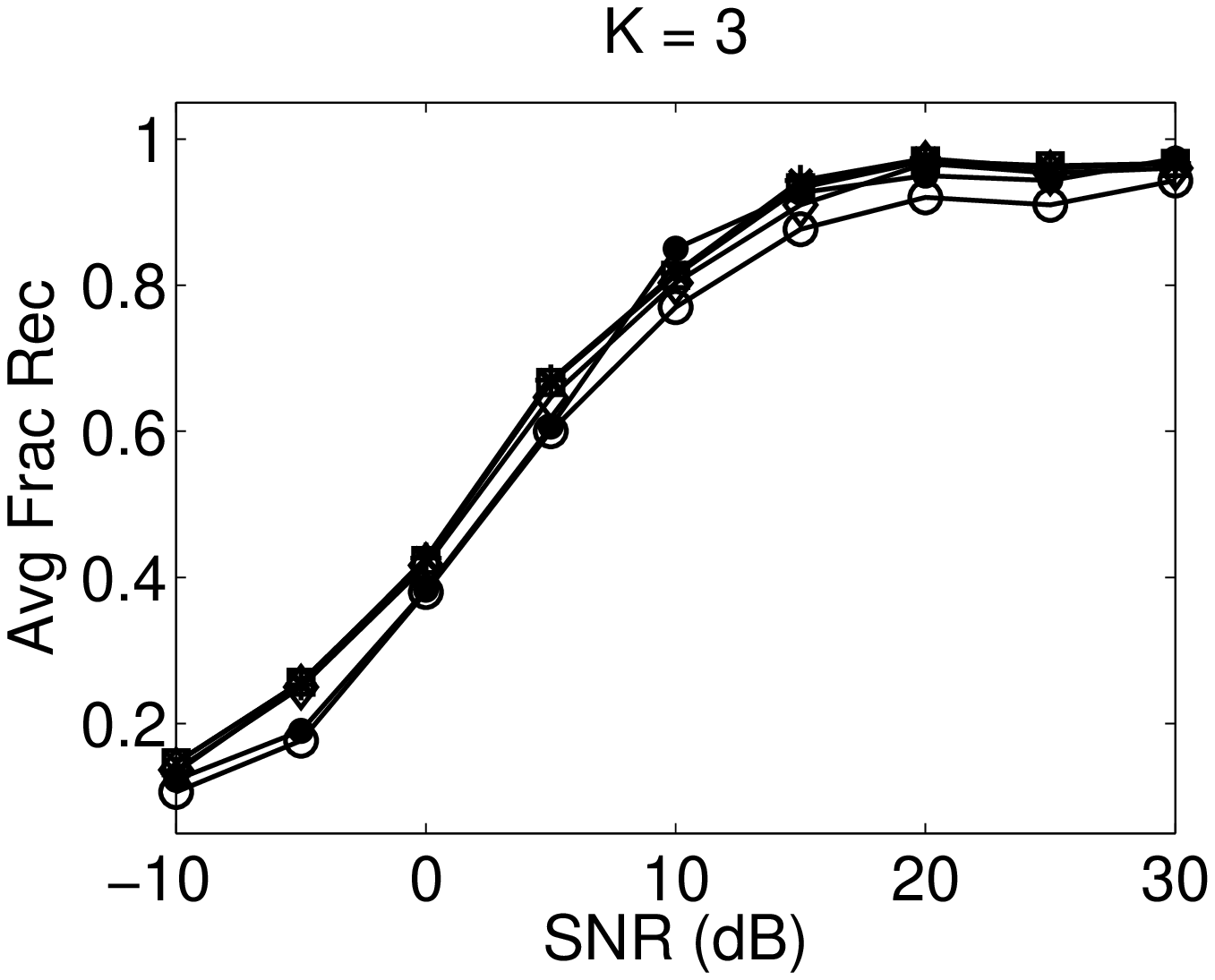,width=\widthC} &
       \epsfig{figure=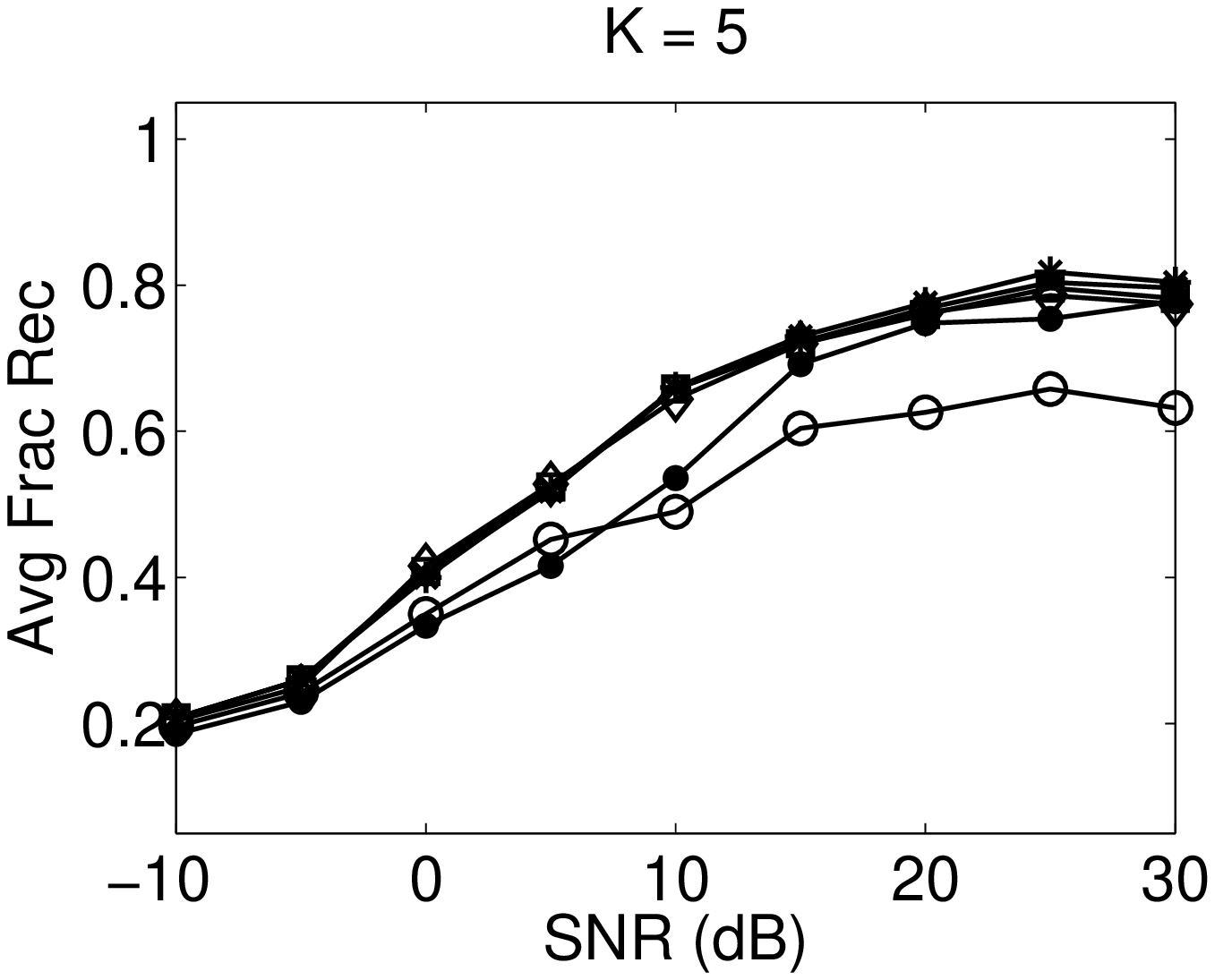,width=\widthC}
     \end{tabular}

     \begin{tabular}{c}
       \epsfig{figure=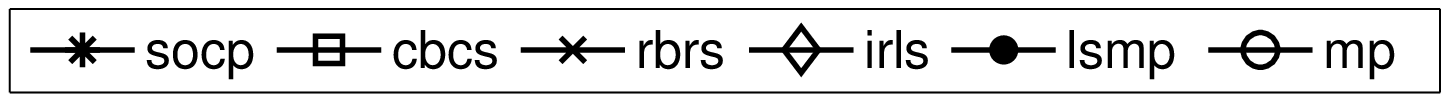,width=2.5in}
     \end{tabular}

     \begin{tabular}{cc}
       \epsfig{figure=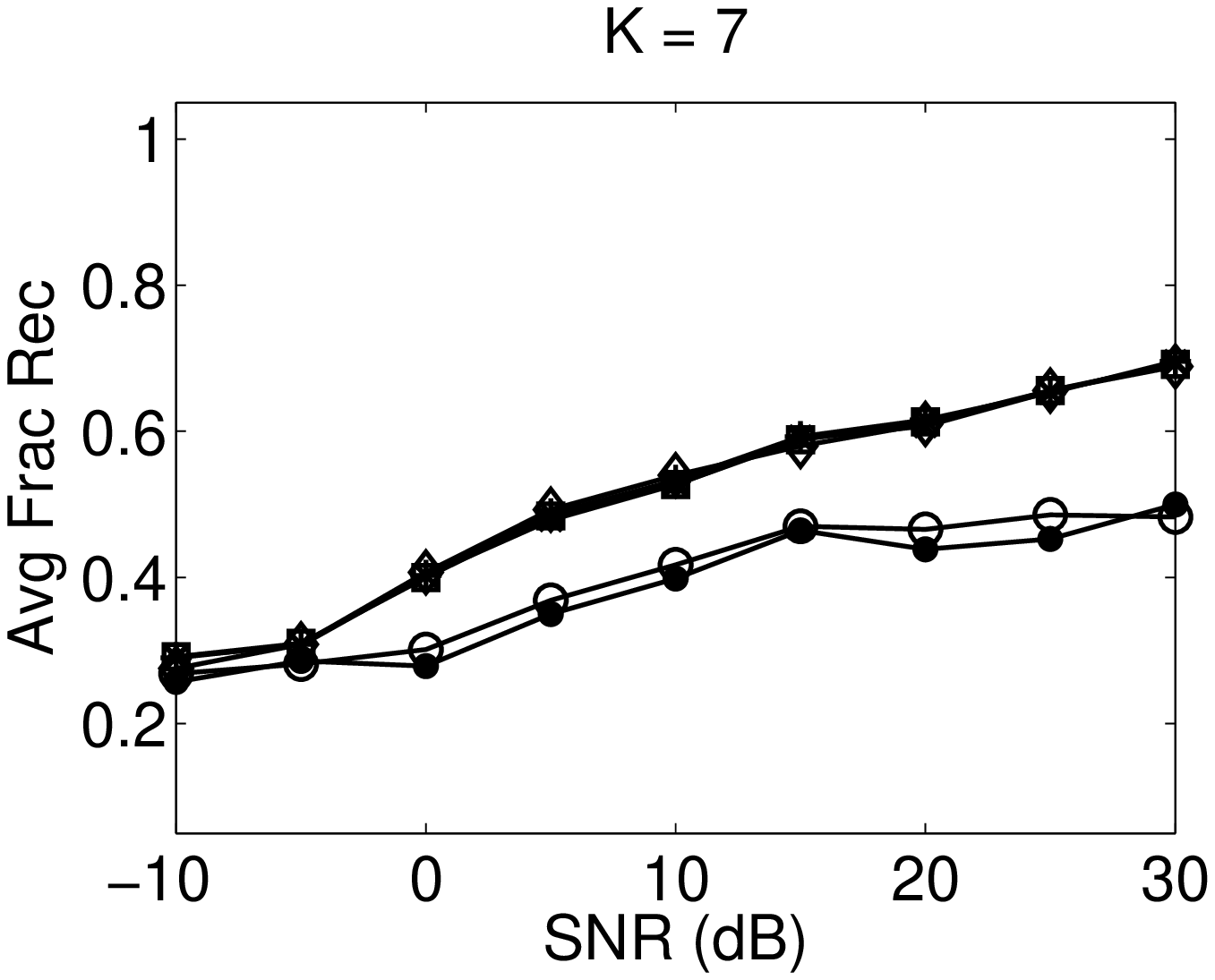,width=\widthC} &
       \epsfig{figure=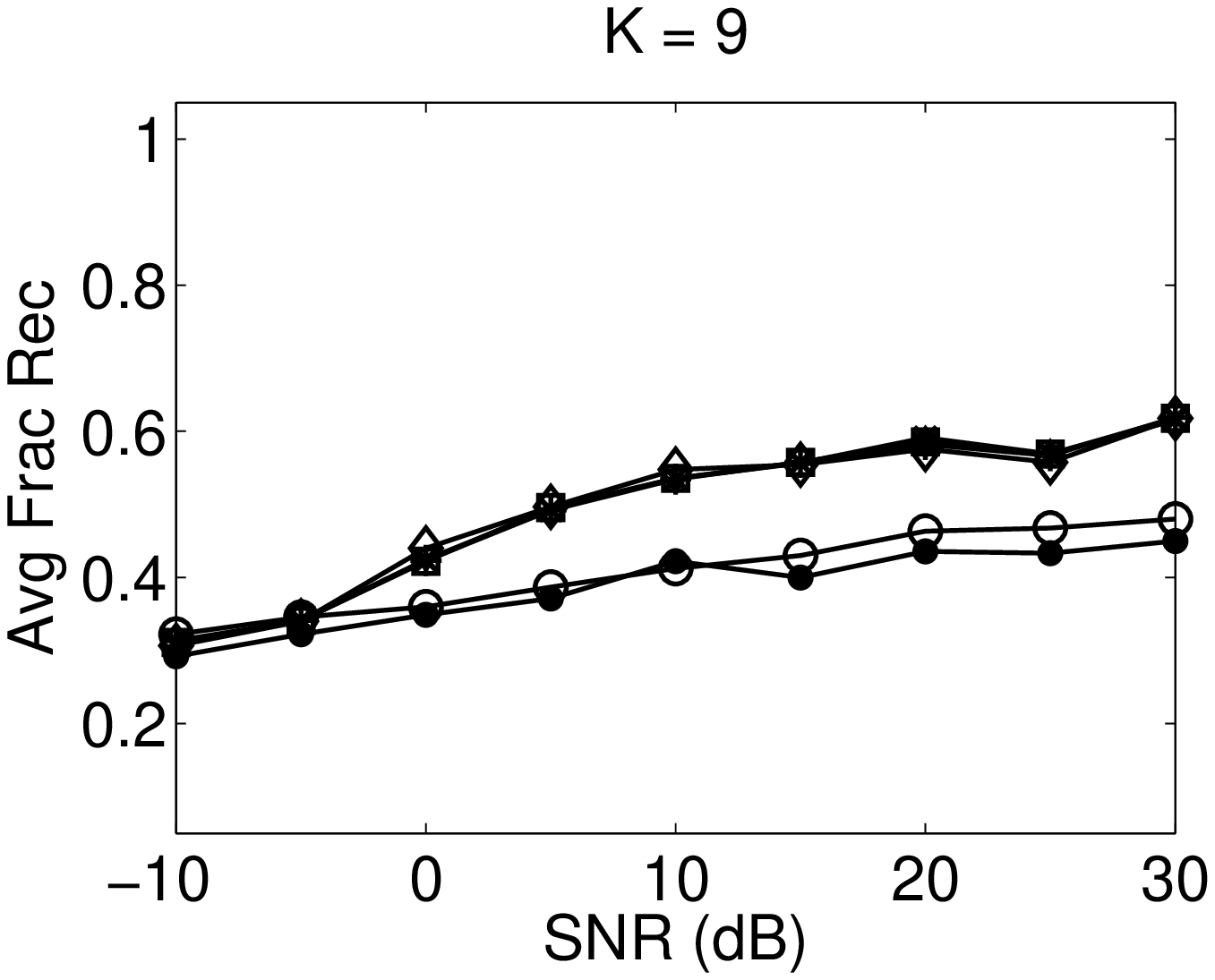,width=\widthC}
     \end{tabular}

     \caption{{\bf{Sparsity Profile Estimation in the Presence of
     Noise}.}  Each subplot depicts the average fraction of recovered
     sparsity profile elements versus SNR for a fixed $K$, revealing
     how well the algorithms recover the $K$ elements of the sparsity
     profile amidst noise in the observation.  Each data point is the
     average fraction recovered across 100 trials; data is randomly
     generated as described in Sec.~\ref{subsubsec:e2}.  $N, M$ and
     $P$ are always fixed at 30, 25, and 3, respectively. For each
     $(\mbox{SNR}, K)$ pair, a ``good'' lambda is chosen by denoising
     a few cases by hand and then using this fixed $\lambda$ for 100
     fresh denoising trials.  Performance degrades with increasing $K$
     and decreasing SNR\@.  For large $K$, the greedy algorithms perform
     worse than IRLS, SOCP, RBRS, and CBCS, whereas the latter four
     methods perform essentially identically across all
     $(\mbox{SNR},K)$ combinations.}

     \label{fig:e2_K} 
     \end{center}
   \end{figure}

Each subplot of
   Fig.~\ref{fig:e2_K} depicts the average fraction of recovered
   sparsity profile elements versus SNR for a fixed $K$, revealing how
   well the six algorithms are able to recover the $K$ elements of the
   sparsity profile amidst noise in the observation.  Each data point
   is the average fraction recovered across 100 trials.

   {\em{Recovery Trends}.} When $K=1$, we see from the upper-left
   subplot of Fig.~\ref{fig:e2_K} that all algorithms have essentially
   equal performance for each SNR\@.  Recovery rates improve
   substantially with increasing SNR, which is sensible.  For each
   algorithm, we see across the subplots that performance generally
   decreases with increasing $K$; in other words, estimating a large
   number of sparsity profile terms is more difficult than estimating
   a small number of terms.  This trend is evident even at high SNRs.
   For example, when SNR is 30 dB and $K=7$, SOCP is only able to
   recover $\sim 70\%$ of sparsity profile terms.  When $K=9$, the
   recovery rate falls to $\sim 60\%$.  For low SNRs, e.g., -5 dB, all
   algorithms tend to perform similarly, but the greedy algorithms
   perform increasingly worse than the others as $K$ goes from
   moderate-to-large and SNR surpasses zero dB\@.  In general, MP
   performs worse than LSMP, and LSMP in turn performs worse than
   IRLS, SOCP, RBRS, and CBCS, while the latter four methods exhibit
   essentially the same performance across all SNRs and $K$s.  For
   $K=3$, MP's performance falls off relative to IRLS, SOCP, RBRS,
   and CBCS, but LSMP's does not.  As $K$ transitions from 3 to 5,
   however, LSMP performs as badly as MP at low SNRs, but its
   performance picks up as SNR increases.  As $K$ continues to
   increase beyond 5, LSMP's performance is unable to surpass that of
   MP, even when SNR is large.  Overall, Fig.~\ref{fig:e2_K} shows
   that convex programming algorithms are superior to greedy methods
   when estimating sparsity profiles in noisy situations; this
   coincides with data collected in the noiseless experiment in
   Sec.~\ref{subsec:e1}, as well as the empirical findings of
   \cite{Cot1999, Cot2005}.

   {\em{Convergence}.} For many denoising trials, CBCS typically
   requires more iterations than the other techniques in order to
   converge.  At times, it fails to converge to within the specified
   $\delta = 10^{-5}$, similarly to how it behaves during the
   noiseless experiment of Sec.~\ref{subsec:e1}.

   {\em{Runtimes}.} Across all denoising trials, MP, LSMP, IRLS, RBRS,
   CBCS, SOCP have average runtimes of 3.1, 25.1, 57.2, 247.0, 148.5,
   and 49.2 milliseconds.  It seems SOCP is best for denoising given
   that it is the fastest algorithm among the four methods that
   outperforms the greedy ones.  IRLS is nearly as fast as SOCP and
   thus is a close second choice for sparsity profile estimation.

\begin{figure}
   \begin{center}
   \small
     \begin{tabular}{cc}
       \epsfig{figure=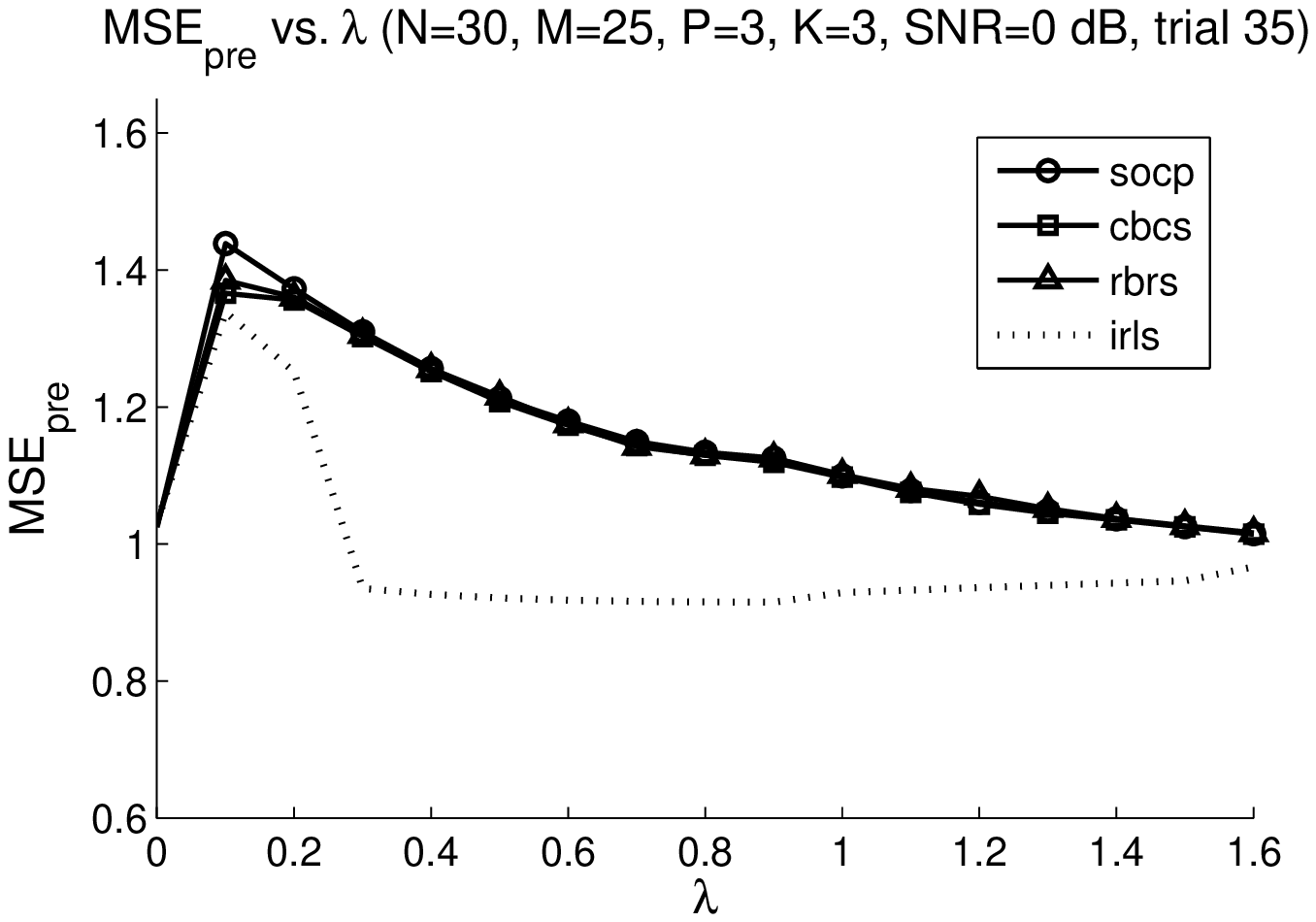,width=2.5in} &
       \epsfig{figure=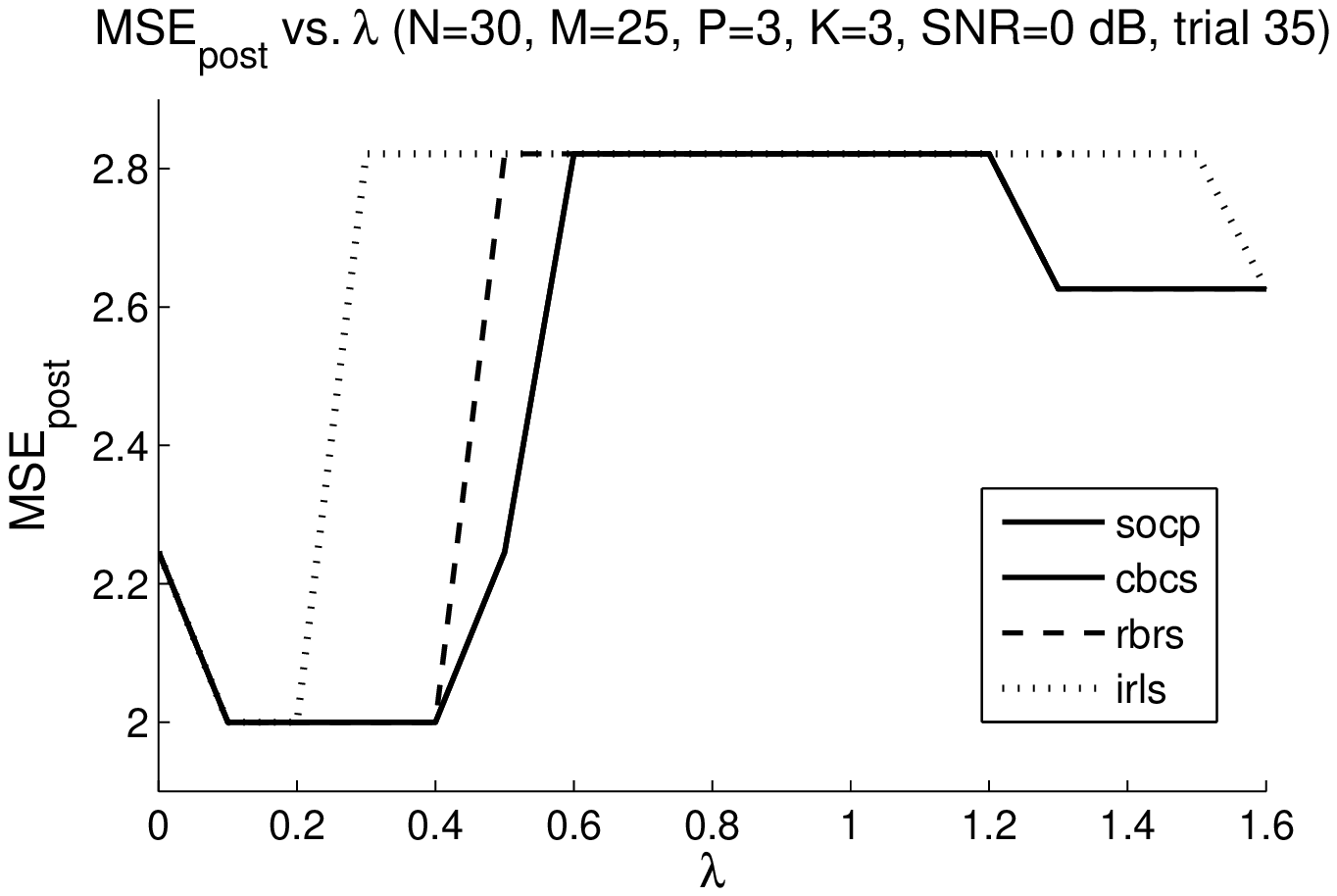,width=2.5in}
     \end{tabular}

     \caption{{\bf{MSEs of Convex Minimization Methods Before and
     After Estimating the Sparsity Profile and Retuning the
     Solution}.}  Here MSE vs.~$\lambda$ is studied during the
     $35^{\sst{th}}$ trial of the $(\mbox{SNR}=0\mbox{ dB}, K=3)$
     denoising series.  Fixing $\lambda$ and applying a given
     algorithm yields the solution
     $\widehat{\la{h}}(\mbox{alg},\lambda)$.  Left subplot: MSEs of
     the $\widehat{\la{h}}(\mbox{alg},\lambda)$s vs.~the true solution
     $\la{h}_{\sst{tot}}$.  Right subplot: MSEs of the solution
     estimates after they undergo retuning to be strictly and
     simultaneously $K$-sparse.  (Sec.~\ref{subsubsec:e2_res} outlines
     the retuning process.) For all algorithms and $\lambda$s, MSE
     increases substantially relative to the left subplot.  None of
     the methods correctly estimate the true $K$-term sparsity profile
     and thus the retuning step causes the estimated solution to
     branch further away (in the MSE sense) from the actual one.}

     \label{fig:e2_zoom} 
   \end{center}
   \end{figure}

   {\em{Closer Look: Mean Square Errors of Convex Minimization Methods
   Before and After Estimating the Sparsity Profile and Retuning the
   Solution}.}  Now let us consider the $35$th trial of the
   $(\mbox{SNR} = 0 \mbox{ dB}, K = 3)$ pair.  We do away with the
   fixed $\lambda$ assumption and now assume we care (to some extent)
   not only about estimating the sparsity profile, but the true
   solution $\la{h}_{\sst{tot}}$ as well.  To proxy for this, we study
   how the mean square errors (MSEs) of solutions generated by IRLS,
   SOCP, RBRS, and CBCS behave across $\lambda$ before and after
   identifying the sparsity profile and retuning the solution.
   Figure~\ref{fig:e2_zoom} depicts the results of this investigation.

   Running each algorithm for a particular $\lambda$ yields a solution
   $\widehat{\la{h}}_{\sst{tot}}(\mbox{alg},\lambda)$.  The left
   subplot simply illustrates the MSEs of the
   $\widehat{\la{h}}_{\sst{tot}}(\mbox{alg},\lambda)$s with respect
   to the true solution. Among SOCP, RBRS, CBCS, and IRLS, only the
   last is able to determine solutions with MSEs less than unity
   (consider the IRLS error curve for $\lambda \geq 0.3$).

   Consider now retuning each of the
   $\widehat{\la{h}}_{\sst{tot}}(\mbox{alg}, \lambda)$s as follows:
   unstack each into $\widehat{\la{h}}_n(\mbox{alg},\lambda)$ for $n
   \in \{1, \ldots, N\}$ and then remember the $K$ vectors whose
   $\ell_2$ energies are largest, yielding an estimate of the
   $K$-element sparsity profile.  Let these estimated indices be
   $\{q_1, \ldots, q_K\}$.  Now, generate a {\em{retuned solution}} by
   using the $K$ matrices associated with the estimated sparsity
   profile and solving $\la{d}_{\sst{noisy}} = [\la{C}_{q_1}, \ldots,
   \la{C}_{q_K} ] \la{x}_{\sst{tot}}$ for $\la{x}_{\sst{tot}} \in
   \field{R}^{KP}$.  This latter vector consists of $KP$ elements and
   by unstacking it we obtain a retuned estimate of the
   $\widehat{\la{h}}_n(\mbox{alg},\lambda)$s, e.g.,
   $\widehat{\la{h}}_{q_1}(\mbox{alg},\lambda)$ equals the first $K$
   elements of $\la{x}_{\sst{tot}}$, and so forth, while the other
   $\widehat{\la{h}}_n(\mbox{alg},\lambda)$s for $n \notin \{q_1,
   \ldots, q_K\}$ are now simply all-zeros vectors.  Reshuffling the
   retuned $\widehat{\la{h}}_n(\mbox{alg},\lambda)$s yields
   $\widehat{\la{g}}_p(\mbox{alg},\lambda)$s that are strictly and
   simultaneously $K$ sparse whose weightings yield the best match to
   the noisy observation in the $\ell_2$ sense.  Unlike the original
   solution estimate, which is not necessarily simultaneously
   $K$-sparse, here we have enforced true simultaneous $K$-sparsity.
   We may or may not have improved the MSE with respect to the true
   solution: for example, if we have grossly miscalculated the
   sparsity profile, the MSE of the retuned solution is likely to
   increase substantially, but if we have estimated the true sparsity
   profile exactly, then the retuned solution will likely be quite
   close (in the $\ell_2$ sense) to the true solution, and MSE will
   thus decrease.
   
   The MSEs of these {\em{retuned solutions}} with respect to the true
   $\la{h}_{\sst{tot}}$ are plotted in the right subplot of
   Fig.~\ref{fig:e2_zoom}.  For all algorithms and $\lambda$s, MSE has
   increased relative to the left subplot, which means that in every
   case our estimate of the true underlying solution has worsened.
   This occurs because across all algorithms and $\lambda$s in
   Fig.~\ref{fig:e2_zoom}, the true $K$-term sparsity profile is
   incorrectly estimated and thus the retuning step makes the
   estimated solution worse.  The lesson here is that if one is
   interested in minimizing MSE in low-to-moderate SNR regimes it may
   be best to simply keep the original estimate of the solution rather
   than detect the sparsity profile and retune the result.  If one is
   not certain that the sparsity profile estimate is accurate,
   retuning is likely to increase MSE by fitting the estimated
   solution weights to an incorrect set of generating matrices.  On
   the other hand, if one is confident that the entire sparsity
   profile will be correctly identified with sufficiently high
   probability, retuning will be beneficial; see \cite{Ela2008,
   Goy2008, Fle2006} for related ideas.

\subsection{MRI RF Excitation Pulse Design}
\label{subsec:e3}

   \subsubsection{Overview} For the final experiment we study how well
   the six algorithms design MRI RF excitation pulses.  In the
   interest of space and because the conversion of the physical
   problem into an MSSO format involves MRI physics and requires
   significant background, we only briefly outline how the system
   matrices arise and why simultaneously sparse solutions are
   necessary.  A complete formulation of the problem for engineers and
   mathematicians is given in \cite{Zel2008_CISS}; MRI pulse designers
   may refer to \cite{Zel2008_TMI}.

   \subsubsection{Formulation} For the purposes of this paper, the
   design of an MRI RF excitation pulse reduces to the following
   problem: assume we are given $M$ points in the 2-D spatial domain,
   $\la{r}_1, \ldots, \la{r}_M$, along with $N$ points in a 2-D
   ``Fourier-like'' domain, $\la{k}_1, \ldots, \la{k}_N$.  Each
   $\la{r}_m$ equals $[x_m, y_m]^{\sst{T}}$, a point in space, while
   each $\la{k}_n$ equals $[k_{x,n}, k_{y,n}]^{\sst{T}}$, a point in
   the Fourier-like domain, referred to as ``$k$-space''.  The
   $\la{r}_m$s and $\la{k}_n$s are in units of centimeters (cm) and
   inverse centimeters ($\mbox{cm}^{-1}$), respectively.  The
   $\la{k}_n$s are Nyquist-spaced relative to the sampling of the
   $\la{r}_m$s and may be visualized as a 2-D grid located at low
   $k_x$ and $k_y$ frequencies (where ``$k_x$'' denotes the frequency
   domain axis that corresponds to the spatial $x$ axis).  Under
   reasonable assumptions, energy placed at one or more points in
   $k$-space produces a pattern in the spatial domain; this pattern is
   related to the $k$-space energy via a ``Fourier-like'' transform
   \cite{Pau1989}.  Assume we place an arbitrary complex weight $g_n
   \in \field{C}$ (i.e., both a magnitude and phase) at each of the
   $N$ locations in $k$-space.  Let us represent these weights using a
   vector $\la{g} = [g_1, \ldots, g_N]^{\sst{T}} \in \field{C}^{N}$.
   In an ideal (i.e., physically-unrealizable) setting, the following
   holds:
   \begin{equation}\label{mri1} 
        \la{m} = \la{A} \la{g},
   \end{equation}
   where $\la{A} \in \field{C}^{M \times N}$ is a known dense Fourier
   matrix\footnote{Formally, $\la{A}(m,n) = j \gamma e^{i
   \mbox{\scriptsize{\la{r}}}_m \cdot \mbox{\scriptsize{\la{k}}}_n}$,
   where $j = \sqrt{-1}$ and $\gamma$ is a known lumped gain
   constant.} and the $m$th element of $\la{m} \in \field{C}^{M}$ is
   the image that arises at $\la{r}_m$, denoted $m(\la{r}_m)$, due to
   the energy deposition along the $N$ points on the $k$-space grid as
   described by the weights in the $\la{g}$ vector.
   
   The goal now is to form a desired (possibly complex-valued)
   spatial-domain image $d(\la{r})$ at the given set of spatial domain
   coordinates (the $\la{r}_m$s) by placing energy at some of the
   given $\la{k}_n$ locations while obeying a special constraint on
   how the energy is deposited.  To produce the spatial-domain image,
   we will use a ``$P$-channel MRI parallel excitation system''
   \cite{Kat2003, Set2006}---each of the system's $P$ channels is able
   to deposit energy of varying magnitudes and phases at the $k$-space
   locations and is able to influence the resulting spatial-domain
   pattern $m(\la{r})$ to some extent. Each channel $p$ has a known
   ``profile'' across space, $S_p(\la{r}) \in \field{C}$, that
   describes how the channel is able to influence the magnitude and
   phase of the resulting image at different spatial locations.  For
   example, if $S_3(\la{r}_5) = 0$, then the 3rd channel is unable to
   influence the image that arises at location $\la{r}_5$, regardless
   of how much energy it deposits along $\la{k}_1, \ldots, \la{k}_N$.
   The special constraint mentioned above is as follows: {\em{the
   system's channels may only visit a small number of points in
   $k$-space---they may only deposit energy at $K \ll N$ locations.}}

   We now finalize the formulation of problem: first, we construct $P$
   diagonal matrices $\la{S}_p \in \field{C}^{M \times M}$ such that
   $\la{S}_p(m,m) = S_p(\la{r}_m), m = 1, \ldots, M$.  Now we assume that each channel
   deposits arbitrary energies at each of the $N$ points in $k$-space
   and describe the weighting of the $k$-space grid by the $p$th channel with the vector
   $\la{g}_p \in \field{C}^{N}$.  Based on the physics of the
   $P$-channel parallel excitation system, the overall image
   $m(\la{r})$ that forms is the {\em{superposition}} of the
   profile-scaled subimages produced by each channel:
   \begin{equation}\label{mri2}
   \begin{split}
      \la{m} & = \la{S}_1 \la{A} \la{g}_1 + \cdots + \la{S}_P \la{A} \la{g}_P \\
             & = \la{F}_1 \la{g}_1 + \cdots + \la{F}_P \la{g}_P \\
             & = \la{F}_{\sst{tot}} \la{g}_{\sst{tot}},
   \end{split}
   \end{equation}
   where $\la{m} = [m(\la{r}_1), \ldots, m(\la{r}_M)]^{\sst{T}}$.
   Essentially, (\ref{mri2}) is the real-world version of (\ref{mri1})
   for $P$-channel systems with profiles $S_p(\la{r})$ that are
   not constant across \la{r}.

   Recalling that our overall goal is to deposit energy in $k$-space
   to produce the image $d(\la{r})$, and given the special constraint that we
   may only deposit energy among a small subset of the $N$ points in
   $k$-space, we arrive at the following problem:
   \begin{equation}\label{mri3}
         \mathop{\mbox{min}}_{\sst{$\la{g}_1, \ldots, \la{g}_P$}} 
            \mbox{ } \Vert \la{d} - \la{m} \Vert_2 \,\,
            \mbox{ s.t. the simultaneous $K$-sparsity of the $\la{g}_p$s},
   \end{equation}
   where $\la{d} \in \field{C}^{M} = [d(\la{r}_1), \ldots,
   d(\la{r}_M)]^{\sst{T}} \in \field{C}^{M}$ and $\la{m}$ is given by
   (\ref{mri2}).  That is, we seek out $K < N$ locations in $k$-space
   at which to deposit energy to produce an image $m(\la{r})$ that is
   close in the $\ell_2$ sense to the desired image $d(\la{r})$.
   Strictly and simultaneously $K$-sparse $\la{g}_p$s are the
   only valid solutions to the problem.

   One sees that (\ref{mri3}) is precisely the MSSO system given in
   (\ref{msso_nphard}) and thus the algorithms posed in
   Sec.~\ref{sec:algorithms} are applicable to the pulse design
   problem.  In order to apply the convex minimization techniques
   (IRLS, SOCP, RBRS, and CBCS) to this problem, the only additional
   step needed is to retune any given solution estimate
   $\widehat{\la{g}}_{\sst{tot}}(\mbox{alg},\lambda)$ into a strictly
   and simultaneously $K$-sparse set of vectors; this retuning step is
   computationally tractable and was described in
   Sec.~\ref{subsubsec:e2_res}'s ``{\em{Closer Look}}'' subsection.

   {\em{Aside}.} An alternative approach to decide where to place
   energy at $K$ locations in $k$-space is to compute the Fourier
   transform of $d(\la{r})$ and decide to place energy at $(k_x, k_y)$
   frequencies where the Fourier coefficients are largest in magnitude
   \cite{Yip2006}.  This method does yield valid $K$-sparse energy
   placement patterns, but empirically it is always outperformed by
   the convex minimization approaches \cite{Zel2007, Zel2008_CISS,
   Zel2008_TMI} so we do not delve into the Fourier-based method in
   this paper.


   \subsubsection{Experimental Setup} 
\begin{figure}
   \begin{center}
   \small
      \begin{tabular}{c}
         \epsfig{figure=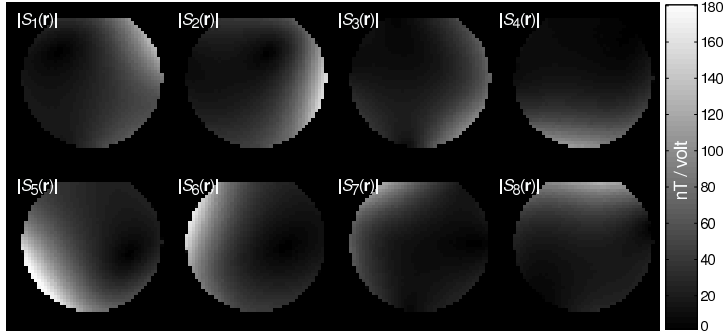,width=3.75in}
      \end{tabular}

      \caption{{\bf{Profile Magnitudes of the Eight-Channel Parallel
      Excitation MRI System}.} Here the magnitudes of the
      $S_p(\la{r})$s are depicted for $p = 1, \ldots, 8$; 356 samples
      of each $S_p(\la{r})$ are taken within the nonzero region of
      influence and stacked into the diagonal matrix $\la{S}_p$ used
      in (\ref{mri2}). Across space, the $S_p(\la{r})$s are not
      orthogonal---their regions of influence overlap each other to
      some extent.}

      \label{fig:e3_profiles} 
   \end{center}
   \end{figure}

\begin{figure}
   \begin{center}
   \small
     \begin{tabular}{c}
       \epsfig{figure=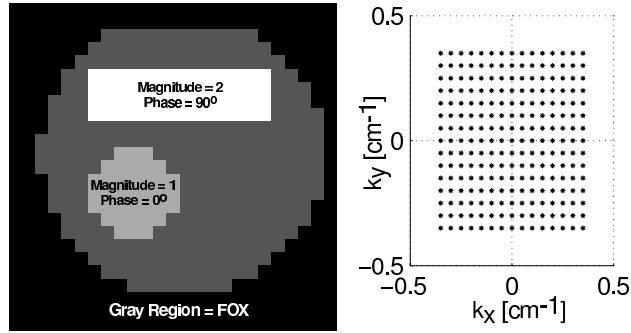,width=3.25in}
     \end{tabular}

     \caption{{\bf{Desired Image and $k$-Space Grid}.} Left image:
     desired complex-valued image, $d(\la{r})$.  Medium-gray region =
     FOX; other regions denote locations where we want image to be
     nonzero with the given magnitudes and phases.  Sampling
     $d(\la{r})$ at the 356 locations within the FOX allows us to
     construct \la{d} in (\ref{mri2}). Right subplot: $15 \times 15$
     grid of $N=225$ candidate $k$-space locations, $\la{k}_1, \ldots,
     \la{k}_{225}$, at which the $P$ channels may deposit energy and
     thus influence the resulting image.  The physical constraints of
     the MRI excitation process force us to place energy at only a
     small number of grid locations.}

     \label{fig:e3_tgt} 
   \end{center}
   \end{figure}

Using an eight-channel system
   (i.e., $P=8$) whose profile magnitudes (the $S_p(\la{r})$s) are
   depicted in Fig.~\ref{fig:e3_profiles}, we will design pulses to
   produce the desired complex-valued image shown in the left subplot
   of Fig.~\ref{fig:e3_tgt}.  We sample the spatial $(x,y)$ domain at
   $M=356$ locations within the region where at least one of the
   profiles in Fig.~\ref{fig:e3_profiles} is active---this region of
   interest is referred to as the {\em{field of excitation}} (FOX) in MRI
   literature.\footnote{Sampling points outside of the FOX where no
   profile has influence is unnecessary because an image can never
   be formed at these points no matter how much energy any given
   channel places throughout $k$-space.} The spatial samples are
   spaced by 0.8 cm along each axis and the FOX has a diameter of
   roughly 20 cm.  Given our choice of $\la{r}_1, \ldots,
   \la{r}_{356}$, we sample the $S(\la{r})$s and $d(\la{r})$ and
   construct the $\la{S}_p$s and $\la{d}$.  Next, we define a grid of
   $N=225$ points in $(k_x, k_y)$-space that is $15 \times 15$ in size
   and extends outward from the $k$-space origin.  The points are
   spaced by $\frac{1}{20} \mbox{ cm}^{-1}$ along each $k$-space axis
   and the overall grid is shown in the right subplot of
   Fig.~\ref{fig:e3_tgt}.  Finally, because we know the 356
   $\la{r}_m$s and 225 $\la{k}_n$s, we construct the $356 \times 225$
   matrix $\la{A}$ in (\ref{mri1}, \ref{mri2}) along with the
   $\la{F}_p$s in (\ref{mri2}).  We now have all the data we need to
   apply the algorithms and determine simultaneously $K$-sparse
   $\la{g}_p$s that (approximately) solve (\ref{mri3}).

   We apply the algorithms and evaluate designs where the use of $K
   \in \{1, \ldots, 30\}$ candidate points in $k$-space is permitted
   (in practical MRI scenarios, $K$ up to 30 is permissible).
   Typically, the smallest $K$ possible that produces a version of
   $d(\la{r})$ to within some $\ell_2$-fidelity is the design that the
   MRI pulse designer will use on a real system.

   To obtain simultaneously $K$-sparse solutions with MP and LSMP, we
   set $K=30$, run each algorithm once, remember the ordered list of
   chosen indices, and back out every solution for $K=1$ through
   $K=30$ via the retuning technique of Sec.~\ref{subsubsec:e2_res}.

   For each convex minimization method (IRLS, SOCP, RBRS, and CBCS),
   we apply the following procedure: first, we run the algorithm for
   14 values of $\lambda \in \left [0, \frac{1}{4} \right ]$, storing
   each resulting solution,
   $\widehat{\la{g}}_{\sst{tot}}(\mbox{alg},\lambda)$.  Then for fixed
   $K$, to determine a simultaneously $K$-sparse deposition of energy
   on the $k$-space grid, we apply the retuning process of
   Sec.~\ref{subsubsec:e2_res} to each of the 14 solutions, obtaining
   14 strictly simultaneously $K$-sparse retuned sets of solution
   vectors, $\widehat{\la{g}}^{(K)}_{\sst{tot}}(\mbox{alg},\lambda)$.
   The one solution among the 14 that best minimizes the $\ell_2$
   error between the desired and resulting images, $\Vert \la{d} -
   \la{F}_{\sst{tot}} \la{g}^{(K)}_{\sst{tot}}(\mbox{alg},\lambda)
   \Vert_2$, is chosen as the solution for the $K$ under
   consideration.  Essentially, we again assume we know a good value
   for $\lambda$ when applying each of the convex minimization
   methods.  To attempt to avoid convergence problems, RBRS and CBCS
   are permitted 5,000 and 10,000 maximum outer iterations,
   respectively (see below).

   \subsubsection{Results} 

\begin{figure}
   \begin{center}
   \small
       \begin{tabular}{cc}
         \epsfig{figure=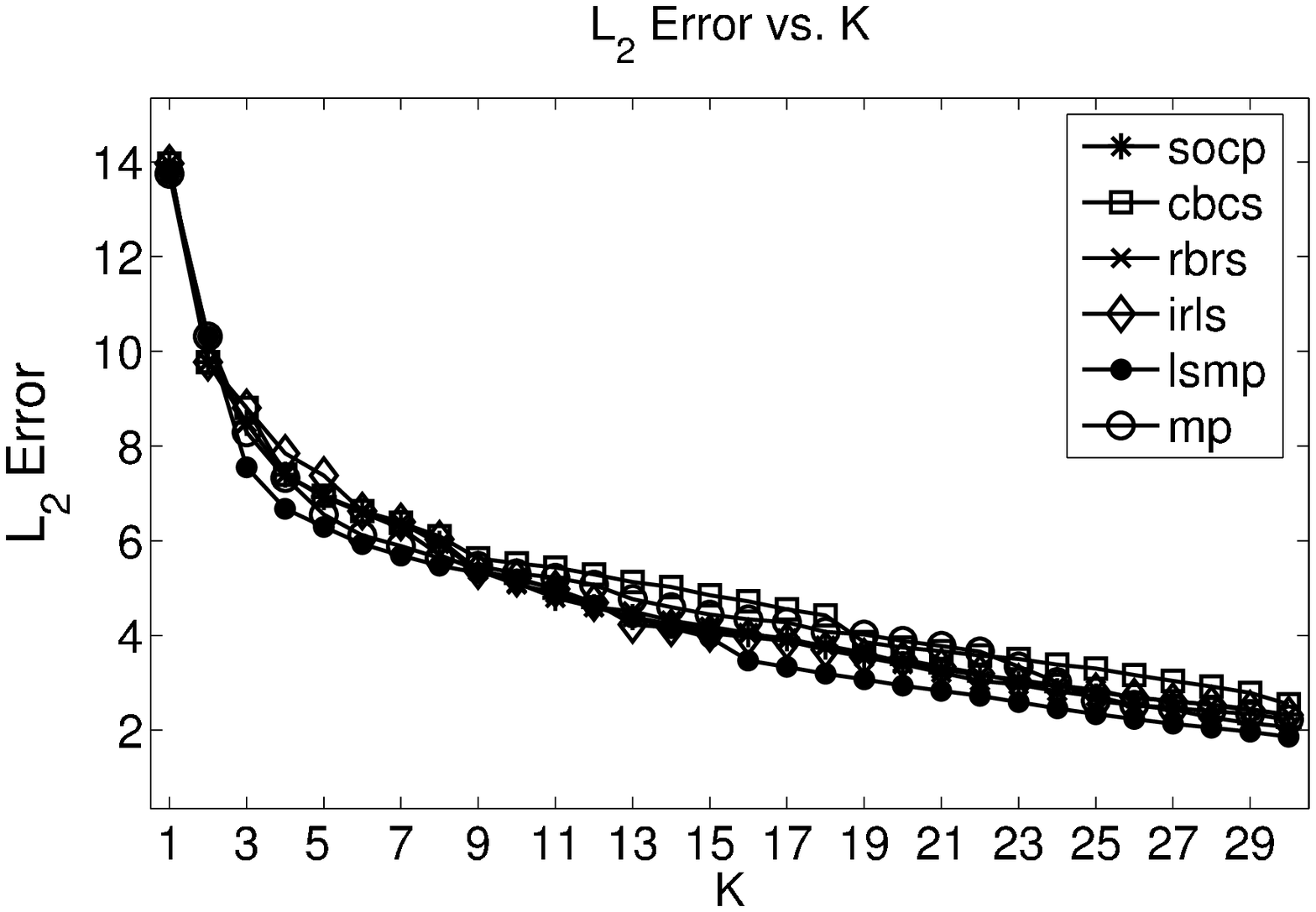,width=2.5in} &
         \epsfig{figure=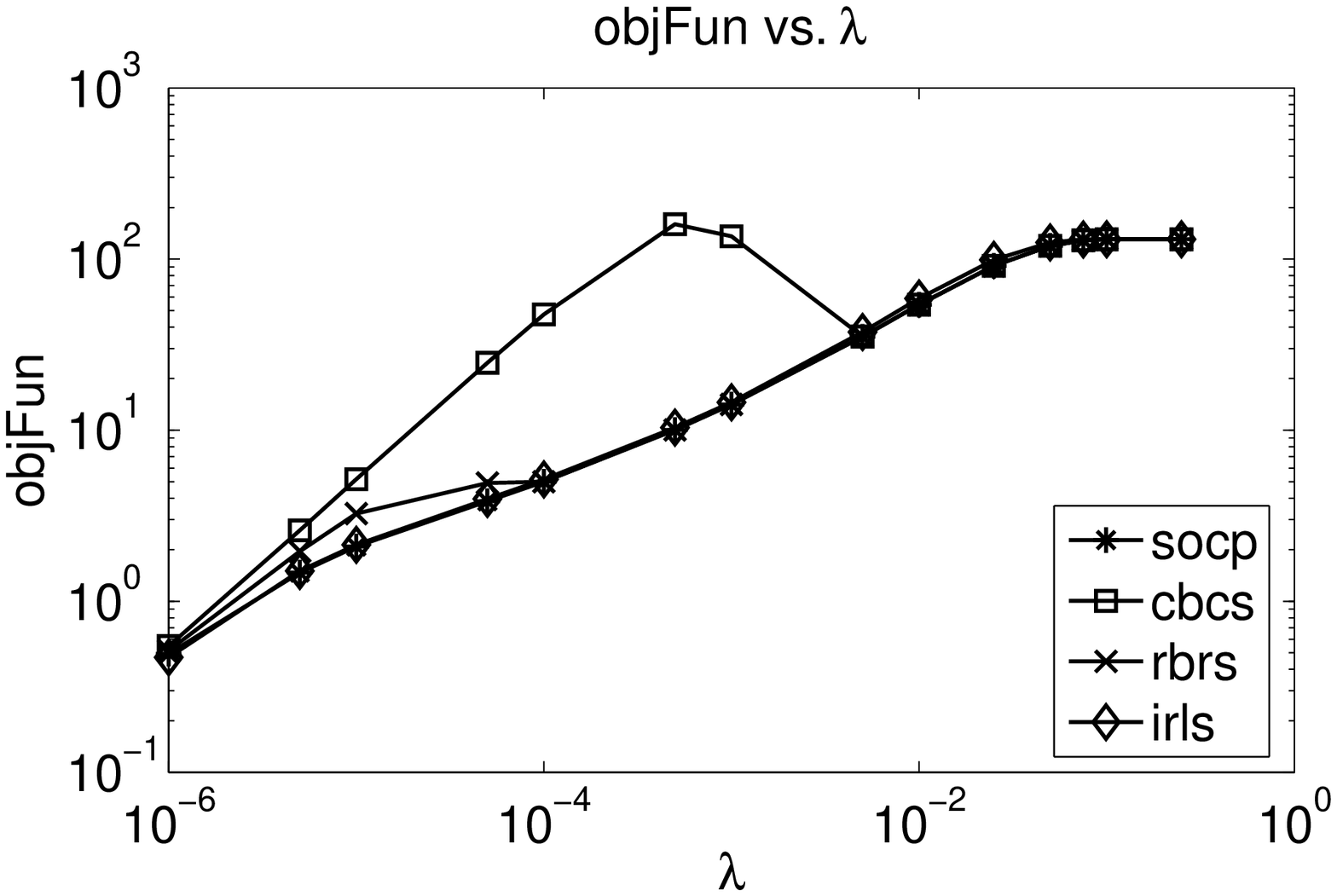,width=2.5in}
       \end{tabular}

       \caption{{\bf{MRI Pulse Design Results}.} Left subplot:
       $\ell_2$ error vs.~$K$ is given for MP, LSMP, IRLS, RBRS,
       CBCS, and SOCP\@.  For fixed $K$, LSMP consistently outperforms
       the other algorithms.  Right subplot: objective function values
       vs. $\lambda$ when SOCP, CBCS, RBRS, and IRLS attempt to
       minimize (\ref{mmv2_eq2}, \ref{mmv2_eq4}).  SOCP and IRLS
       converge and seem to minimize the objective; RBRS does so as
       well for most $\lambda$s.  CBCS routinely fails to converge
       even after 10,000 iterations and thus its solutions yield
       extremely large objective function values.}

       \label{fig:e3_L2} 
   \end{center}
   \end{figure}

   {\em{Image $\ell_2$ Error vs.~Number of Energy Depositions in $k$-Space}.}
   Figure~\ref{fig:e3_L2}'s left subplot shows the $\ell_2$ error
   versus $K$ curves for each algorithm.  As $K$ is increased, each
   method produces images with lower $\ell_2$ error, which is
   sensible: depositing energy at more locations in $k$-space gives
   each technique more degrees of freedom with which to form the
   image.  In contrast to the sparsity profile estimation experiments
   in Sec.~\ref{subsec:e1} and Sec.~\ref{subsec:e2}, however, here
   LSMP is the best algorithm: for each fixed $K$ considered in
   Fig.~\ref{fig:e3_L2}, the LSMP technique yields a simultaneously
   $K$-sparse energy deposition that produces a higher-fidelity image
   than all other techniques.  For example, when $K=17$ LSMP yields a
   solution that leads to an image with $\ell_2$ error of 3.3.  In
   order to produce an image with equal or better fidelity, IRLS,
   RBRS, and SOCP need to deposit energy at $K=21$ points in
   $k$-space, and thus produce less useful designs from an MRI
   perspective.  CBCS fares the worst, needing to deposit energy at
   $K=25$ grid points in order to compete with the fidelity of LSMP's
   $K=17$ image.

   {\em{Closer Look: Objective Function vs.~$\lambda$}.} The right
   subplot of Fig~\ref{fig:e3_L2} shows how well the four convex
   minimization algorithms minimize the objective function
   (\ref{mmv2_eq2}, \ref{mmv2_eq4}) before retuning any solutions and
   enforcing strict simultaneous $K$-sparsity. For each fixed
   $\lambda$, SOCP and IRLS find solutions that yield the same
   objective function value.  RBRS's solutions generally yield
   objective function values equal to those of SOCP and IRLS, but at
   times lead to higher values: in these cases RBRS almost converges
   but does not do so completely.  Finally, for most $\lambda$s CBCS's
   solutions yield extremely large objective function values; in these
   cases CBCS completely fails to converge.

\begin{figure}
   \begin{center}
   \small
       \begin{tabular}{ccc}
         \epsfig{figure=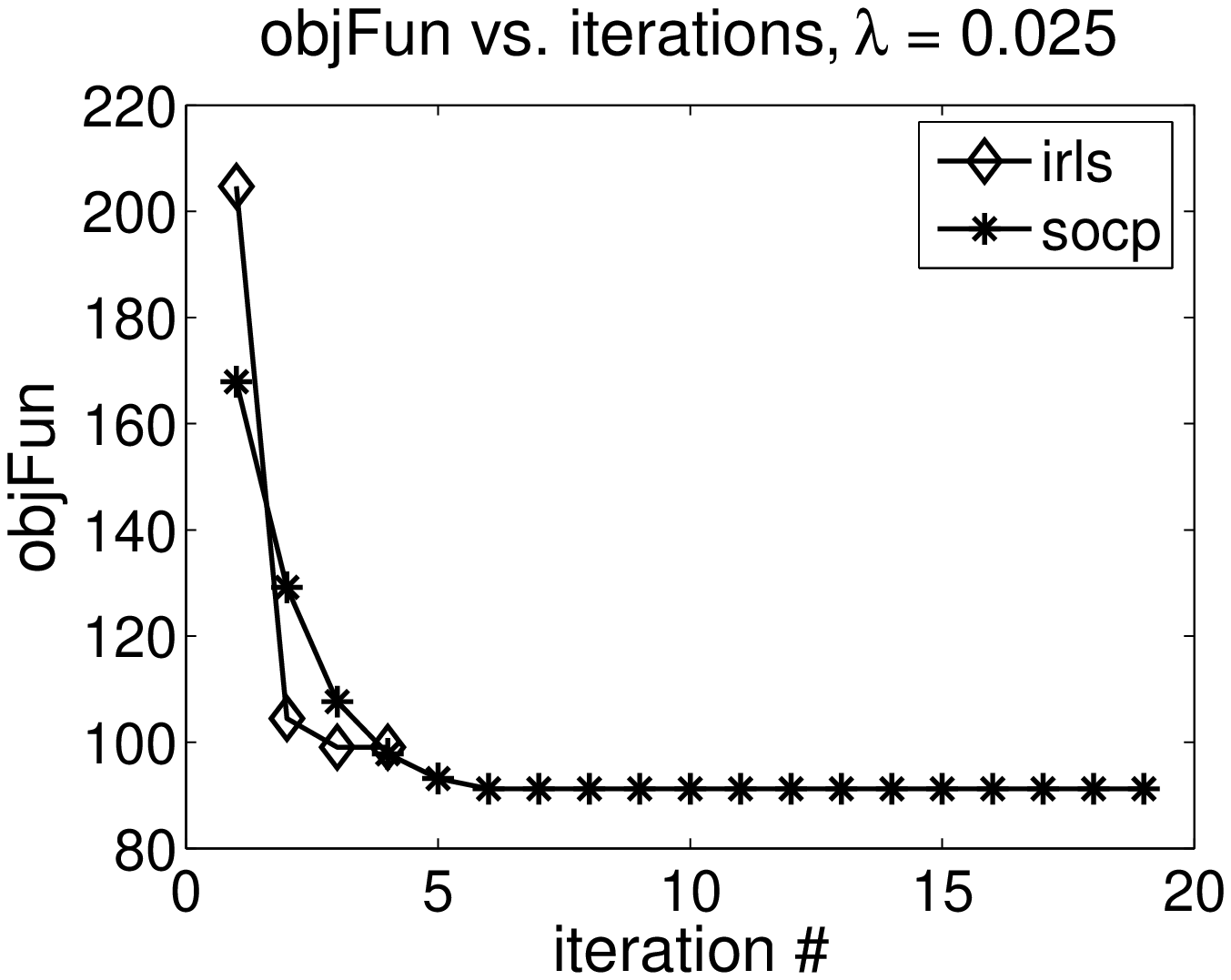,width=\widthC} &
         \epsfig{figure=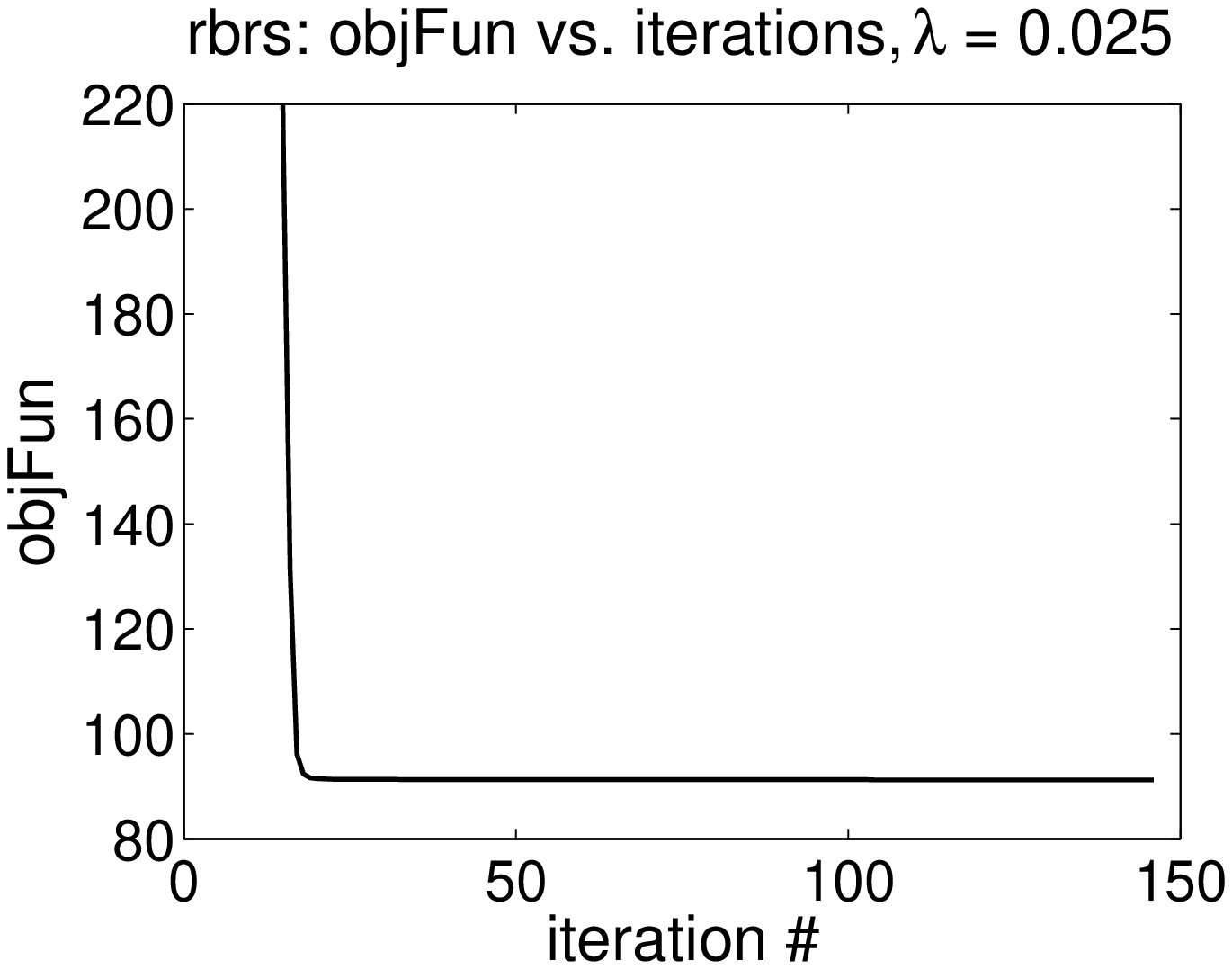,width=\widthC} &
         \epsfig{figure=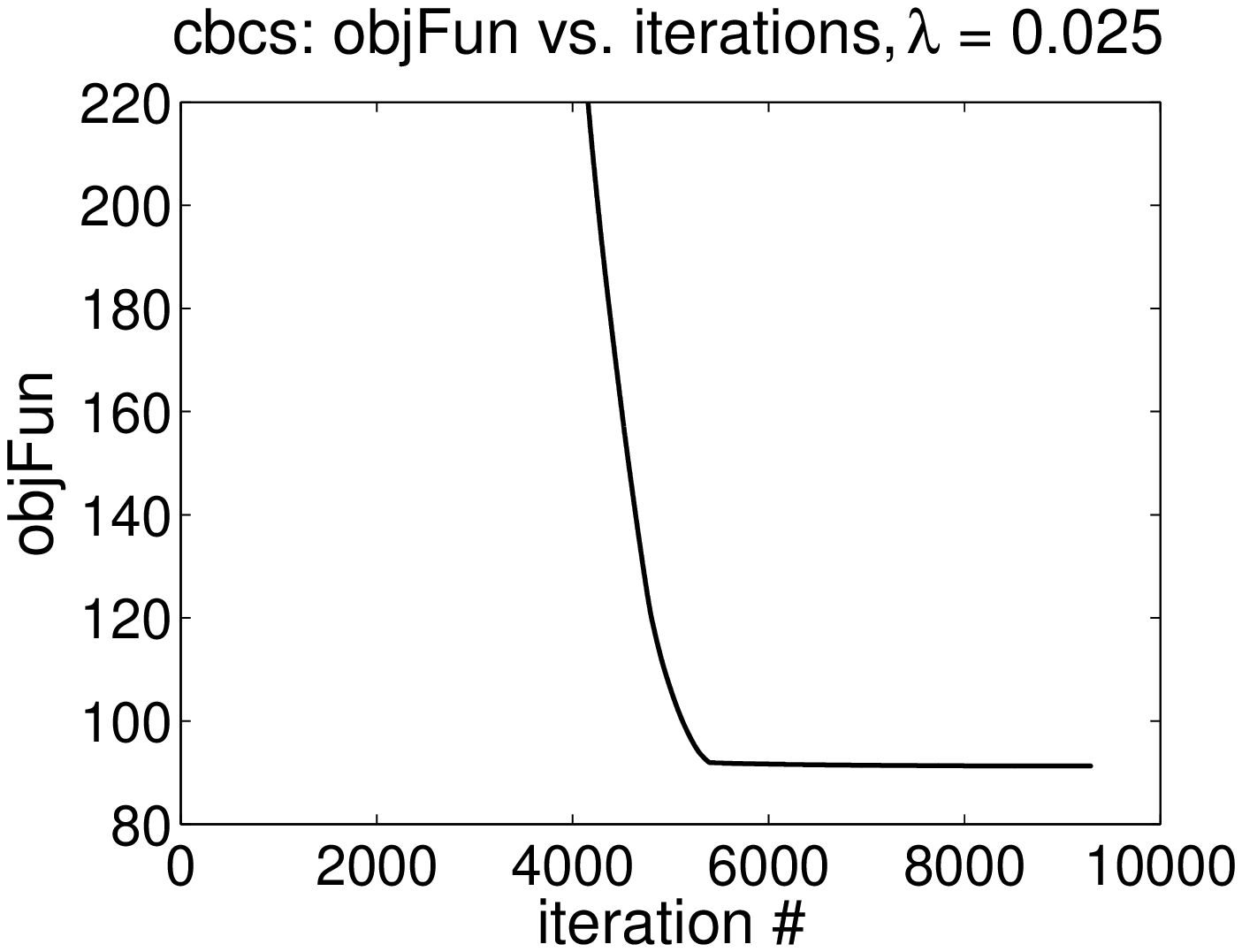,width=\widthC} 
         \\ &
         \epsfig{figure=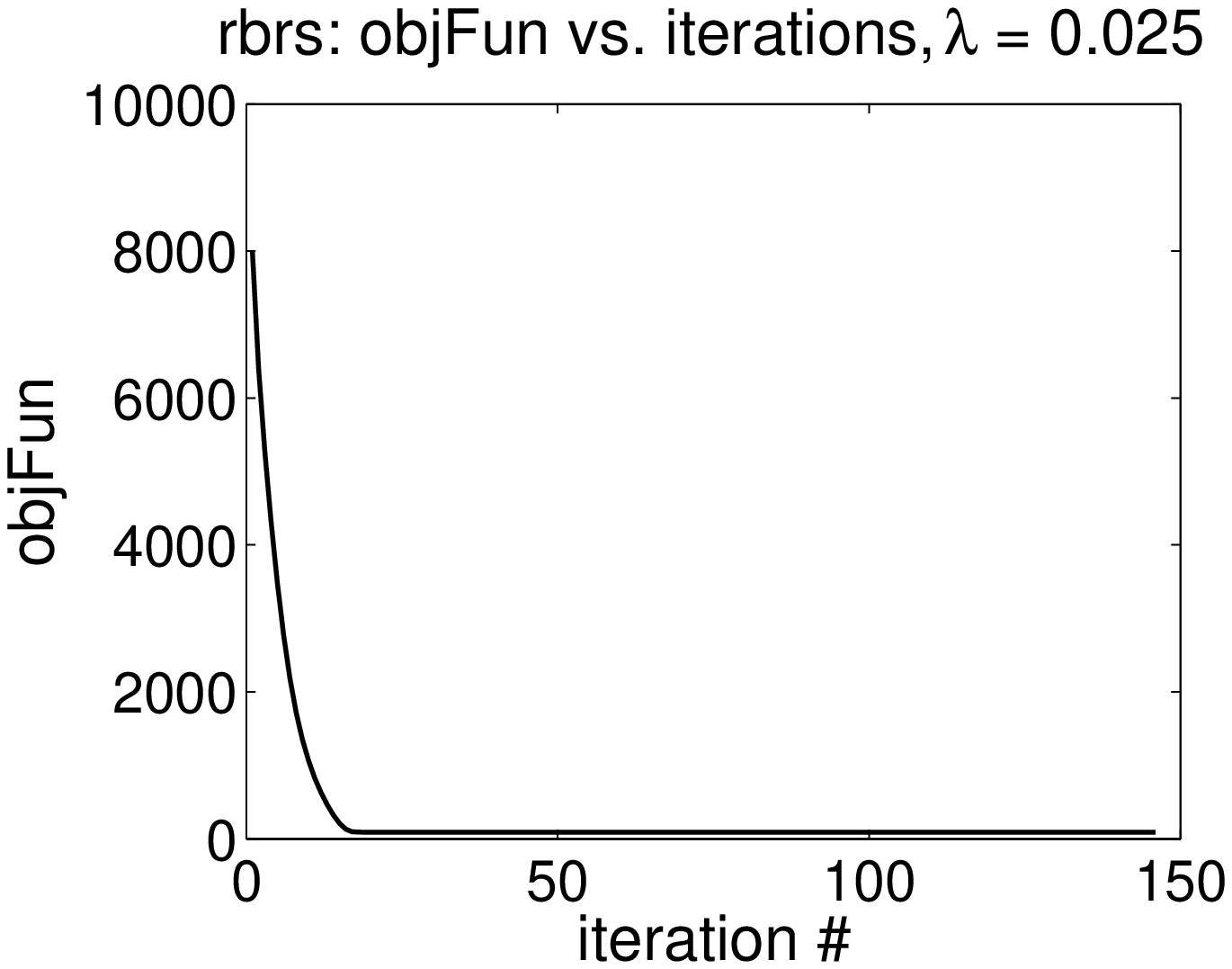,width=\widthC} &
         \epsfig{figure=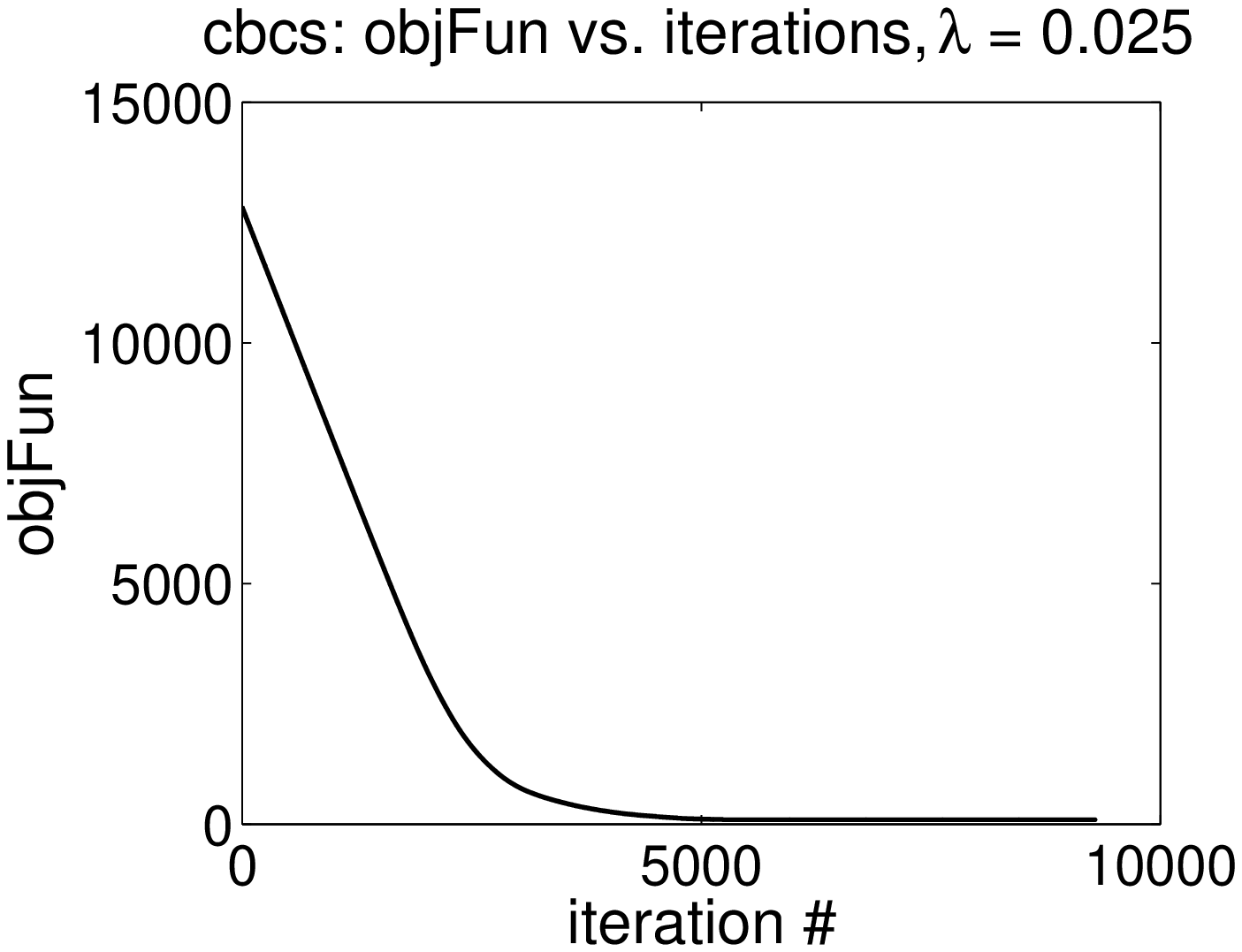,width=\widthC}
       \end{tabular}

       \caption{{\bf{Convergence Behavior of IRLS, SOCP, RBRS, CBCS}.}
       Each subplot illustrates the value of each algorithm's
       objective function (\ref{mmv2_eq2}, \ref{mmv2_eq4}) as the
       algorithm iterates.  Upper row subplots are scaled to have the
       same $y$ axis, whereas the bottom row subplots are ``zoomed
       out'' to illustrate the overall behavior of RBRS and CBCS\@.
       IRLS and SOCP converge rapidly, within 4 and 19 iterations,
       respectively.  RBRS and CBCS require roughly 150 and 10,000
       iterations, respectively.  The runtimes of IRLS, SOCP, RBRS,
       and CBCS in this case are 29, 121, 450, and 5,542 seconds.}

       \label{fig:e3_objFun_vs_iters} 
   \end{center}
   \end{figure}

   {\em{Closer Look: Objective Function Convergence for $\lambda =
   0.025$}.}  The right subplot of Fig~\ref{fig:e3_L2} shows that for
   $\lambda=0.025$, IRLS, SOCP, RBRS, and CBCS generate solutions that
   yield the same objective function value, suggesting that each
   method succeeds at minimizing the objective function.
   Figure~\ref{fig:e3_objFun_vs_iters} illustrates how the algorithms
   converge in this specific case: each subplot tracks the value of an
   algorithm's objective function as it iterates.  Subplots along the
   top row all have the same $y$ axis, giving a close look at how the
   algorithms behave.  The two subplots along the bottom row ``zoom
   out'' along the $y$ axis to show RBRS's and CBCS's total behavior.
   IRLS and SOCP converge rapidly, within 4 and 19 iterations,
   respectively.  RBRS requires roughly 150 outer iterations, while
   CBCS requires nearly 10,000.

   {\em{Runtimes and Peak Memory Usage}.} Setting $K=30$, we run MP
   and LSMP and record the runtime of each.  Across the 14 $\lambda$s,
   IRLS, RBRS, CBCS, and SOCP's runtimes are recorded and averaged.
   The peak memory usage of each algorithm is also noted; these
   statistics are presented in Table~\ref{tab:r3}.  In distinct
   contrast to the smaller-variable-size experiments in
   Sec.~\ref{subsec:e1} and Sec.~\ref{subsec:e2} where all four convex
   minimization methods have relatively short runtimes (under one
   second), here RBRS and CBCS are much slower, leaving IRLS and SOCP
   as the only feasible techniques among the four.  Furthermore, while
   LSMP does indeed outperform IRLS and SOCP on an $\ell_2$ error
   basis (as shown in Fig.~\ref{fig:e3_L2}), the runtime statistics
   here show that LSMP requires order-of-magnitude greater runtime to
   solve the problem---therefore, in some real-life scenarios where
   designing pulses in less than 5 minutes is a necessity, IRLS and
   SOCP are superior.  Finally, in contrast to Sec.~\ref{subsec:e1}'s
   runtimes given in Table~\ref{tab:r1}, here IRLS is 1.9 times faster
   than SOCP and uses 1.4 times less peak memory, making it the
   superior technique for MRI pulse design since IRLS's $\ell_2$ error
   performance and ability to minimize the objective function
   (\ref{mmv2_eq2}, \ref{mmv2_eq4}) essentially equal that of SOCP\@.
   \begin{table}
   \begin{center}
   \small
   \begin{tabular}{|l|c|c|}
   \hline
    {\bf{Algorithm}}   & {\bf{Runtime}} & {\bf{Peak Memory Usage (MB)}} \\ \hline
    MP               & 11 sec            & 704      \\
    LSMP             & 46 min            & 304      \\
    IRLS             & 50 sec            & 320      \\
    RBRS             & 87 min            & 320      \\
    CBCS             & 3.3 hr            & 320      \\
    SOCP             & 96 sec            & 432      \\ \hline
   \end{tabular}
   \end{center}

   \caption{{\bf{Algorithm Runtimes and Peak Memory Usage for MRI
   Pulse Design.}}  Each algorithm's runtime and peak memory usage is
   listed.  The runtimes of the latter four algorithms are averaged
   over the fourteen $\lambda$s per trial.  MP is again faster than
   the other techniques, but consumes more memory because of its
   precomputation step (see {\em{Algorithm \ref{alg:MP}}}).  IRLS and
   SOCP are quite similar performance-wise and minimize the convex
   objective function equally well (see Fig.~\ref{fig:e3_L2}), but we
   see here that IRLS is approximately 1.9 times faster and uses 1.4
   times less peak memory than SOCP, making the former the superior
   technique among the convex methods.}  \label{tab:r3}

   \end{table}

\begin{figure}
   \begin{center}
   \small
       \begin{tabular}{c}
          \epsfig{figure=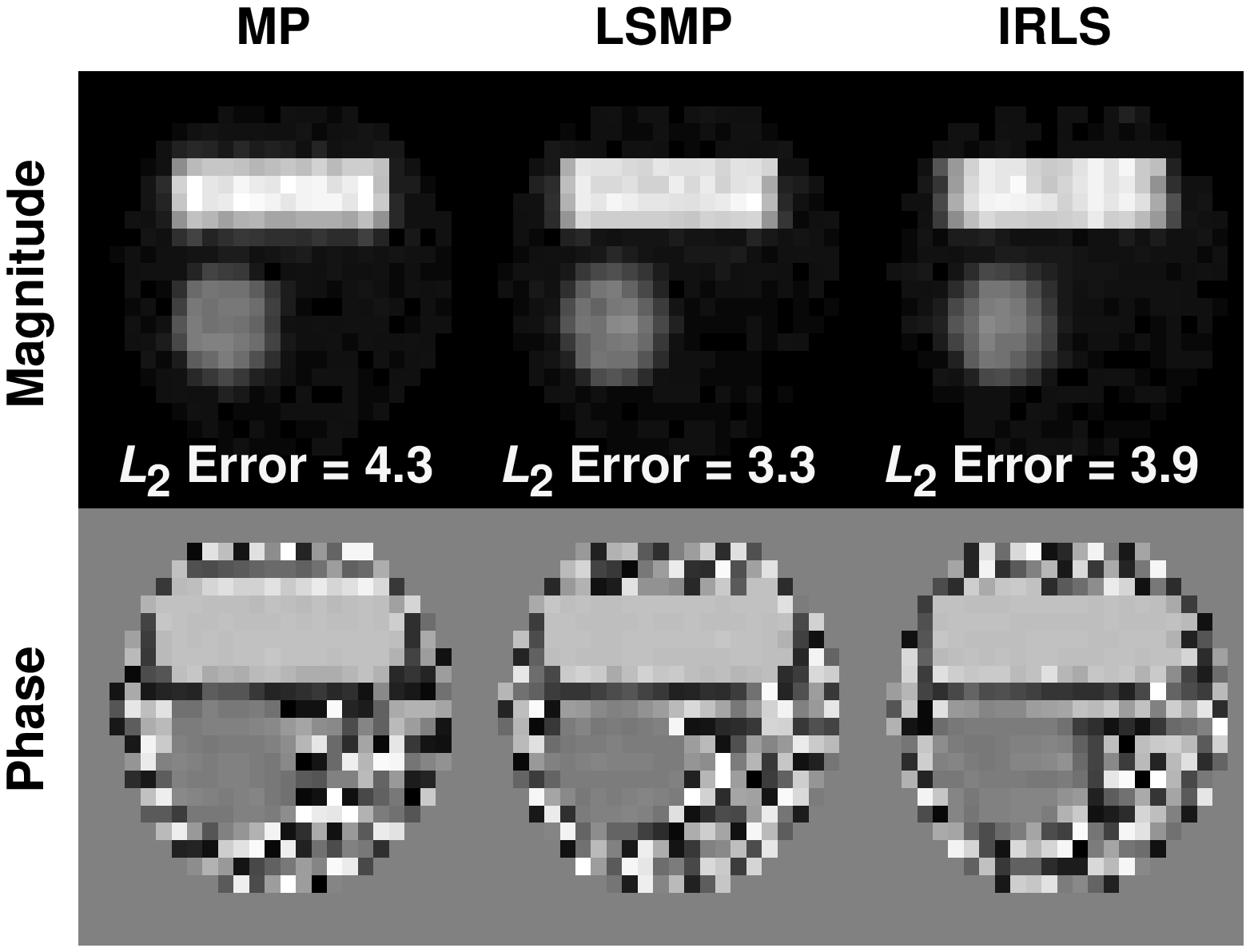,width=4.0in} \\
          \epsfig{figure=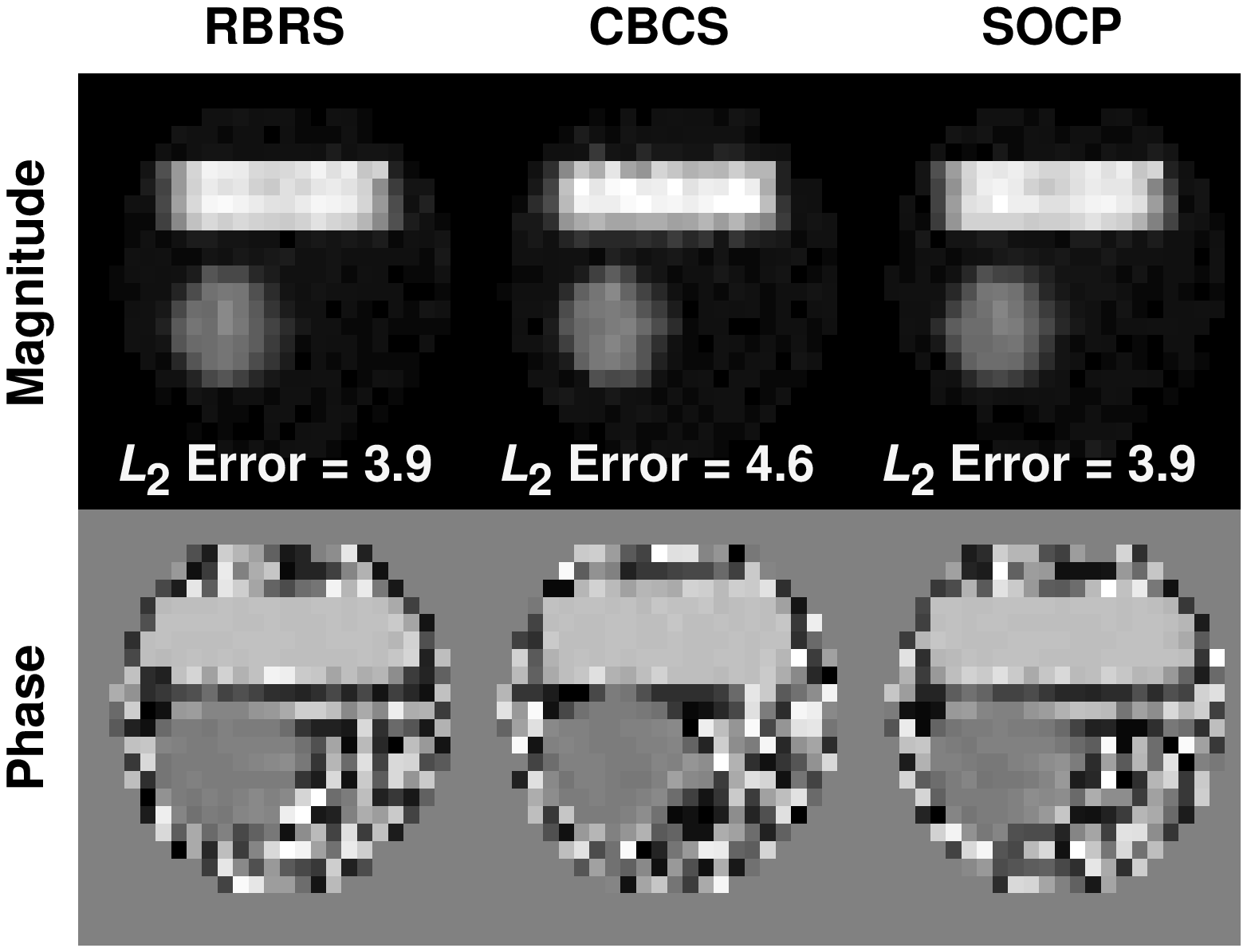,width=4.0in}
       \end{tabular}

       \caption{{\bf{MRI Pulse Design: Images per Algorithm for $K=17$}.}
       Each algorithm is used to solve the MRI pulse design problem
       using 17 energy depositions along the $k$-space grid,
       attempting to produce an image close to the desired one,
       $d(\la{r})$, given in the left subplot of
       Fig.~\ref{fig:e3_tgt}.  From each set of simultaneously
       17-sparse solution vectors, we calculate the resulting image
       via (\ref{mri2}) and display both its magnitude and phase.
       LSMP's image best resembles the desired one; IRLS's, RBRS's,
       and SOCP's images are nearly as accurate; MP's and CBCS's
       images lack crisp edges, coinciding with their larger $\ell_2$
       errors.}

       \label{fig:e3_excits} 

   \end{center}
   \end{figure}

\begin{figure}
   \begin{center}
   \small
       \begin{tabular}{c}
          \epsfig{figure=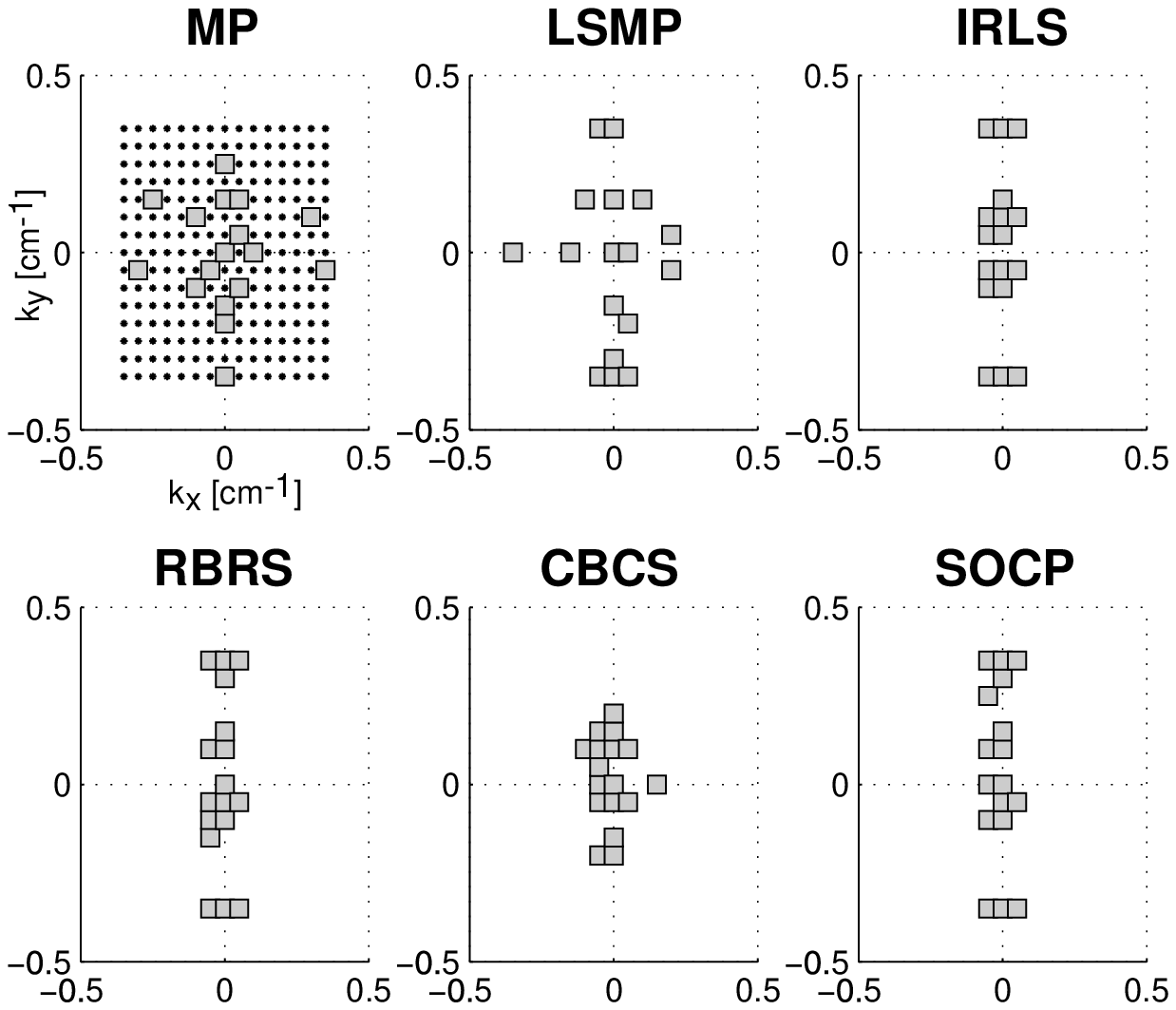,width=4.0in}
       \end{tabular}

       \caption{{\bf{MRI Pulse Design: Energy Deposition Patterns per
       Algorithm for $K=17$}.}  Each algorithm's placement of energy
       in $k$-space is displayed.  LSMP branches out to moderate $k_x$
       frequencies and high $k_y$ frequencies, partly explaining the
       superiority of its image in Fig.~\ref{fig:e3_excits}.  IRLS,
       RBRS, and SOCP succeed in branching out to higher $k_y$
       frequencies but do not place energy at $\vert k_x \vert \gg 0$.
       MP and CBCS fail to spread their energy to high spatial
       frequencies, and thus their images in Fig.~\ref{fig:e3_excits}
       lack distinct edges and appear as ``low-pass filtered''
       versions of $d(\la{r})$.}

       \label{fig:e3_patterns} 

   \end{center}
   \end{figure}

   {\em{Closer Look: Images and Chosen $k$-Space Locations for
   $K=17$}.}  To conclude the experiment, we fix $K = 17$ and observe
   the images produced by the algorithms along with the points at
   which each algorithm chooses to deposit energy along the grid of
   candidate points in $(k_x, k_y)$-space.  Figure~\ref{fig:e3_excits}
   illustrates the images (in both magnitude and phase) that arise due
   to each algorithm's simultaneously 17-sparse set of solution
   vectors,\footnote{Each image is generated by taking the
   corresponding solution $\la{g}_{\sst{tot}}$, computing $\la{m}$ in
   (\ref{mri2}), unstacking the elements of $\la{m}$ into $m(\la{r})$,
   and then displaying the magnitude and phase of $m(\la{r})$.} while
   Fig.~\ref{fig:e3_patterns} depicts the placement pattern chosen by
   each method.  From Fig.~\ref{fig:e3_excits}, we see that each
   algorithm forms a high-fidelity version of the desired image
   $d(\la{r})$ given in the left subplot of Fig.~\ref{fig:e3_tgt}, but
   among the images, LSMP's most accurately represents $d(\la{r})$
   (e.g., consider the sharp edges of the LSMP image's rectangular
   box).  MP's and CBCS's images are noticeably fuzzy relative to the
   others.  The placements in Fig.~\ref{fig:e3_patterns} give insight
   into these performance differences.  Here, LSMP places energy at
   several higher frequencies along the $k_y$ and $k_x$ axes, which
   ensures the resulting rectangle is narrow with sharp edges along
   the spatial $y$ and $x$ axes.  In contrast, CBCS fails to place
   energy at moderate-to-high $(k_x, k_y)$-space frequencies and thus
   cannot produce a rectangle with desirable sharp edges, while MP
   branches out to some extent but fails to utilize high $k_y$
   frequencies.  IRLS, RBRS, and SOCP branch out to higher $k_y$
   frequencies but not to high $k_x$ frequencies, and thus their
   associated rectangles in Fig.~\ref{fig:e3_excits} are sharp along
   the $y$ axis but exhibit less distinct transitions (more fuzziness)
   along the spatial $x$ axis.  In general, each algorithm has
   determined 17 locations at which to place energy that yield a
   fairly good image and each has avoided the computationally
   impossible scenario of searching over all $N$-choose-$K$
   (225-choose-17) possible placements.

\section{Discussion}
\label{sec:discussion}

   \subsection{MRI Pulse Design vs.~Denoising and Source Localization}

   The MRI pulse design problem in Sec.~\ref{subsec:e3} differs
   substantially from the source localization problem in
   Sec.~\ref{subsec:e1}, the denoising experiment in
   Sec.~\ref{subsec:e2}, and other routine applications of sparse
   approximation (e.g. \cite{Don1995, Che1998, Fle2006, Ela2006,
   Cot1999, Cot2005, Mal2005}).  It differs not only in purpose but in
   numerical properties, the latter of which are summarized in
   Table~\ref{tab:diffs}.  While this list will not always hold true
   on an application-by-application basis, it does highlight general
   differences between the two problem classes.\\
   \begin{table}
   \begin{center}
   \small
   \begin{tabular}{|l|l|}
   \hline
    {\bf{MRI Pulse Design}} & {\bf{Denoising and Source Localization}} \\ \hline
    $\bullet$ $\la{F}_p$s overdetermined & $\bullet$ $\la{F}_p$s underdetermined  \\
    $\bullet$  No concept of noise: given \la{d} is $\la{d}_{\sst{true}}$ & 
      $\bullet$  Noisy: given \la{d} is not $\la{d}_{\sst{true}}$ \\
    $\bullet$ Sweep over $\lambda$ useful & 
      $\bullet$ Ideal $\lambda$ unknown \\
    $\bullet$  Metric: $\Vert \la{d} - \la{m} \Vert_2$ & 
      $\bullet$ Metrics: $\Vert \la{g}_{\sst{tot}} - \widehat{\la{g}}_{\sst{tot}} \Vert_2$, and/or \\
        & $\,\,\,$ fraction of rec.~sparsity profile terms \\ \hline
   \end{tabular}
   \end{center}

   \caption{{\bf{Unique Trends of the MRI Pulse Design Problem}.} This
   table highlights differences between the MRI problem and standard
   denoising and source localization applications.  Items here will
   not always be true, instead providing general highlights about
   each problem class.}
   \label{tab:diffs}
   \end{table}

   \subsection{Merits of Row-by-Row and Column-by-Column Shrinkage}
   Even though LSMP, IRLS, and SOCP tend to exhibit superior
   performance across different experiments in this manuscript, RBRS
   and CBCS are worthwhile because unlike the former methods
   that update all $PN$ unknowns concurrently, the shrinkage
   techniques update only a subset of the total variables during each
   iteration.

   For example, RBRS as given in {\em{Algorithm \ref{alg:RBRS}}}
   updates only $P$ unknowns at once, while CBCS as given in
   {\em{Algorithm \ref{alg:CBCS}}} updates but a single scalar at a
   time.  RBRS and CBCS are thus applicable in scenarios where $P$ and
   $N$ are exceedingly large and tuning all $PN$ parameters during
   each iteration is not possible.  If storing and handling $M \times
   PN$ or $PN \times PN$ matrices exceeds a system's available memory
   and causes disk thrashing, RBRS and CBCS, though they require far
   more iterations, might still be better options than LSMP, IRLS, and
   SOCP in terms of runtime.

   \subsection{Future Work}

   \subsubsection{Efficient Automated Control Parameter Selection} A
   fast technique for finding ideal values of $\lambda$ is an open
   problem.  It might be useful to investigate several approaches to
   automated selection such as the ``L-curve'' method \cite{Han1994},
   universal parameter selection \cite{Don1994}, and min-max parameter
   selection \cite{Joh1994}.

   \subsubsection{Runtime, Memory, and Complexity Reduction} As noted
   in Sec.~\ref{sec:algorithms}, LSMP's computation and runtime could
   be improved upon by extending the projection based recursive
   updating schemes of \cite{Cot1999, Cot2005} to MSSO LSMP\@.  In
   regards to the MRI design problem, runtime for any method might be
   reduced via a multi-resolution approach (as in \cite{Mal2005}) by
   first running the algorithm with a coarse $k$-space frequency grid,
   noting which energy deposition locations are revealed, and then
   running the algorithm with a grid that is finely sampled around the
   locations suggested by the coarse result.  This is faster than
   providing the algorithm a large, finely-sampled grid and attempting
   to solve the sparse energy deposition task in one step.  An initial
   investigation has shown that reducing the maximum frequency extent
   of the $k$-space grid (and thus the number of grid points, $N$) may
   also reduce runtime without significantly impacting image fidelity
   \cite{Zel2008_TMI}.

   \subsubsection{Shrinkage Algorithm Convergence Improvements} When
   solving the MRI pulse design problem, both RBRS and CBCS required
   excessive iterations and hence exhibited lengthy runtimes.  The
   latter was especially problematic as illustrated in
   Fig.~\ref{fig:e3_objFun_vs_iters}.  To mitigate these problems, one
   may consider extending parallel coordinate descent (PCD) shrinkage
   techniques used for single-system single-output sparse
   approximation (as in \cite{Ela2006_TransIT, Ela2006}).  Sequential
   subspace optimization (SESOP) \cite{Ela2007} might also be employed
   to reduce runtime.  Combining PCD with SESOP and adding a line
   search after each iteration would yield sophisticated versions of
   RBRS and CBCS\@.

   \subsubsection{Multiple-System Multiple-Output (MSMO) Simultaneous
   Sparse Approximation} In the future it may be useful to consider a
   combined problem where there are multiple observations as well as
   multiple system matrices.  That is, assume we have a series of $J$
   observations, $\la{d}_1, \ldots, \la{d}_J$, each of which arises
   due to a set of $P$ simultaneously $K$-sparse unknown vectors
   $\la{g}_{1,j}, \ldots, \la{g}_{P,j}$\footnote{The $K$-term
   simultaneous sparsity profile of each set of $\la{g}_{p,j}$s may or
   may not change with $j$.}  passing through a set of $P$ system
   matrices $\la{F}_{1,j}, \ldots, \la{F}_{P,j}$ and then undergoing
   linear combination, as follows:
   \begin{equation}\label{mmv3_eq1}
       \la{d}_j = \la{F}_{1,j} \la{g}_{1,j} + \cdots + \la{F}_{P,j} \la{g}_{P,j}
           = \sum_{p=1}^{P} \la{F}_{p,j} \la{g}_{p,j}
              \mbox{ for } j = 1, \ldots, J.
   \end{equation}
   If $\la{F}_{p,j}$ is constant for all $J$ observations then the
   problem reduces to
   \begin{equation}\label{mmv3_eq2}
       \la{d}_j = \la{F}_{1} \la{g}_{1,j} + \cdots + \la{F}_{P} \la{g}_{P,j}
            = \la{F}_{\sst{tot}} \la{g}_{\sst{tot}, j}
              \mbox{ for } j = 1, \ldots, J,
   \end{equation}
   and we may stack the matrices and terms as follows:
   \begin{equation}\label{mmv3_eq3}
       \left[ \la{d}_1, \ldots, \la{d}_J \right] = 
           \la{F}_{\sst{tot}} \left[ \la{g}_{\sst{tot},1}, \ldots, \la{g}_{\sst{tot},J} \right].
   \end{equation}
   Having posed (\ref{mmv3_eq1}, \ref{mmv3_eq2}, \ref{mmv3_eq3}), one
   may formulate optimization problems similar to (\ref{mmv1},
   \ref{mmv2_eq2}) to determine simultaneously sparse $\la{g}_{p,j}$s
   that solve (\ref{mmv3_eq3}).  Algorithms to solve such problems may
   arise by combining the concepts of SSMO algorithms
   \cite{Cot2005,Mal2005,Tro2006_I,Tro2006_II} with those of the MSSO
   algorithms posed earlier.

\section{Conclusion}
\label{sec:conclusion}

  We defined the linear inverse multiple-system, single-output (MSSO)
  simultaneous sparsity problem where simultaneously sparse sets of
  unknown vectors are required as the solution.  This problem differed
  from prior problems involving multiple unknown vectors because in
  this case, each unknown vector was passed through a different system
  matrix and the outputs of the various matrices underwent linear
  combination, yielding only one observation vector.

  To solve the proposed MSSO problem, we formulated three greedy
  techniques, matching pursuit, orthogonal matching pursuit, and least
  squares matching pursuit, along with algorithms based on iteratively
  reweighted least squares, iterative shrinkage, and second-order cone
  programming methodologies. The MSSO algorithms were evaluated across
  noiseless and noisy sparsity profile estimation experiments as well
  as a magnetic resonance imaging pulse design experiment; for
  sparsity profile recovery, algorithms that minimized the relaxed
  convex objective function outperformed the greedy methods, whereas
  in the noiseless magnetic resonance imagine pulse design experiment,
  greedy LSMP exhibited superior performance.

  Finally, when deriving CBCS for complex-valued data, we proved that
  seeking a single sparse complex-valued vector is equivalent to
  seeking two simultaneously sparse real-valued vectors---we mapped
  single-vector sparse approximation of a complex vector to the MSSO
  problem, increasing the applicability of algorithms that solve the
  latter.

  Overall, while improvements upon these seven algorithms (and new
  algorithms altogether) surely do exist, this manuscript has laid the
  groundwork of the MSSO problem and conducted an initial exploration
  of algorithms with which to solve it.

\section*{Acknowledgments} 

   The authors thank D.~M.~Malioutov for assistance with the
   derivation step that permitted the transition from (\ref{socp2}) to
   (\ref{socp3}) in Sec.~\ref{subsec:socp}, as well as K.~Setsompop,
   B.~A.~Gagoski, V.~Alagappan, and L.~L.~Wald for collecting the
   experimental coil profile data in Fig.~\ref{fig:e3_profiles}.

   The authors gratefully acknowledge the
   following sponsors and associated grants: National Institute of Health
   1P41RR14075, 1R01EB000790, 1R01EB006847, and 1R01EB007942; NSF
   CAREER Award CCF-643836; United States Department of Defense
   National Defense Science and Engineering Graduate Fellowship
   F49620-02-C-0041; the MIND Institute; the A.~A.~Martinos Center for
   Biomedical Imaging; and R.~J.~Shillman's Career Development Award.

\bibliographystyle{siam}
\bibliography{paper}

\end{document}